\documentclass[10pt, a4paper]{article}

\usepackage[T1]{fontenc} 
\usepackage[utf8]{inputenc}
\usepackage[english]{babel}
\usepackage{graphicx}
\usepackage{caption}
\usepackage[usenames,dvipsnames]{xcolor}
\usepackage{color}
\usepackage{amsmath,amssymb,amsthm}
\usepackage{float}
\usepackage{subfigure}
\usepackage{verbatim}
\usepackage{algpseudocode,algorithm,algorithmicx}
\usepackage{multirow}
\usepackage{authblk}
\usepackage[top=1in, bottom=1.25in, left=1.25in, right=1.25in]{geometry}
\usepackage{url}

\algdef{SE}[DOWHILE]{Do}{doWhile}{\algorithmicdo}[1]{\algorithmicwhile\ #1}%

\graphicspath{%
	{figures/}%
}

\newcommand{\reviewerB}[1]{#1}
\DeclareMathOperator*{\argmin}{\arg\!\min}

\theoremstyle{remark}
\newtheorem{remark}{Remark}

\begin{document}


\title{A POD-selective inverse distance weighting method for fast parametrized shape morphing}

\author[1]{Francesco Ballarin}
\author[1,2]{Alessandro D'Amario}
\author[3]{Simona Perotto}
\author[1]{Gianluigi Rozza}

\affil[1]{mathLab, Mathematics Area, SISSA, via Bonomea 265, I-34136 Trieste, Italy}
\affil[2]{Politecnico di Milano, via Giuseppe La Masa, 34, I-20156 Milano, Italy}
\affil[3]{MOX - Modeling and Scientific Computing, Dipartimento di Matematica, Politecnico di Milano, piazza Leonardo da Vinci 32, I-20133 Milano, Italy}

\maketitle

\begin{abstract}
Efficient shape morphing techniques play a crucial role in the approximation of partial differential equations defined in parametrized domains, such as for fluid-structure interaction or shape optimization problems. In this paper, we focus on Inverse Distance Weighting (IDW) interpolation techniques, where a reference domain is morphed into a deformed one via the displacement of a set of control points. We aim at reducing the computational burden characterizing a standard IDW approach without significantly compromising the accuracy. To this aim, first we propose an improvement of IDW based on a geometric criterion which automatically selects a subset of the original set of control points. Then, we combine this new approach with a dimensionality reduction technique based on a Proper Orthogonal Decomposition of the set of admissible displacements. This choice further reduces computational costs. We verify the performances of the new IDW techniques on several tests by investigating the trade-off reached in terms of accuracy and efficiency.

\end{abstract}


\section{Introduction}\label{sec:introduction}
Shape morphing plays a meaningful role in several engineering and life science fields, such as, for instance, aero-elasticity \cite{samareh2004}, high performance boat design \cite{lombardiparoliniquarteronirozza2012,lombardi2012}, modelling of the cardiovascular system \cite{Ballarin2015,Shontz1a}. 
On one hand, the same physical problem may involve the deformation of the domain, for instance, when dealing with fluid-structure interaction problems \cite{bazilevs2013computational}. On the other hand, several applications imply  iterative procedures where the same problem is solved in different geometric configurations. 
This is the case, e.g., of a multi-query context, such as shape optimization problems, where the shape of the domain is varied until a prescribed cost functional is minimized (or maximized) \cite{SokolowskiZolesio1992}.

Our interest is focused on partial differential equations solved in a domain which changes in time. In such a case, the goal becomes twofold, since we aim to accurately approximate the domain as well as the partial differential equation.
In general, it is not computationally affordable to generate a new discretization (mesh) of the domain at each deformed configuration.
Indeed, mesh generation may be a time consuming procedure (both in terms of CPU time and assembling)
and, sometimes, it is not directly integrated in the solver at hand. Therefore, it is often more convenient to generate a mesh for a reference configuration and then to morph it into the deformed grid. This work can be framed in such a context.

Several shape morphing techniques are available in the computer graphics community.
A reference mesh is deformed by displacing some points (the so-called control points); then the shape morphing map associated with these points is applied to the whole discretized domain, thus avoiding any remeshing. 
These techniques have been recently applied to parametrized partial differential equations (PPDEs). 
For instance, shape morphing techniques based on \emph{Free Form Deformation} (FFD) \cite{sederbergparry1986} and \emph{Radial Basis Functions} (RBF) interpolation \cite{buhmann2003} have been successfully applied to PPDEs \cite{LassilaRozza2010,ManzoniRBF,TezzeleBallarinRozza2017}, shape optimization problems \cite{BallarinManzoniRozzaSalsa2013,Shontz2a,Marco}, reconstruction of scattered geometrical data \cite{CrivellaroRBF}, mesh motion and interface coupling for fluid-structure interaction problems \cite{beckertwendland2001,Farhat2a,fortirozza2014}, interpolation between non-conforming meshes \cite{deparisfortiquarteroni2013}, sensitivity analysis studies in complex geometrical configurations \cite{BallarinJCP}.

In this paper we focus on a different approach, the \emph{Inverse Distance Weighting} (IDW) interpolation
that relies on the inverse of a weighted sum of distances between mesh nodes, some of which will be chosen as control points \cite{shepard1968,witteveen2010,witteveenbijl2009,forti2012}.

A proper choice of the control points is a crucial step independently of the adopted morphing technique.
A first requirement is to keep the number of control points as small as possible, since the complexity 
of the morphing map evaluation increases as the number of control points becomes larger. 
Such selection is usually driven by a prior knowledge of the phenomenon at hand. 
In some cases an automatic selection can be carried out. 
For instance, a sensitivity analysis to the control points is proposed for shape optimization problems 
in \cite{BallarinManzoniRozzaSalsa2013}, where the control points selected are the ones providing the largest variation of a cost functional. Local minima and maxima of structure eigenmodes are employed in \cite{fortirozza2014} as control points to perform mesh motion in fluid-structure interaction problems.

Here, we propose a variant of the IDW method that combines the shape morphing procedure 
with a new criterion to select the control points
based on geometrical considerations and, consequently, independently of the underlying partial differential equations.
We will refer to this variant as to Selective IDW (SIDW) approach.We remark that the Selective methodology proposed in this work for IDW can be applied to RBF and to improved variants of FFD, such as the one based on NURBS functions\footnote{Control points are constrained to a lattice in the original FFD approach \cite{sederbergparry1986}, while later FFD variants, such as the NURBS-FFD method, allow a greater freedom in the control point selection.} \cite{LamousinWaggenspack}.
We further reduce the
computational cost of morphing by applying dimensionality reduction to the deformations associated with shape morphing.
In particular, thanks to the employment of a reference domain, we identify any shape deformation 
with a certain parametric configuration of the original structure, the parameter being strictly related to 
the control point displacement. This 
suggests us to resort to a dimensionality reduction approach, and we choose, in particular, the Proper Orthogonal Decomposition (POD).
We will refer to this combined technique as to POD-SIDW. POD, also known as Principal Component Analysis (PCA), is used in the literature, for instance, in 
statistical shape analysis \cite{Drydmard98}, as well as for mesh deformation and optimization \cite{Anttonen2003,Walton2013a}.
Manifold learning techniques are combined with PCA for efficient structural shape optimization in \cite{raghavanbreitkopftourbiervillon2013}. An equivalent shape representation by means of currents is sought in \cite{Vaillant2005} in order to introduce an Hilbert space over shapes and apply PCA. Application to life sciences of 
PCA has been proposed in \cite{GuibertMcLeodCaiazzoMansiFernandezSermesantPennecVignon-ClementelBoudjemlineGerbeau2014,McLeodCaiazzoFernandezMansiVignon-ClementelSermesantPennecBoudjemlineGerbeau2010}.

The paper is organized as follows. In Section~\ref{sec:standard_inverse_distance_weighting}, after formalizing the standard approach, we set the selective variant of IDW interpolation. Both the approaches are customized in a shape morphing setting and numerically compared on structural configurations of interest in aeronautic and naval engineering.
Section~\ref{sec:model_reduction} deals with dimensionality reduction. POD is directly applied to the selective IDW interpolation
and the numerical benefits due to such a merging of techniques are investigated when dealing with shape morphing. 
Some conclusions are finally drawn in Section~\ref{sec:ROMIDW conclusions} together with perspectives for a future development of the present work.

\section{Inverse Distance Weighting techniques}\label{sec:standard_inverse_distance_weighting}
In this section, we consider two techniques to drive a shape morphing process, with particular interest in fluid-structure interaction (FSI) problems. 
In this context, a mesh motion is performed to extend the structural displacement (computed with an appropriate structural model) into the fluid, resulting in a morphed fluid domain. In particular, we first introduce 
the original Inverse Distance Weighting (IDW) interpolation, and then we propose a new
variant of this method to overcome some of its intrinsic computational limits.

\subsection{The Inverse Distance Weighting (IDW) approach}
The IDW method has been originally proposed in \cite{shepard1968} to deal with 
two-dimensional interpolation problems, and successively extended to higher-di\-men\-sional
and applicative settings 
(see, e.g., \cite{CebralLohner,luwong2008,witteveenbijl2009,witteveen2010}).

Let $\Omega_r \subset \mathbb{R}^d$, with $d=1,\ 2, 3$, be the reference 
domain, and let $u: \Omega_r \to \mathbb{R}$ be a continuous function. 
We select a prescribed set ${\mathcal C}=\{\mathbf{c}_k\}_{k=1}^{\mathcal{N}_c}$
of points in $\Omega_r$, known as \emph{control points}. 
Then, the \emph{IDW interpolant} of $u$, $\Pi_{\rm IDW}(u)$, coincides with the continuously differentiable function
\begin{equation}\label{IDW_formula}
\Pi_{\rm IDW}(u)(\mathbf{x}) =
\sum _{k=1}^{\mathcal{N}_c}  w_k(\mathbf{x}) \, u(\mathbf{c}_k)\quad \mathbf{x}\in \Omega_r,
\end{equation}
where the local weight functions $w_k: \Omega_r \to \mathbb{R}$, for $k=1, \ldots, {\mathcal{N}_c}$, are defined by
\begin{equation}\label{weighting_function}
w_k(\mathbf{x}) = \begin{cases}
\frac{\displaystyle\| \mathbf{x} - \mathbf{c}_k \|^{-p}}{\displaystyle\sum_{j=1}^{\mathcal{N}_c}{\| \mathbf{x} - \mathbf{c}_j \|^{-p}} } & \text{ if } \mathbf{x} \neq \mathbf{c}_k,\\[2mm]
1 & \text{ if } \mathbf{x} = \mathbf{c}_k,\\
0 & \text{ otherwise},
\end{cases}
\end{equation}
with $\| \cdot\|$ the standard Euclidean norm, and for some positive integer $p$. Notice that
weights are automatically selected such that the influence of the $k$-th control point 
$\mathbf{c}_k$ on $\mathbf{x}$ diminishes 
as the distance between $\mathbf{x}$ and $\mathbf{c}_k$ increases.
Power $p$ tunes such an inverse dependence, so that the weight assigned to a point $\mathbf{x}$ far from 
the control points becomes smaller and smaller as $p$ increases.
With reference to FSI problems, the choice of parameter $p$ is crucial to avoid compenetration effects between fluid 
and structure \cite{witteveen2010}. 

\subsubsection{IDW for shape morphing}
To set shape morphing into the IDW interpolation framework, we identify the reference domain $\Omega_r$
with the initial configuration of the physical domain. Then, we consider 
a discretization $\Omega_h=\{\mathbf{x}_i\}_{i=1}^{\mathcal{N}_h} \subset \Omega_r$, with $\mathbf{x}_i$ being all the nodes of the mesh, and we identify function $u$ in \eqref{IDW_formula} 
with the deformation $d$ of points $\mathbf{x}_i$ and the set ${\mathcal C}$ of the control points (also denoted as parameters)
with a subset of $\Omega_h$. In such a context, the IDW interpolant, $\Pi_{\rm IDW}(d)$, represents the so-called \emph{shape parametrization map}. \\
The choice of the control points clearly represents a crucial issue.
In general, ${\mathcal C} \subset \Omega_h$ coincides with the set of the
boundary points of $\Omega_h$.
This is a common practice in several engineering applications, in particular, in
FSI problems, where the displacement is prescribed at the
interface between fluid and structure and, successively, extended to the interior of the fluid to identify the deformed fluid domain, $\Omega_d$ 
\cite{Bazilevs2013,BungartzSchafer2006,Richter2017}.  

To formalize shape morphing in a more computationally practical way, we collect the displacements\footnote{In the following we will use the word \emph{displacement} to address the input to the deformation process, and the word \emph{deformation} to address the output of the deformation process.} assigned at 
the control points in the 
vector $\mathbf{d}_c\in \mathbb{R}^{\mathcal{N}_c}$, with 
$[\mathbf{d}_c]_i=d(\mathbf{c}_i)$ and $i=1, \ldots, \mathcal{N}_c$.
In practice, values $d(\mathbf{c}_i)$ are often constrained to 
suitable ranges in order to satisfy admissible shape configurations (see, e.g., \cite{LassilaRozza2010}).
Then, the deformation $\mathbf{d}\in \mathbb{R}^{\mathcal{N}_h}$
of points $\mathbf{x}_i$ in $\Omega_h$, with $[\mathbf{d}]_i=d(\mathbf{x}_i)$ and $i=1, \ldots, \mathcal{N}_h$,
is computed via the IDW interpolant \eqref{IDW_formula}-\eqref{weighting_function} as
\begin{equation}\label{IDW_matrix_form}
{\mathbf{d}} = W \mathbf{d}_c,
\end{equation}
with $W \in \mathbb{R}^{\mathcal{N}_h \times \mathcal{N}_c}$ the IDW matrix
of generic component
\begin{equation}\label{IDW_matrix_comp}
W_{ik} = w_k(\mathbf{x}_i) \quad \text{for } i=1, \hdots \mathcal{N}_h \text{ and } k = 1, \hdots, \mathcal{N}_c.
\end{equation}
The matrix $W$ keeps track of the internal structure of $\Omega_h$, recording
the reciprocal distance between control and internal points. For this reason, the IDW matrix is computed before morphing takes place, once and for all. The actual motion is imposed by vector $\mathbf{d}_c$, which, vice versa,
varies during the morphing process in order to follow
the shape deformation. Finally, vector ${\mathbf{d}}$ identifies a discretization of the deformed domain $\Omega_d$.

Relation \eqref{IDW_matrix_form} highlights the ease of implementation of a morphing strategy driven by IDW
interpolation. In particular, when applied to FSI configurations, this approach 
provides a sharp description of the interface displacement, by properly tackling also portions of the domain characterized by a null displacement and by avoiding compenetrations effects. 
Additionally, a good mesh quality is usually guaranteed (we refer to Section~\ref{sec:app} for more details), even in the presence of large deformations. \\
Nevertheless, the standard IDW approach does not prescribe, a priori, any smart criterion to select the
control points. The value $\mathcal{N}_c$ may consequently become very large, especially in the presence of practical configurations, so that 
the computational effort required by the
assembly and by the storage of matrix $W$ may be very massive. Matrix $W \in \mathbb{R}^{\mathcal{N}_h \times \mathcal{N}_c}$ is actually dense, and its storage may lead to extremely high memory requirements.
This justifies the proposal of the new IDW formulation in the next section, where a new set $\widehat {\mathcal C}$ of control points is properly selected via an automatic procedure, so that 
${\rm card}(\widehat {\mathcal C})\ll \mathcal{N}_c$, 
${\rm card}({\mathcal S})$ denoting the cardinality of the generic set ${\mathcal S}$.
Sparsification of matrix $W$ provides an
alternative way to reduce the computational burden of IDW~\cite{ripepi15,romanelli2012}, though this is beyond the scope of this work. Notice that the method proposed here is complementary to a sparsification approach, so that the two approaches can be, in principle, combined to further reduce the computational effort.

\subsection{The Selective Inverse Distance Weighting (SIDW) approach}\label{sec:sidw}
The new procedure proposed in this section aims at reducing the computational
effort demanded by the standard IDW interpolation.
In particular, since the most computationally demanding part of the IDW algorithm is the memory storage, we aim
at reducing the number of control points by automatically selecting the most relevant points in ${\mathcal C}$
to sharply describe the initial configuration $\Omega_h$ as well as the deformed domain $\Omega_d$. 
For this reason, we refer to the new approach as to the Selective IDW (SIDW) formulation.

The starting point is the approach proposed in~\cite{shepard1968} which is essentially suited to deal
with regimes of small deformations. Our goal is to improve this procedure to tackle also 
large displacement configurations, without violating the no-compenetration constraint.
Additionally, we aim at guaranteeing a uniform distribution in the reduced set of control points.

The complete SIDW procedure is provided by \textbf{Algorithm 1}.
To set notation, we denote by: $\widehat {\mathcal C}=\{ {\widehat{\mathbf{c}}}_j \}\subset {\mathcal C}$ the 
subset of the selected control points; 
$\omega\subseteq \Omega_r$ the region where the selection is applied; $R > 0$ the selection radius; 
$B(\textbf{c}; r)$ the ball of center $\textbf{c}$ and radius $r > 0$; 
$B(\textbf{c}; r_i, r_e)$ the circular annulus of center $\textbf{c}$, inner radius $r_i$ and
outer radius $r_e$, with $r_e > r_i > 0$.\\
We first exemplify the action of the SIDW procedure starting from a configuration characterized by $\Omega_r=(0,10)^2$ and $\omega\equiv \Omega_r$.
Then, we will particularize such a procedure to a shape morphing setting.\\
\begin{figure*}
\centering
    \includegraphics[width=0.3\textwidth]{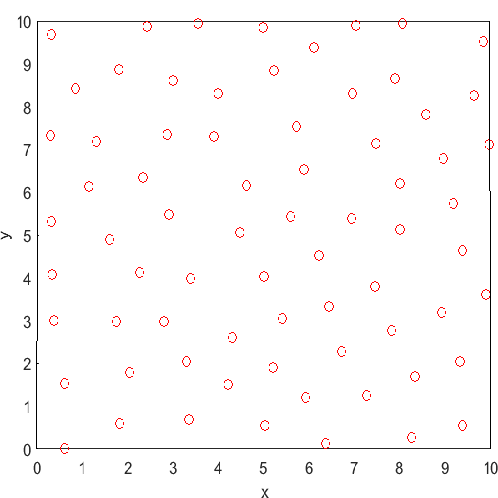}
\caption{SIDW algorithm: distribution of control points after the selection process on the square $\Omega_r=(0,10)^2$, being $\omega \equiv \Omega_r$.}\label{es1_prepost}
\end{figure*}
SIDW algorithm consists essentially of two phases, i.e., an \emph{initialization} phase and the actual \emph{selection} phase.
Now, let $\widehat{\mathbf{c}}_1\in {\mathcal C}$ be a randomly selected control point. 

\begin{figure*}
\centering
\subfigure[Initialization phase.]{
\includegraphics[width=0.4\textwidth]{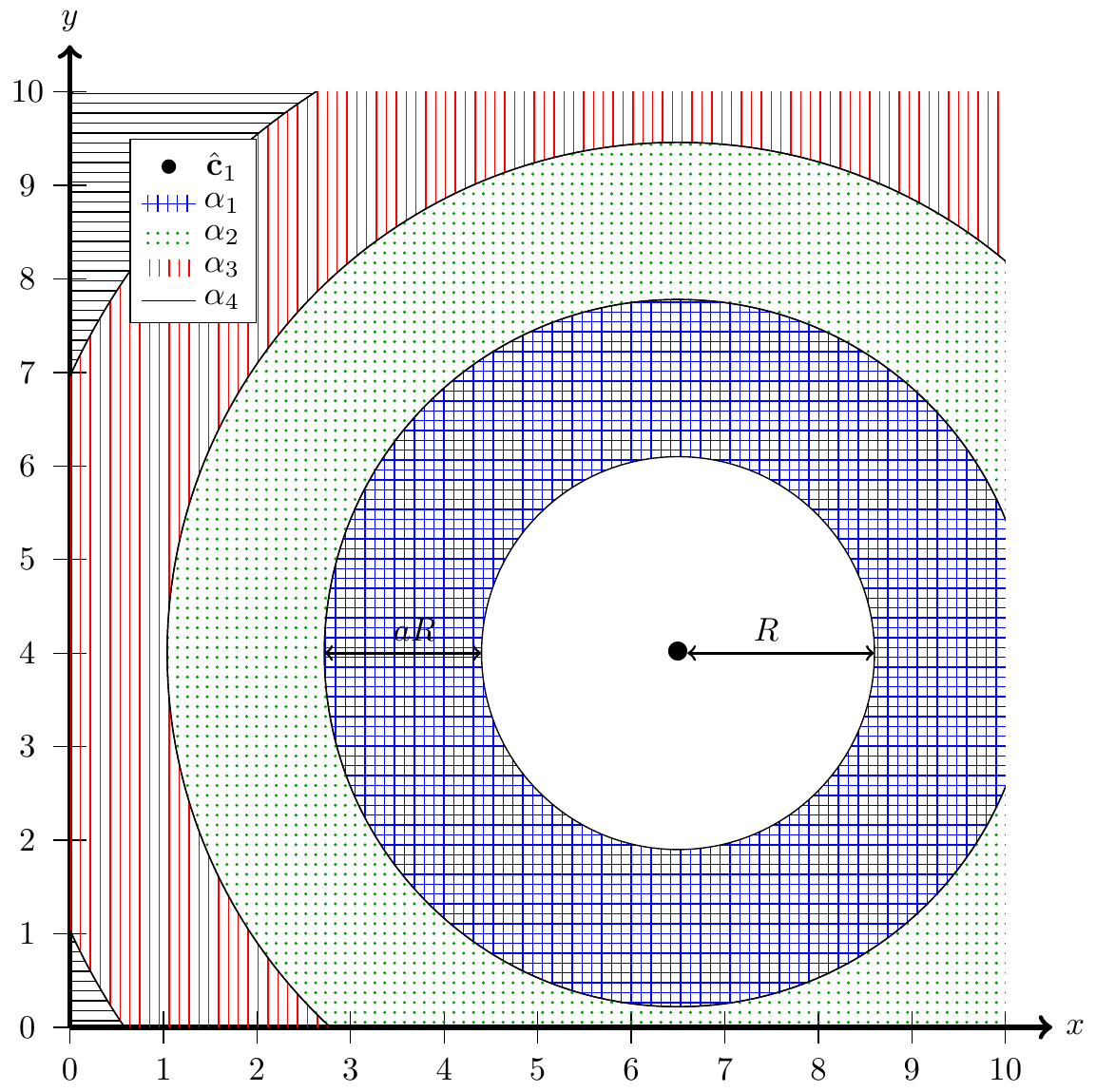}\label{CT1_1}
}
\subfigure[Selection phase: identification of $\beta$ and of $\widehat{\mathbf{c}}_{2}$.]{
\includegraphics[width=0.4\textwidth]{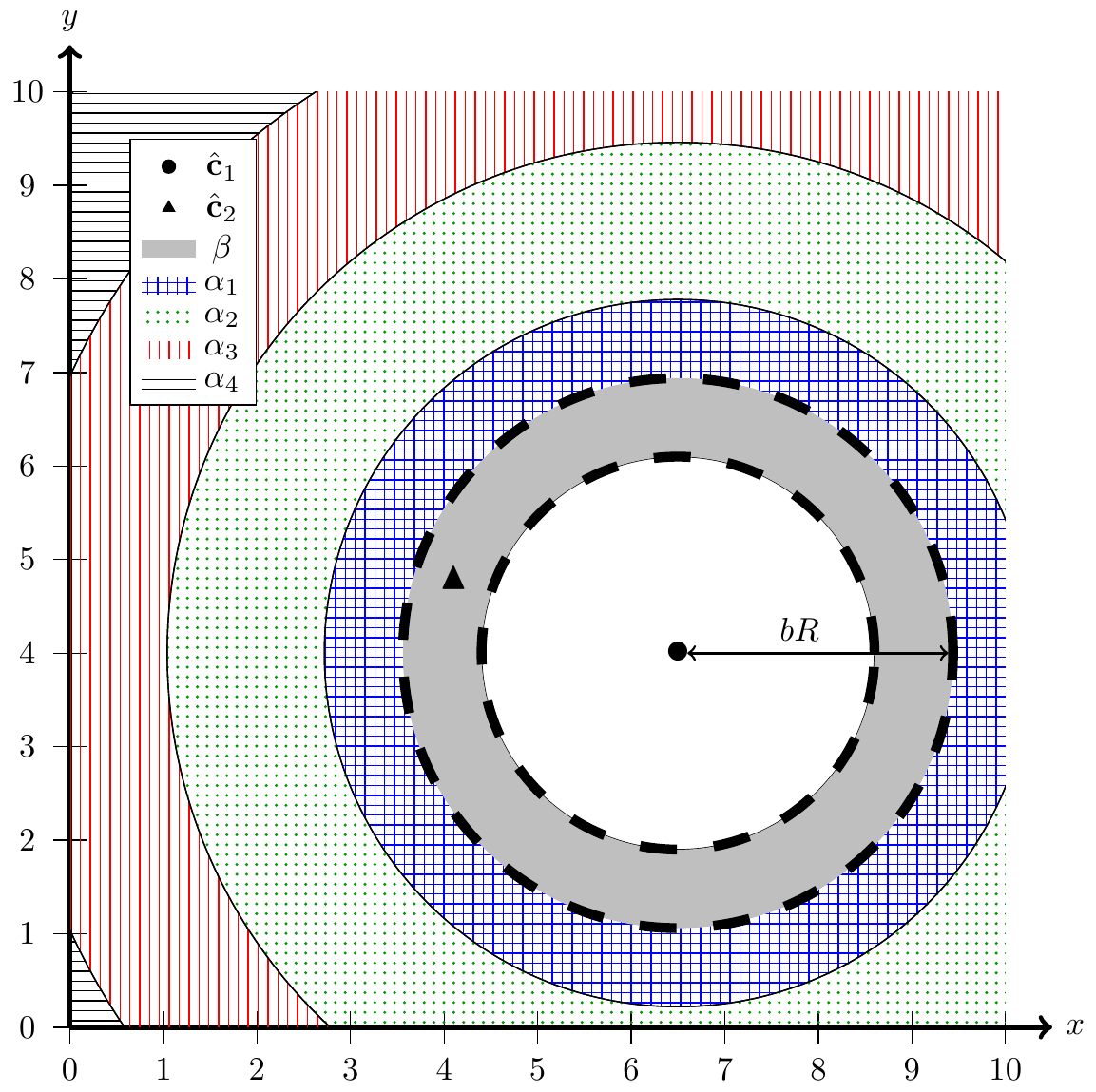}\label{CT1_2}
}
\subfigure[Selection phase: first removal associated with $\alpha_1$.]{
\includegraphics[width=0.4\textwidth]{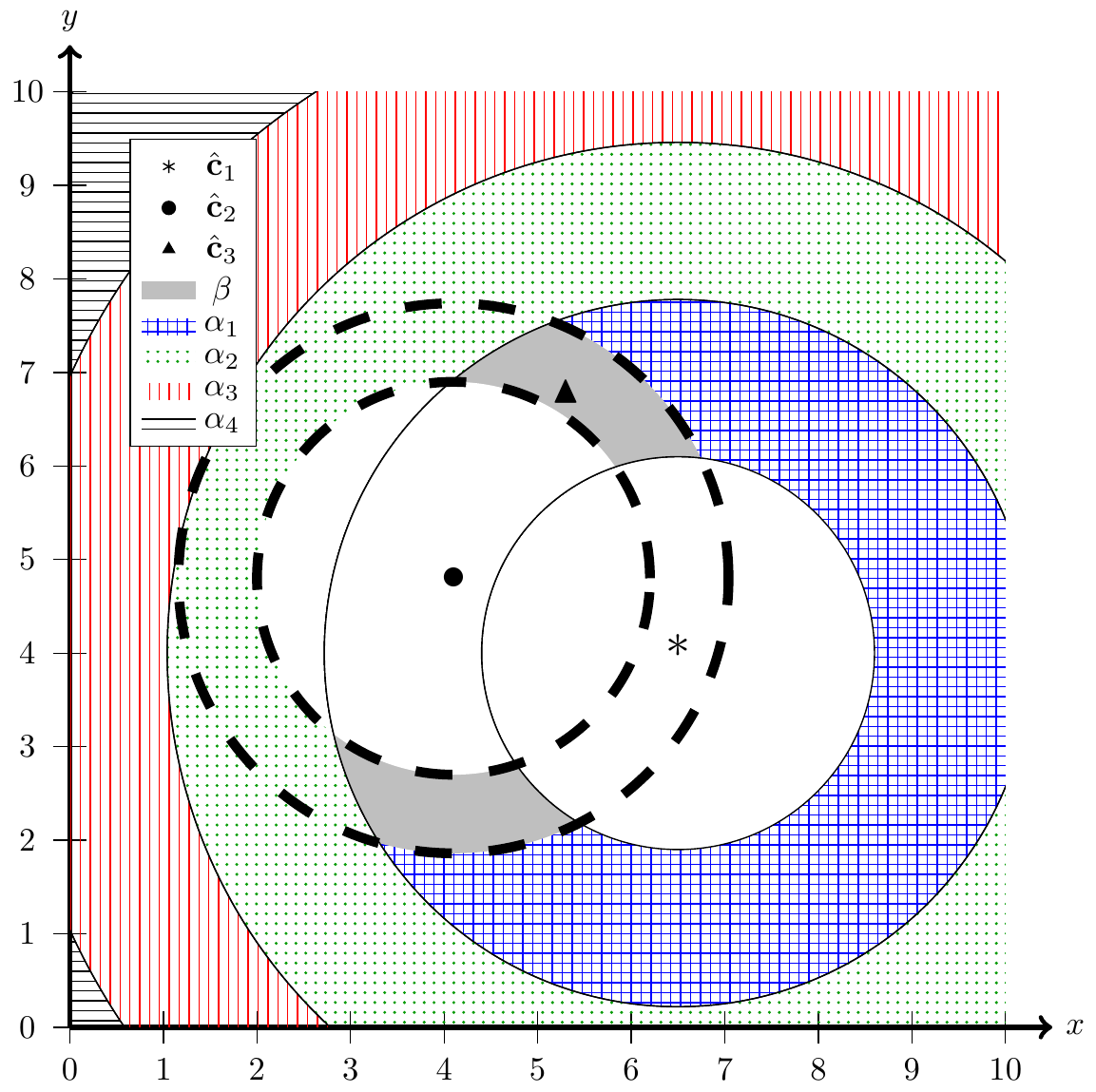}\label{CT1_3}
}
\subfigure[Selection phase: second removal associated with $\alpha_1$.]{
\includegraphics[width=0.4\textwidth]{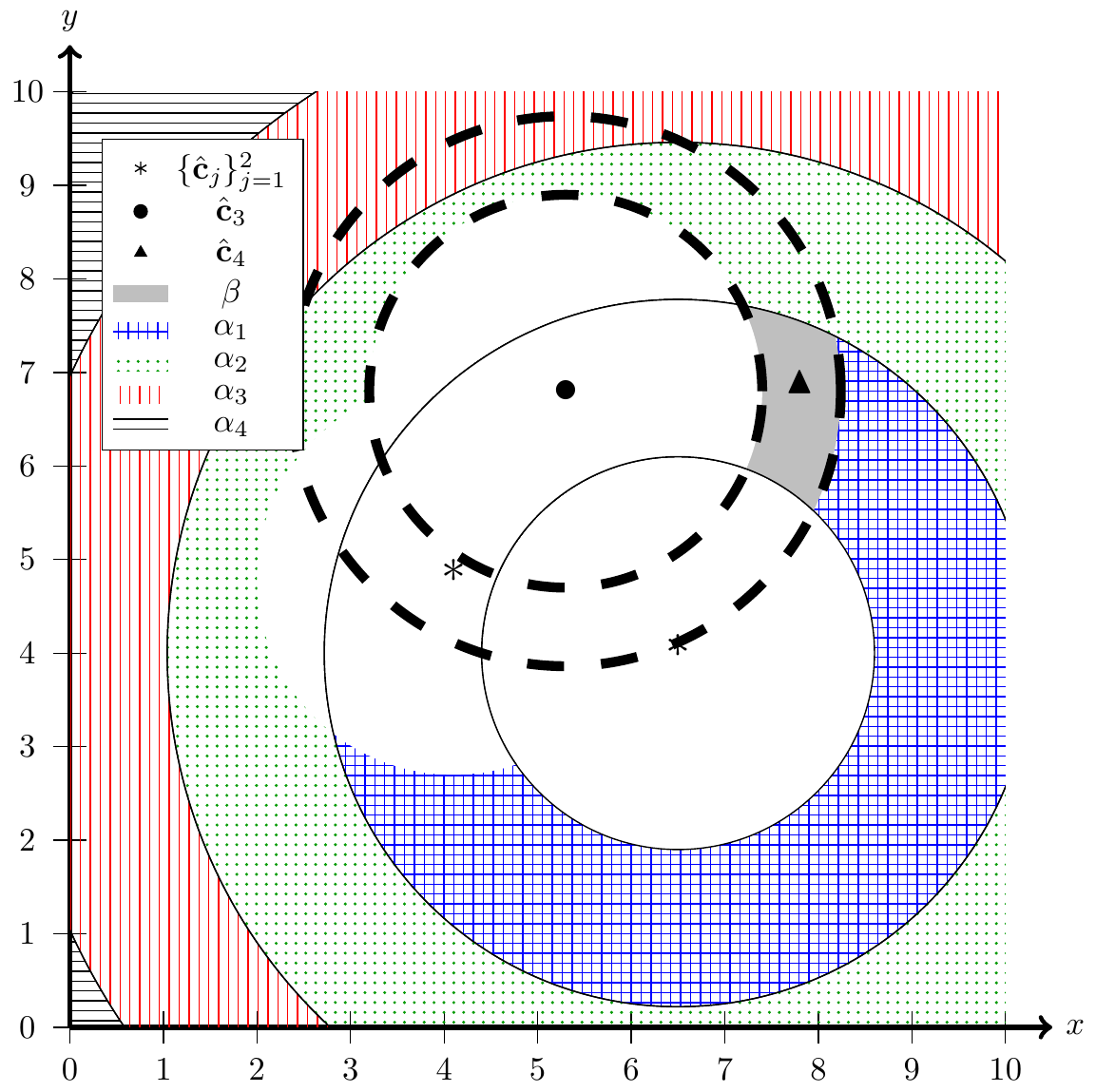}\label{CT1_4}
}
\subfigure[Selection phase: sixth removal associated with $\alpha_1$.]{
\includegraphics[width=0.4\textwidth]{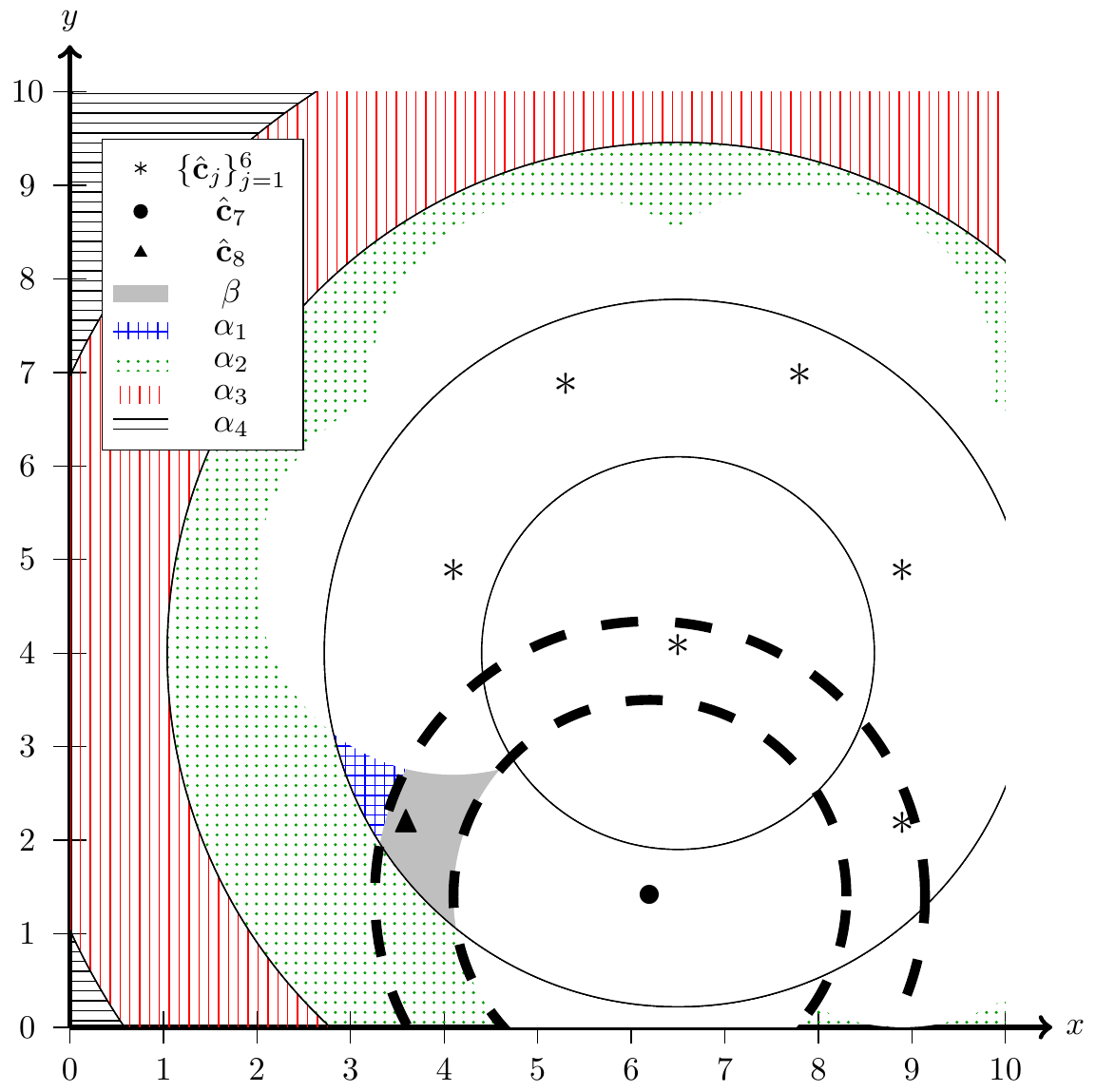}\label{CT1_8}
}
\subfigure[Selection phase: first removal associated with $\alpha_2$.]{
\includegraphics[width=0.4\textwidth]{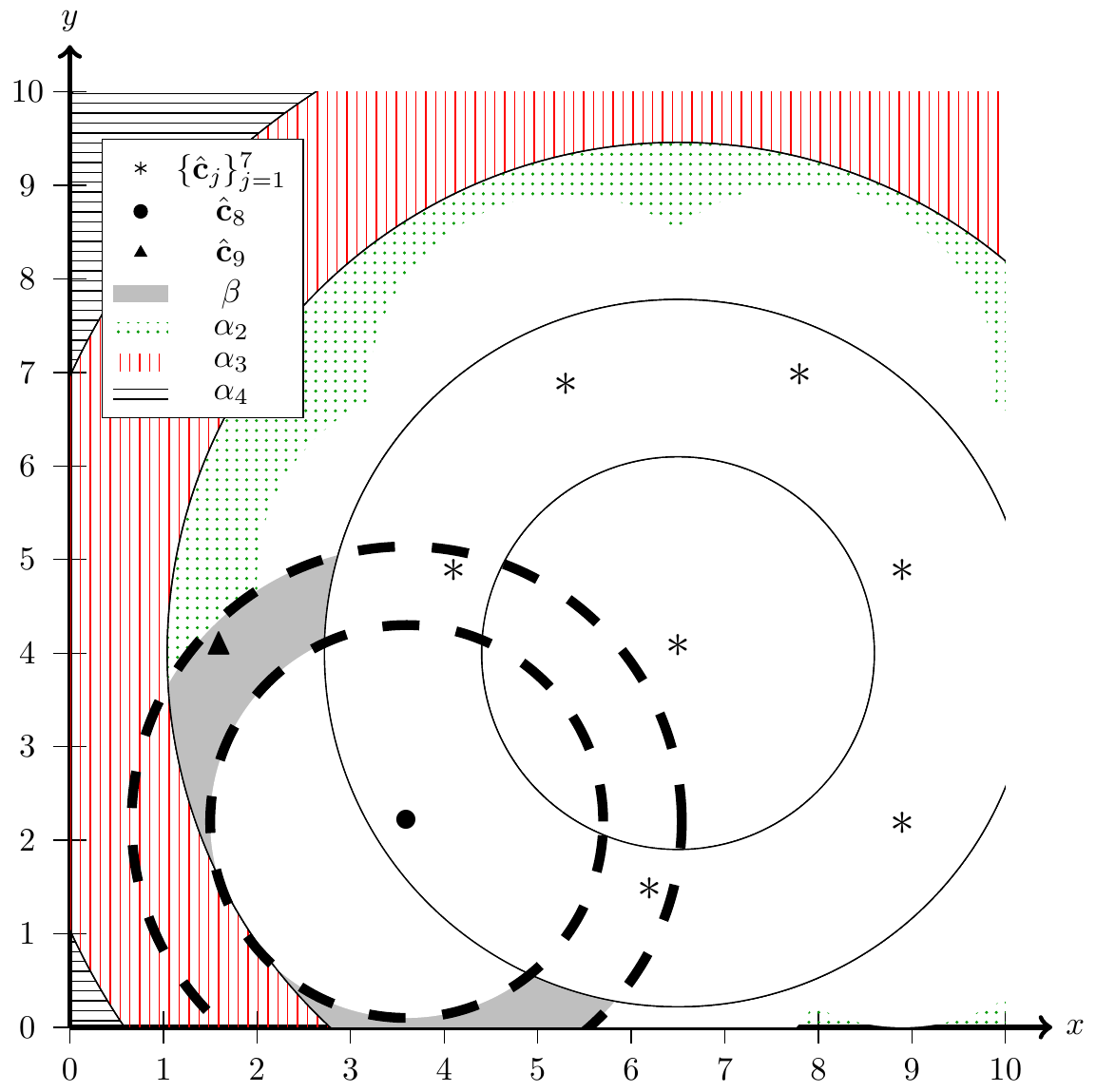}\label{CT1_9}
}
\caption{SIDW algorithm: initialization (a) and first iterations (b)-(f) of the selection phase applied to the
square $\Omega_r=(0,10)^2$, being $\omega \equiv \Omega_r$.}\label{CT1}
\end{figure*}

The goal of the initialization phase is to perform a tessellation 
of the region $\omega \setminus B(\widehat{\textbf{c}}_1; R)$ via $n$ concentric annuli,
$\{\alpha_m\}_{m=1}^{n}$, of thickness $a R$ ({\tt lines 2--7}), being $a\in \mathbb{R}$ a positive constant picked by the user, with $a < 1$. Concerning the specific example on the square $\Omega_r=(0,10)^2$, we refer to Figure~\ref{CT1}(a), where the tessellation corresponding to $R=2.1$ and $a=0.8$
is shown, the control point $\widehat{\mathbf{c}}_1$ being marked in the plot. To simplify the graphical representation,
we highlight only the control points constituting $\widehat {\mathcal C}$.
Four annular regions are identified by the choice done for the input parameters, thus inducing a partitioning of the original control points in $\mathcal C$. Then, 
all the control points in $\mathcal C$ belonging to the \emph{area of influence}, $B(\widehat{\textbf{c}}_1; R)$, 
around $\widehat{\mathbf{c}}_1$ are removed from $\mathcal C$. This operation concludes the first phase of SIDW algorithm.
\begin{algorithm*}
\normalsize
\caption{: SIDW interpolation}\label{alg:selection}
\begin{algorithmic}[1]
\Function{SIDW}{$\omega$, ${\mathcal C}$, $\widehat{\textbf{c}}_1$, $R$, $a$, $b$, $u$}
    \State $n \gets 1$, $r \gets R$, $r_{\omega} \gets \min\{r > 0: \omega \supseteq B(\widehat{\textbf{c}}_1; r)\}$
    \While{$r \leq r_{\omega}$}
        \State $\alpha_n \gets B(\widehat{\textbf{c}}_1; r, r + a R) \cap \mathcal{C}$
        \State $r \gets r + a R$
        \State $n \gets n + 1$
    \EndWhile
    \State $k \gets 1$, $m \gets 1$, $\beta \gets \alpha_1 \cap B(\widehat{\mathbf{c}}_1; R, b R)$
    \While{$m \leq n$}
        \While{$\beta \neq \emptyset$}
            \State $\widehat{\mathbf{c}}_{k + 1} \gets \text{SelectControlPoint}(\beta)$
            \For{$l=m, \hdots, n$}
                \State $\alpha_l \gets \alpha_l \setminus B(\widehat{\mathbf{c}}_{k + 1}; R)$
            \EndFor
            \State $\beta \gets \alpha_m \cap B(\widehat{\mathbf{c}}_{k+1}; R, b R)$
            \State $k \gets k + 1$
        \EndWhile
        \If{$\alpha_m = \emptyset$}
            \State $\beta \gets \alpha_{m+1} \cap B(\widehat{\mathbf{c}}_{k}; R, b R)$
            \State $m \gets m + 1$
        \Else
            \State $\beta \gets \alpha_m$
        \EndIf
    \EndWhile
    \State $\widehat {\mathcal C} \gets \{\widehat{\mathbf{c}}_j\}_{j=1}^{k}$
    \State computation of $\Pi_{\rm SIDW}(u)$
\EndFunction
\end{algorithmic}
\end{algorithm*}

The selection phase begins with a first initialization of the \emph{area of selection}, $\beta$, around the first control point $\widehat{\textbf{c}}_1$, defined as the intersection between $\alpha_1$ and the 
circular annulus $B(\widehat{\mathbf{c}}_1; R, b R)$, with $b > 1$ a user defined positive constant ({\tt line 8}).
Then, a loop on the annular regions $\{\alpha_m\}_{m=1}^n$ of the tessellation is started ({\tt line 9}). 
Inside this loop, until $\beta$ is empty, a new control point is selected in $\beta$ via the {\tt SelectControlPoint} function ({\tt line 11}).
Then, all the control points in $\mathcal C$ belonging to each area of influence $\alpha_l$
around the current control point $\widehat{\mathbf{c}}_{k+1}$ are removed ({\tt lines 12--14}), and a 
new area of selection $\beta$ is computed as the intersection between the current area of influence
$\alpha_m$ and a circular annulus $B(\widehat{\mathbf{c}}_{k+1}; R, b R)$ ({\tt line 15}).\\
Concerning function {\tt SelectControlPoint}, different criteria can be implemented to select the new control point  
$\widehat{\mathbf{c}}_{k+1}$. For instance, a geometric selection can be applied by picking 
the control point closest to the centroid of $\beta$ or the farthest control point with respect 
to the already selected points in $\widehat {\mathcal C}$. To simplify the procedure, in the numerical assessment below, we resort to a random choice of $\widehat{\mathbf{c}}_{k+1}$ in $\beta$. \\
Figure~\ref{CT1}(b) displays the  
firstly initialized area of selection $\beta$, together with the second control point $\widehat{\mathbf{c}}_2$, being 
$b=1.4$. The remaining panels in Figure~\ref{CT1} show the effects of the successive removals operated by the {\tt while} loop.
Six iterations are associated with the first area of influence $\alpha_1$. 
In particular, the new area of selection $\beta$ and the corresponding control point $\widehat{\mathbf{c}}_{k+1}$ associated 
with the first two (Figure~\ref{CT1}(c)-(d)) and the last (Figure~\ref{CT1}(e)) removals are highlighted.\\
When $\beta$ is empty, a re-initialization of the area of selection is performed ({\tt lines 18--23}).
If the current annulus $\alpha_m$ is exhausted, then the area of selection is sought in the next annulus, $\alpha_{m+1}$, 
({\tt lines 19--20}). Figure~\ref{CT1}(f) shows the transition from $\alpha_1$ to $\alpha_2$. Additionally,
 {\tt line 22} handles the very peculiar case when $\alpha_m$ is not actually exhausted, yet it is disconnected from the area of selection around the current control point $\widehat{\mathbf{c}}_{k+1}$. In such a case, $\beta$ is simply reset to $\alpha_m$. 
 
Finally, once all the circular annuli in the tessellation $\{\alpha_m\}_{m=1}^{n}$ are exhausted,
SIDW algorithm ends, and the new set $\widehat {\mathcal C}$ of selected control points is returned ({\tt line 25}).
The final distribution of control points for the considered specific configuration on the square $\Omega_r=(0,10)^2$ is provided in Figure~\ref{es1_prepost}.
The selected points are uniformly distributed as desired. At this point, the \emph{SIDW interpolant} of $u$, $\Pi_{\rm SIDW}(u)$, can be computed 
({\tt line 26}) as
\begin{equation}\label{SIDW_interp}
\Pi_{\rm SIDW}(u)(\mathbf{x}) =
\sum _{k=1}^{\mathcal{N}_{\widehat c}}  w_k(\mathbf{x}) \, u(\widehat{\mathbf{c}}_k)\quad \mathbf{x}\in \Omega_r,
\end{equation}
with ${\rm card}(\widehat {\mathcal C})={\mathcal{N}_{\widehat c}}\ll {\mathcal{N}_c}$, and where the weight functions $w_k$
are defined according to \eqref{weighting_function} having replaced $\mathbf{c}_k$ with $\widehat{\mathbf{c}}_k$.

A few remarks are in order. The three input parameters, $R$, $a$ and $b$, tune the SIDW algorithm and, as a consequence, 
the final selection $\widehat {\mathcal C}$ of control points. In particular, they vary the area of 
the annular regions $\{\alpha_m\}_{m=1}^{n}$, as well as of the areas of influence and selection.
In more details, concerning the radius $R$, it turns out that
the smaller $R$, the larger the number of selected control points. 
As far as $a$ and $b$ are concerned, smaller values of $a$ lead to a larger number $n$ of annular regions $\alpha_m$,
which generally results in a finer control point distribution; larger values of $b$ leads to larger areas of selection $\beta$, resulting in a smaller number of selected control points.
The choice of $R$, $a$ and $b$ as well as the corresponding interplay is 
difficult to be established \emph{a priori} and, clearly, it is problem-dependent.
A sensitivity analysis of SIDW with respect to these parameters will be carried out in Section~\ref{sec:app}.

Finally, we remark that, especially in practical applications, it might be advisable to split $\omega$ into smaller subregions, $\omega_p$, and to use various values for the parameters $R$, $a$ and $b$ to tackle possible local different discretizations.

\subsubsection{SIDW for shape morphing}
Now we customize SIDW algorithm to manage shape morphing.
We consider the discretization $\Omega_h=\{\mathbf{x}_i\}_{i=1}^{\mathcal{N}_h}$
of the initial configuration of the physical domain and we assume to know the displacement $d$ at each control point
$\mathbf{c}_k\in {\mathcal C}$, for $k=1, \ldots, \mathcal{N}_c$. We remind that, in such a case, 
${\mathcal C}$ is a subset of boundary nodes in $\Omega_h$.

Our idea consists in deforming each boundary node of $\Omega_h$
according to the deformation prescribed at $\{ \mathbf{c}_k\}_{k=1}^{\mathcal{N}_c}$, while
selecting a subset of control points to update the deformation of the internal nodes of $\Omega_h$.\\
For this purpose, we resort to the procedure
itemized in \textbf{Algorithm 2}, following a selection-deformation paradigm.\\
During the selection phase, which is performed once and for all, we first filter set ${\mathcal C}$ via the selection
procedure in {\bf Algorithm \ref{alg:selection}},  to obtain the subset $\widehat {\mathcal C}$ of control points
$\widehat{\mathbf{c}}_k$, with $k=1, \ldots, \mathcal{N}_{\widehat c}$ ({\tt item (i)}), using $\omega = \partial\Omega_r$, a random point $\mathbf{c}_1 \in {\mathcal C}$ and considering the displacement as the variable $u$ to be interpolated. Then, we build the 
SIDW matrix $\widehat W$ according \eqref{IDW_matrix_comp} 
associated with $\widehat {\mathcal C}$ and with the set of the internal nodes of 
$\Omega_h$ ({\tt item (iii)}). As an option, we can add to the filtered set $\widehat {\mathcal C}$
extra control points to include possible specific constraints to the problem ({\tt item (ii)}). This occurs, for instance, 
when a null displacement is assigned to a portion of the domain. In such a case, the fixed nodes have to be necessarily included in the set $\widehat {\mathcal C}$. 
User-defined and problem-specific constraints may be added as well, for instance to attempt to avoid poor-quality and inverted elements in case some boundary points
undergo a small deformation, whereas nearby boundary points are deformed significantly.
In the sequel, we refer to this variant of the SIDW algorithm as to 
ESIDW (\emph{Enriched} SIDW) interpolation.\\
The online phase performs the actual shape morphing. Each boundary node of $\Omega_h$ is deformed via the displacement 
assigned at the control points in ${\mathcal C}$ ({\tt item (iv)}). Successively, we deform the internal nodes of 
$\Omega_h$ as
\begin{equation}\label{SIDW_matrix_form}
\mathbf{d}^i= {\widehat W} \mathbf{d}_{\widehat c},
\end{equation}
where $\mathbf{d}^i\in \mathbb{R}^{\mathcal{N}_h^i}$ collects the deformation of the internal nodes, 
$\mathbf{x}_j$ for $j=1, \ldots, \mathcal{N}_h^i$, of $\Omega_h$, 
${\widehat W}\in \mathbb{R}^{\mathcal{N}_h^i\times \mathcal{N}_{\widehat c}}$ is the (E)SIDW matrix computed in the offline phase, and $\mathbf{d}_{\widehat c} \in \mathbb{R}^{\mathcal{N}_{\widehat c}}$ is the vector of the displacements of the selected control points in $\widehat {\mathcal C}$ ({\tt item (v)}).

\textbf{Algorithm 2} can be advantageously exploited to deal with FSI problems. In such a case, we identify domain $\Omega_r$
with the fluid domain. The displacement of the interface is generally provided by a structure solver, so that it suffices 
to compute the deformation of the fluid nodes via the procedure detailed in the algorithm.
\begin{algorithm*}
\normalsize
\caption{SIDW interpolation for shape morphing}\label{alg:sidw_shape_morphing}
\begin{algorithmic}[999]
\State \textsc{Selection phase}: 
    \State \qquad (i) apply the selection procedure in {\bf Algorithm \ref{alg:selection}} to ${\mathcal C}$ to extract the subset 
    $\widehat {\mathcal C}$;
    \State \qquad (ii) \emph{Optional [ESIDW method]} : add to the set $\widehat{\mathcal C}$ extra control points to account for problem constraints;
    \State \qquad (iii) assemble the SIDW matrix $\widehat W$ associated with $\widehat {\mathcal C}$ and with the internal nodes of $\Omega_h$.\\

\State \textsc{Deformation phase}: 
    \State \qquad (iv) deform boundary nodes of $\Omega_h$ via the displacement assigned at the control points in 
    ${\mathcal C}$;
    \State \qquad (v) deform internal nodes of $\Omega_h$ according to \eqref{SIDW_matrix_form}.
\end{algorithmic}
\end{algorithm*}

\subsection{IDW versus SIDW interpolation}\label{sec:app}
In this section we investigate the performances of IDW and SIDW interpolation algorithm.
In particular, we consider configurations of interest in aeronautic and naval engineering. 
We focus on 3D test cases and on tetrahedral meshes, although the procedure can be generalized to
any dimension and to arbitrary meshes.
Concerning the choice of the parameter $p$ in \eqref{weighting_function}, we set $p=4$, essentially driven by numerical considerations. In all the test cases, we will consider a shape morphing process where the deformation is driven by an analytic law provided as input to  \textbf{Algorithm 2}. 
As an alternative, we can employ a structural solver to compute the displacement to be assigned to the structure, but this is beyond the goal of this paper.
To check the improvements led by the new approach, we compare IDW and SIDW interpolations 
in terms of computational effort. Both the procedures are implemented in the {\tt C++} open source library {\tt libMesh} \cite{KirkPetersonStognerCarey2006}, while the visualization software {\tt ParaView} is employed
for the post-processing of the solutions \cite{paraview}. Concerning mesh generation, 
we resort to {\tt SALOME} \cite{Salome}. Finally, 
all the simulations are performed on a laptop with Intel\textsuperscript{\textregistered} Core\texttrademark\ 
i7 CPU and 4GB RAM.

\subsubsection{Structural deformation of a wing}\label{sec:ass1}
\begin{figure*}
\centering
\includegraphics[width=0.65\textwidth]{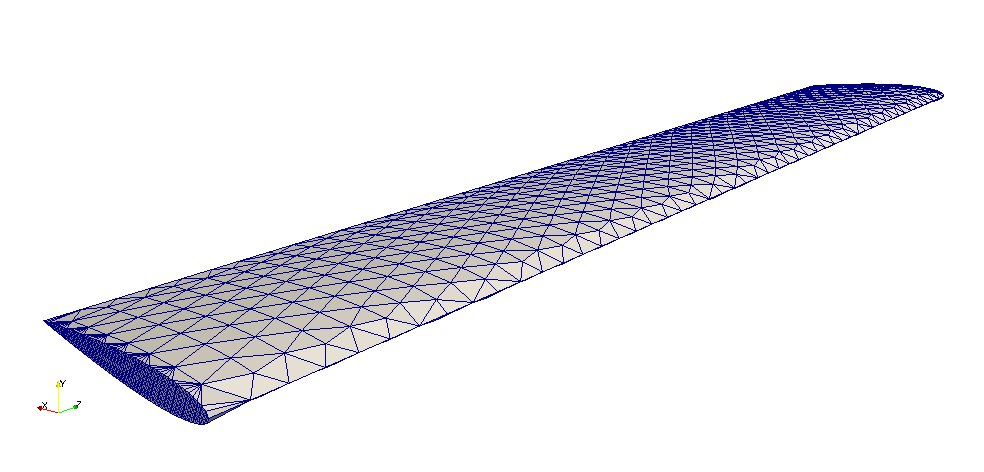}\\
    \includegraphics[width=0.65\textwidth]{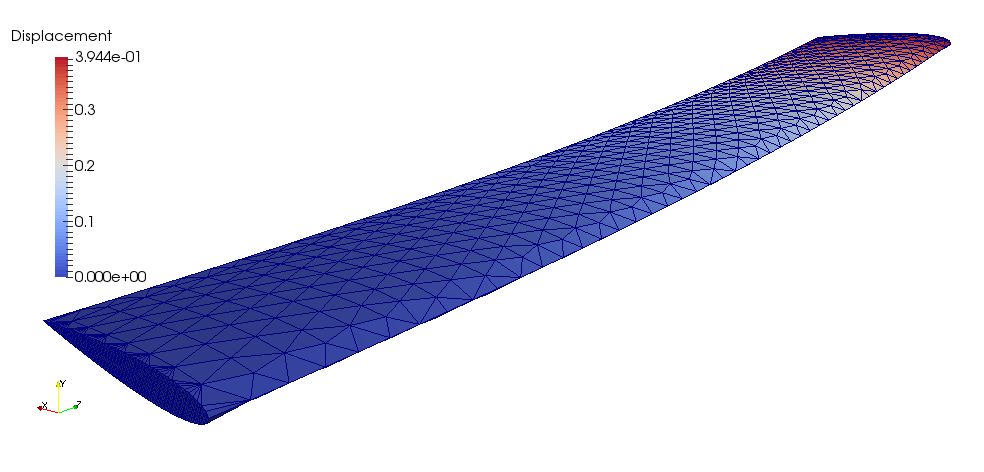}
\caption{Structural deformation of a wing: reference (top) and IDW deformed (bottom) configuration.}\label{wing}
\end{figure*}

We consider a wing characterized by a NACA0012 profile, clamped on the left side (see Figure \ref{wing}, top).
Table \ref{struct_wing_RIDW_table} gathers the main properties
of the reference domain $\Omega_r$ and of the corresponding mesh $\Omega_h$.
We impose a vertical displacement to the boundary of the wing. In particular, 
denoting by $z$ the distance from the clamped side and by $y$ the vertical direction, we assign the displacement
\begin{equation}
\delta y = \delta y(z) = 0.01 z^2
\label{eq:wing_displacement}
\end{equation}
in the $y$ direction  (see Figure \ref{wing}, bottom). The goal of this test case is to check the deformation capabilities 
of IDW and SIDW methods on a simple case of shape morphing, inspired by a possible application to shape optimization problems.
The actual applicability to FSI problems (where mesh motion concerns the fluid domain) will be discussed in Section \ref{sec:ass2}.
\begin{table*}
\centering
\begin{tabular}{|cc|}
\hline
longitudinal dimension & $2\pi$ [m]	\\
\hline
profile chord length & $1.01$ [m]	\\
\hline\hline
\# elements & $8850$	\\
\hline
\# nodes & $2463$ \\
\hline
\# internal nodes & $797$ \\
\hline
\# boundary nodes & $1666$ \\
\hline
\end{tabular}
\caption{Structural deformation of a wing: main properties of $\Omega_r$ and of $\Omega_h$.} \label{struct_wing_RIDW_table}
\end{table*}

First, we apply the standard IDW approach, after identifying ${\mathcal C}$
with the set of all the boundary nodes of $\Omega_h$ (see Figure \ref{wing_idw}, top), so that the IDW 
matrix $W$ in \eqref{IDW_matrix_form} belongs to $\mathbb{R}^{797 \times 1666}$.
The resulting deformed wing, shown in Figure \ref{wing}, bottom, is obtained after\footnote{Here and in the following, unless when referring to a specific phase of the deformation process, the time we refer to includes the time required by \eqref{IDW_matrix_form} as well as the time to update the mesh. For this reason, the actual speedup might be lower than the expected (ideal) one.} 0.25 [s]. 
\begin{figure*}
\centering
    \includegraphics[width=0.65\textwidth]{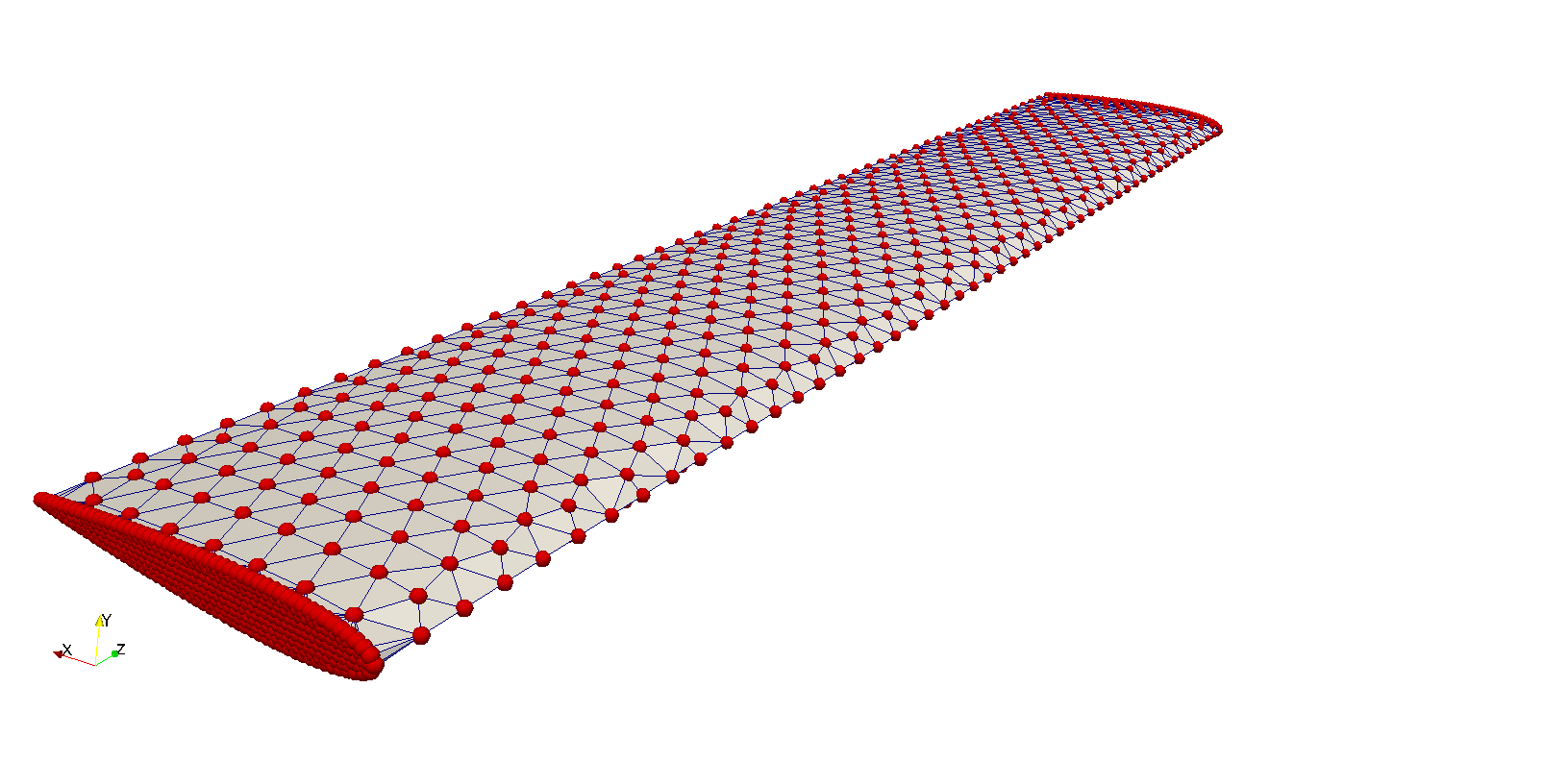}\label{wing_idw}
    \includegraphics[width=0.65\textwidth]{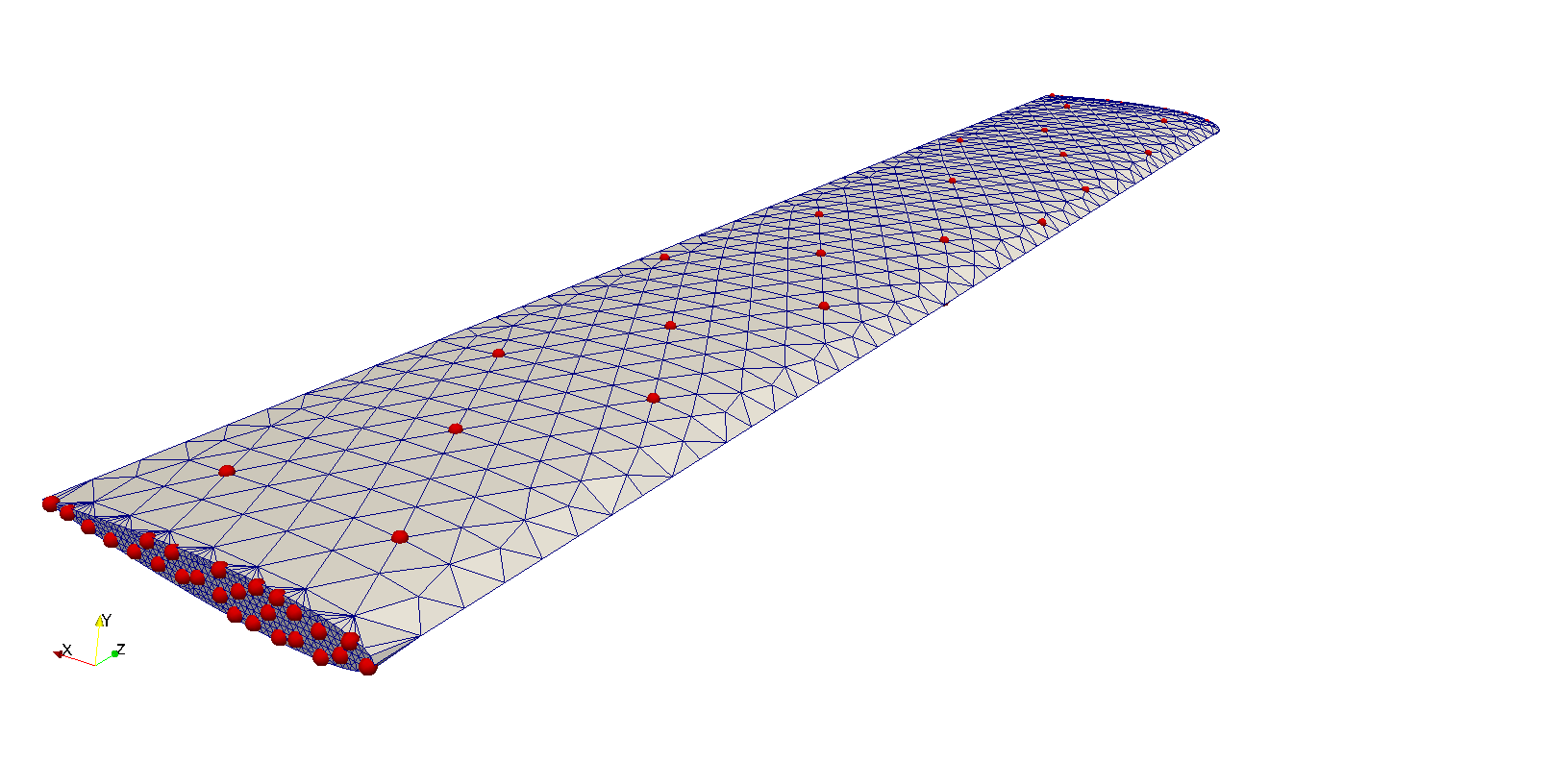}\label{wing_ridwp}
    \includegraphics[width=0.65\textwidth]{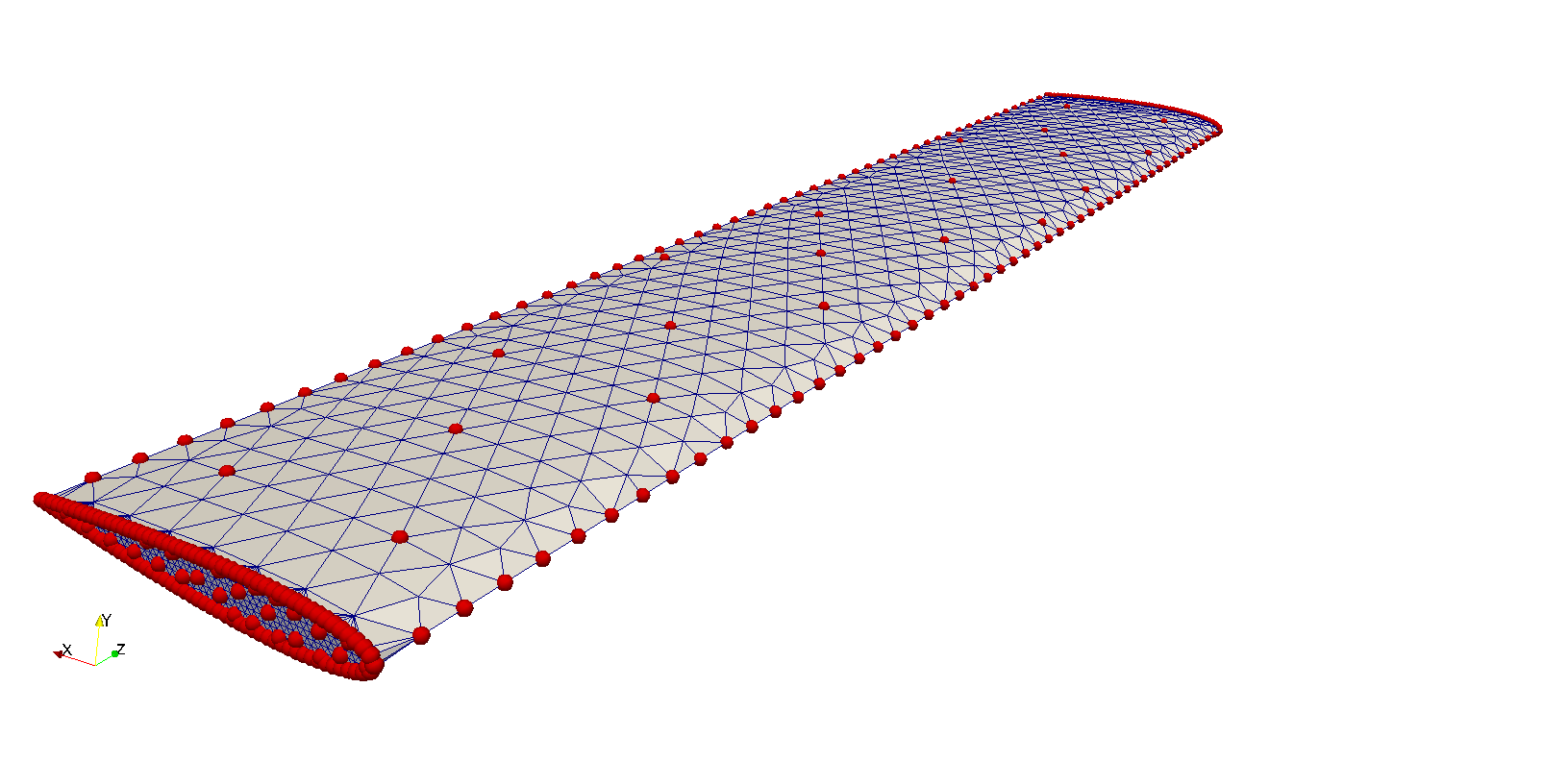}\label{wing_ridw}
\caption{Structural deformation of a wing: control points associated with the IDW (top), 
SIDW (center) and ESIDW (bottom) interpolation.}\label{wing_idw}
\end{figure*}
\begin{figure*}
\centering
    \includegraphics[width=0.313\textwidth]{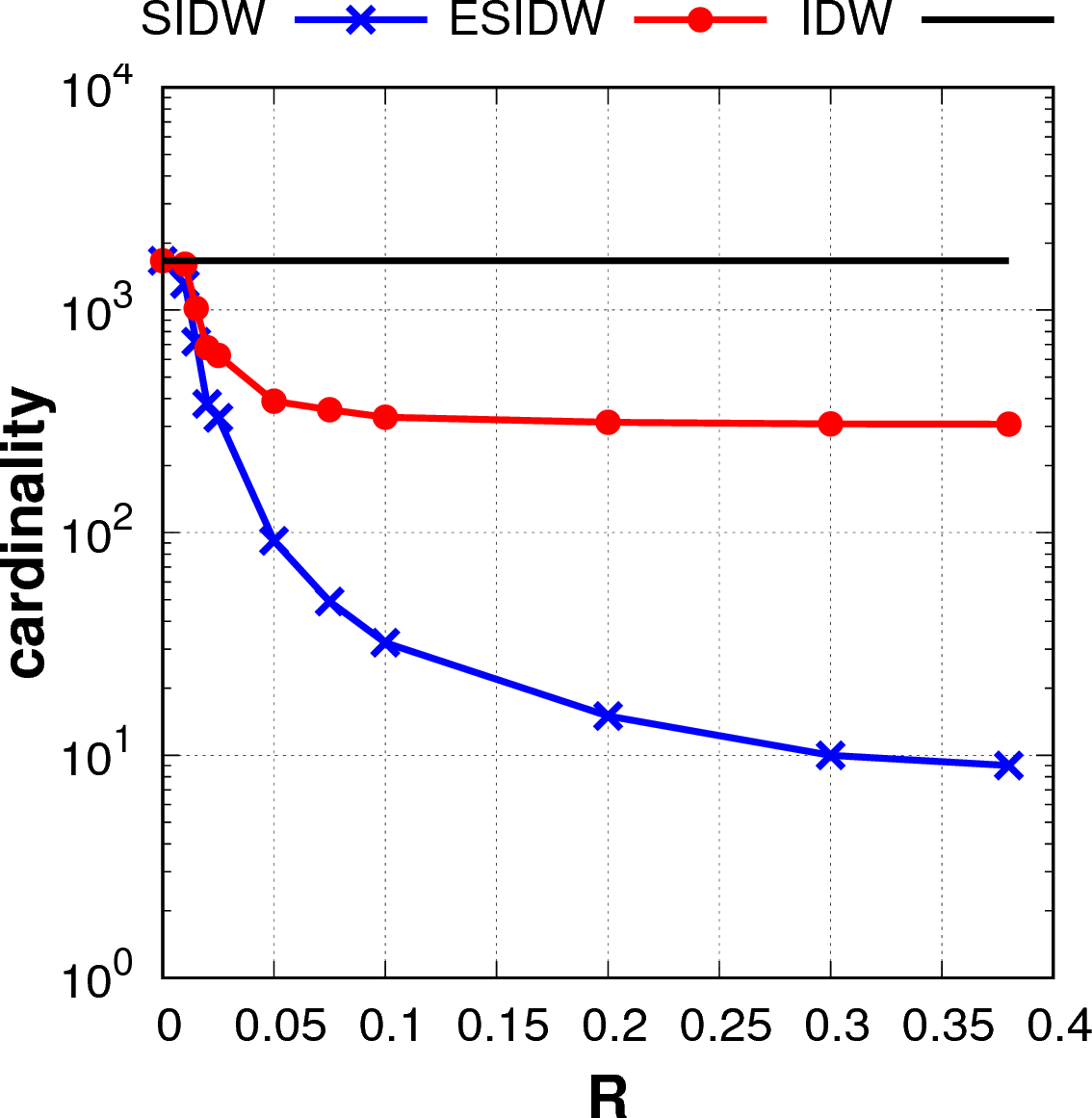}\quad
    \includegraphics[width=0.323\textwidth]{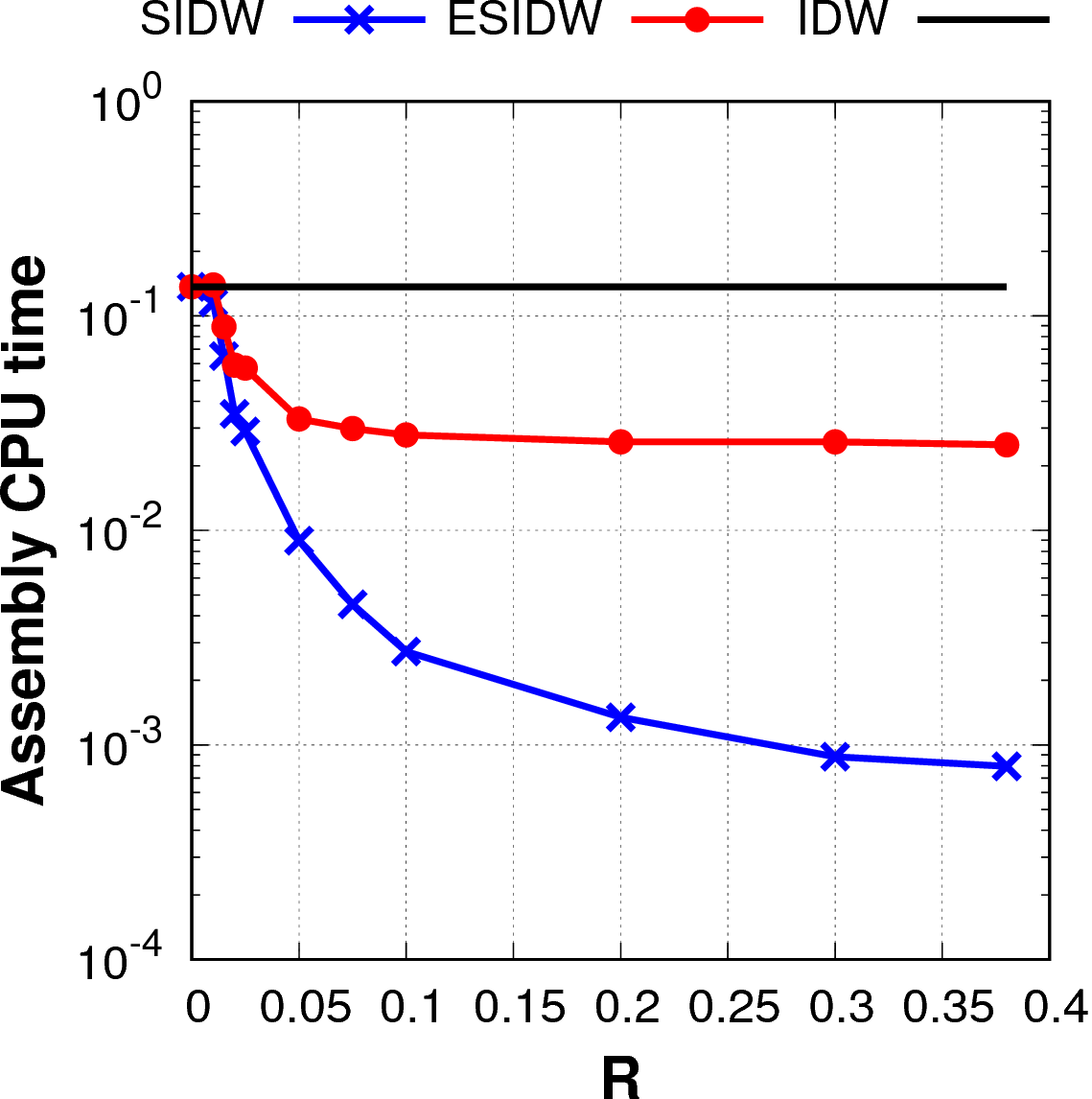}\quad
    \includegraphics[width=0.323\textwidth]{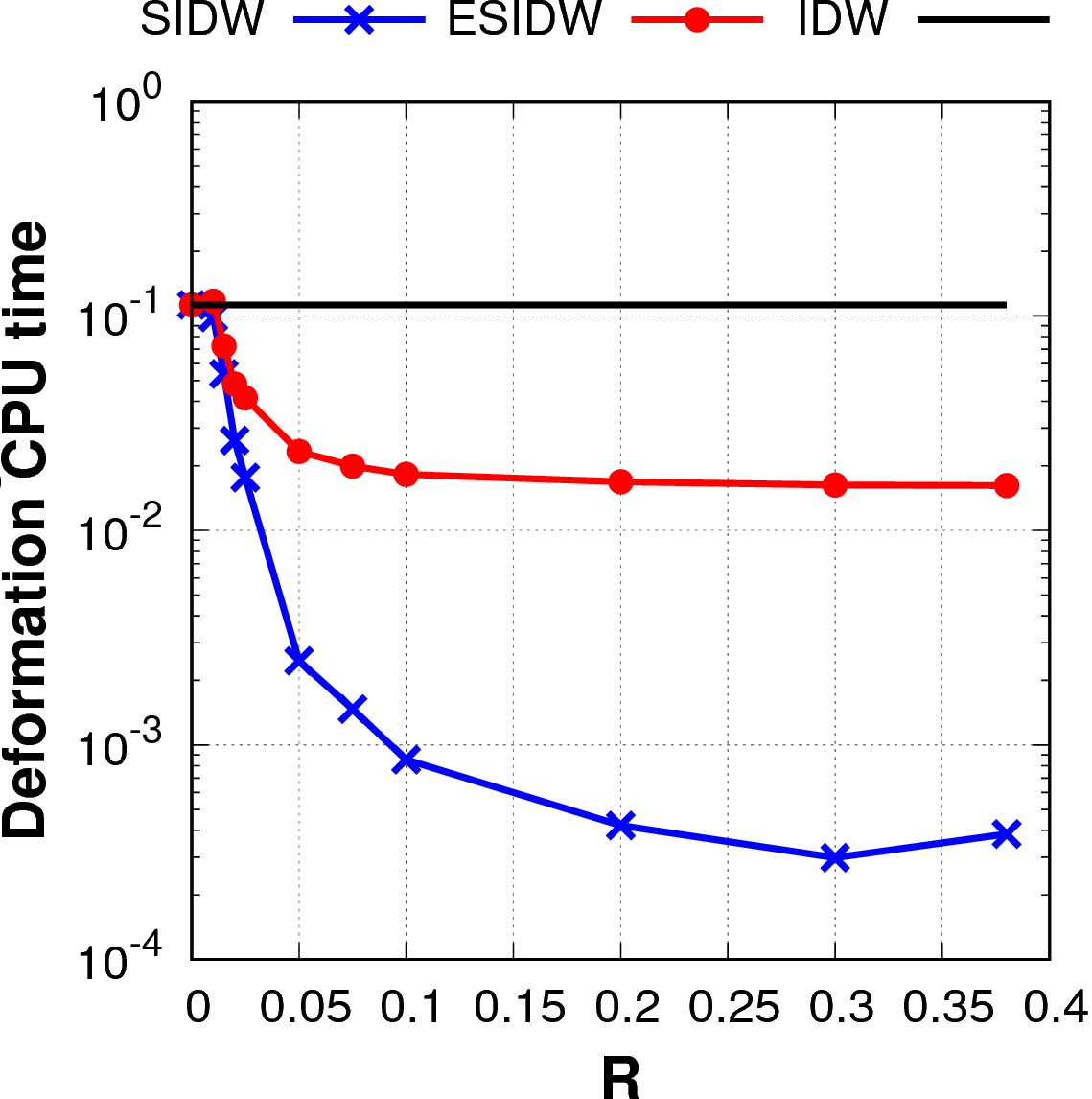}
\caption{Structural deformation of a wing: cardinality of the control point set (left), assembly
(center) and deformation (right) CPU time as a function of $R$.}\label{wing_cpu_analysis}
\end{figure*}

We now resort to the SIDW interpolation algorithm. Set ${\mathcal C}$ still contains all the boundary nodes.
We adopt two different values for the selection radius, i.e., $R=R_{lr}=0.05$ [m]   
for the left and the right lateral surfaces of the wing, and $R=R_{tb}=10\cdot R_{lr}=0.5$ [m]  for the top and the bottom surfaces.
Factor $10$ is approximately the square root of the ratio, $A_{tb}/A_{lr}$,
between the sum, $A_{tb}$, of the areas of the top and of the bottom surfaces of the wing, and
the sum, $A_{lr}$, of the areas of the left and of the right sides. 
Finally, parameters $a$ and $b$ are set to $0.8$ and $1.3$, respectively.\\
We first employ the basic SIDW interpolation procedure, by skipping the enrichment step (ii) in \textbf{Algorithm 2}. 
Figure~\ref{wing_idw}, center highlights the new set $\widehat {\mathcal C}\subset {\mathcal C}$,
consisting of $92$ control points. Notice the very coarse and uniform distribution of control points $\widehat{\mathbf{c}}_k$.
The deformed configuration yielded by the SIDW approach  
essentially coincides with the one in Figure~\ref{wing}, bottom (see the analysis below for more quantitative investigations), despite a considerable reduction of the interpolation matrix, being $\widehat W\in \mathbb{R}^{797 \times 92}$, and of the corresponding computational time 
(the deformation phase takes only 0.015 [s]).\\
We successively resort to the ESIDW variant.
We preserve the same values for $a$, $b$, $R_{lr}$ and $R_{tb}$ as in the previous numerical check. 
Additionally, we constrain the selection procedure to include in  
$\widehat {\mathcal C}$ all the nodes along the left and the right edge profiles of the wing as well as
the nodes along the horizontal edges of the NACA profile. The rationale behind this choice is that control points acting on the left edge profile contribute with zero displacement in the interpolatory procedure. This results in a very small deformation for the neighboring internal points, thus preserving the continuity of the deformation, as the wing is clamped along that side. Moreover, since the displacement depends only on $z$, we expect that enforcing control points along horizontal edges (parallel to the $z$-axis) is instrumental to closely capture (with a small number of control points) the overall motion of the wing. Finally, control points are added along the right edge profile for symmetry.
The enrichment increases the cardinality of $\widehat {\mathcal C}$, now consisting of $390$ control points. As
shown in Figure~\ref{wing_idw}, bottom, all the nodes along the left and the right profiles of the wing are retained as control points, while the distribution of the control points on the top (and on the bottom) of the wing is essentially the same as in Figure~\ref{wing_idw}, center.
Also in this case the deformed configuration is essentially the same as the one in Figure~\ref{wing}, bottom (see the analysis below for more quantitative investigations).
Concerning the computational time, 
the increased cardinality of $\widehat {\mathcal C}$ results in an deformation stage that takes 0.034 [s]. This time, despite 
larger than the one demanded by SIDW, is still one order of magnitude less with respect to the time required by the standard IDW interpolation. 

\paragraph*{Sensitivity to $R.\ $}
We investigate the sensitivity of SIDW and ESIDW interpolation to
the selection radius. In particular, in Figure~\ref{wing_cpu_analysis}, we show the trend, as a function of $R$, of the 
number of the control points (left), and of the CPU time (in seconds) required to assemble the interpolation matrix
(center) and to compute the deformation of the wing (right), respectively\footnote{In order to improve the reliability of these timings, the same deformation has been repeated 100 times. Thus, assembly and deformation CPU times provided in the figures correspond to the average. The standard deviation is less than $0.001$ [s] in any case.}.  \\
The relation between the number of control points and $R$ is nonlinear for both the 
SIDW and ESIDW approaches. 
Nevertheless, while for the SIDW interpolation, the number of control points (and, consequently, 
the associated CPU time)  increases as $R$ becomes
smaller, a low sensitivity to $R$ is shown by the ESIDW variant. 
In particular, the cardinality of $\widehat {\mathcal C}$ remains the same for $R$ greater than $0.1$.\\
Concerning the CPU times, no significant difference distinguishes the trend associated with the matrix assembly
and the deformation step, also quantitatively. 
As expected, ESIDW interpolation is more computationally demanding than 
the basic SIDW approach, whereas the standard IDW approach coincides with the most expensive procedure.

\begin{figure*}
\centering
    \includegraphics[width=0.38\textwidth]{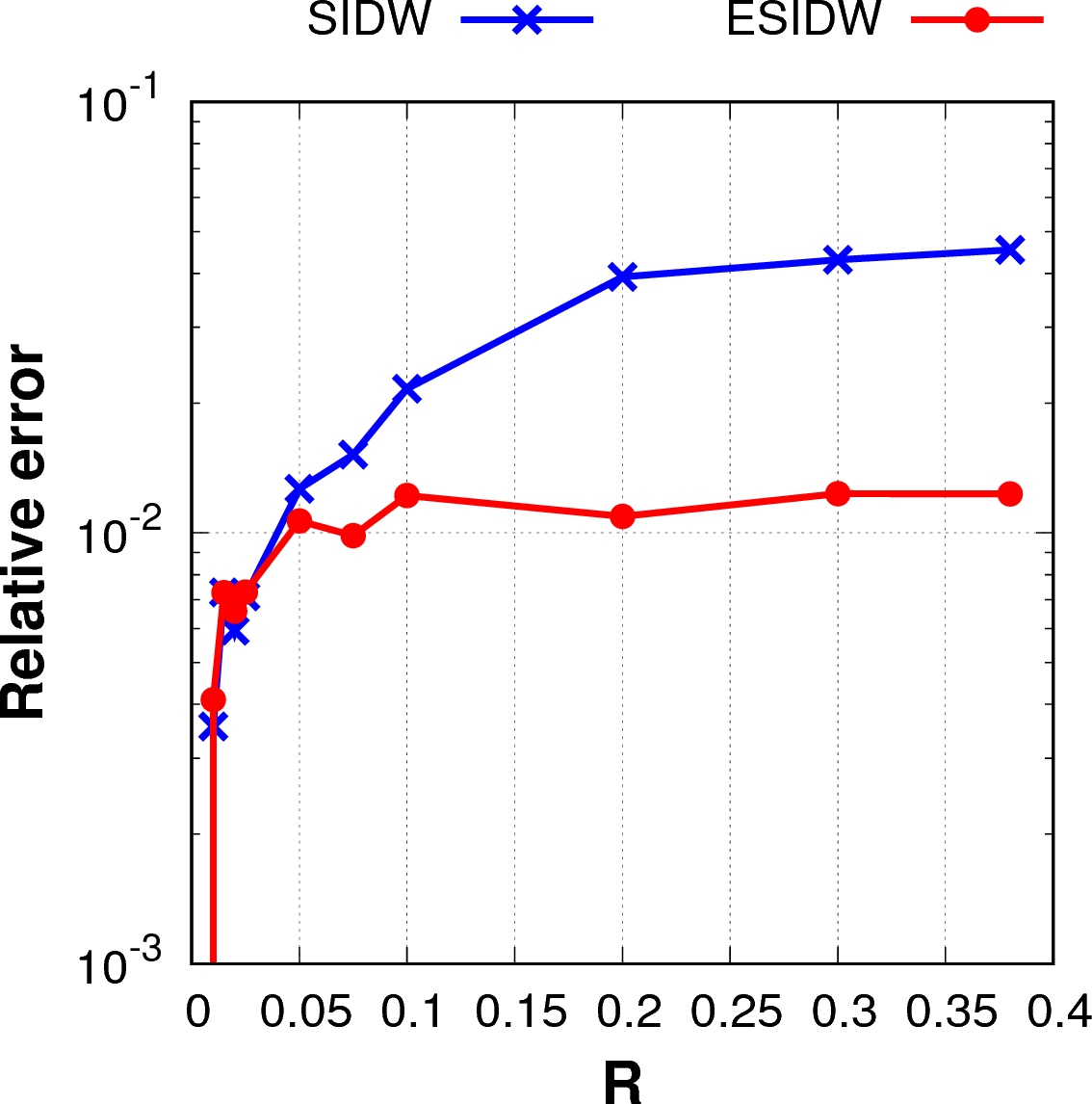}
    \includegraphics[width=0.58\textwidth]{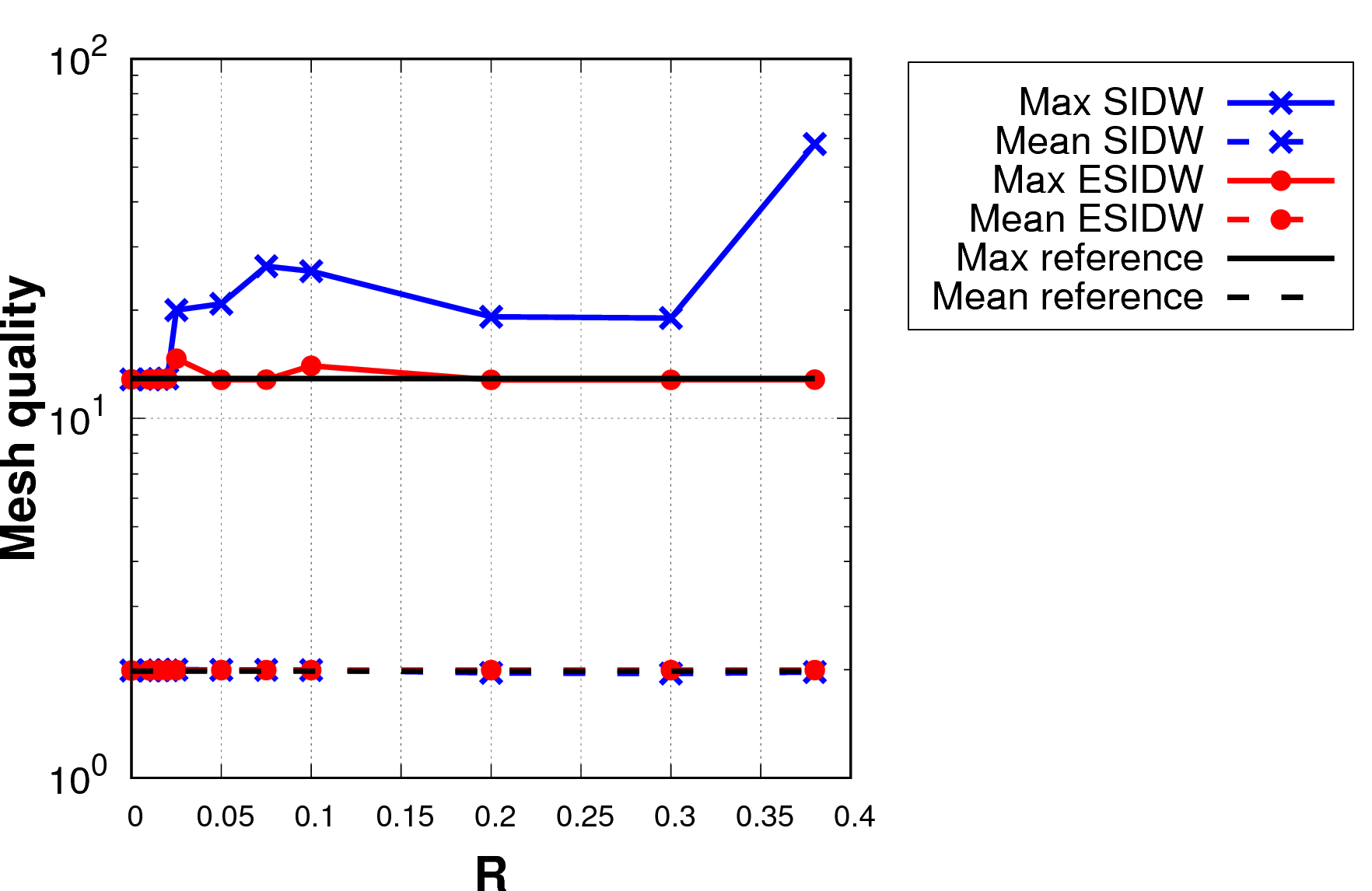}
\caption{Structural deformation of a wing: relative error (left) and maximum and mean mesh quality (right) as a function of $R$.}\label{wing_error_quality_analysis}
\end{figure*}

In Figure~\ref{wing_error_quality_analysis}, left we show the error trend as a function of the selection radius. 
We compute the $L^2(\Omega_h)$-norm of the relative error 
between the SIDW (ESIDW) and the IDW deformation. 
The larger number of control points employed by the ESIDW interpolation yields more accurate deformations
compared with the ones provided by the SIDW approach.
In particular, for large enough values of $R$, the error due to SIDW is approximately $4$\%, while, for the ESIDW
procedure, it is about $1$\%. Convergence to zero is guaranteed by both the methods 
as $R$ decreases.

Finally, in Figure~\ref{wing_error_quality_analysis}, right we investigate the influence of the selected radius $R$  
on the quality of the elements of the deformed meshes.
Different criteria can be employed to quantify the mesh quality $\mathcal Q$. Here, we adopt 
the ratio between the longest and the shortest edge \cite{FreyGeorge}.
In the figure, for different values of $R$, we compare the maximum and the mean value of $\mathcal Q$ on the meshes
generated by SIDW and ESIDW algorithms, with the corresponding values associated with the initial configuration $\Omega_h$. 
While, on average, the mean value of $\mathcal Q$ is essentially independent of 
the deformation procedure and of the selected $R$, more sensitivity is appreciable on the maximum value of $\mathcal Q$.
The sufficiently large number of control points allows ESIDW to preserve about the same value of $\mathcal Q$
as for the initial mesh, also when $R$ increases. On the contrary, a deterioration on the maximum 
mesh quality is evident when dealing with the SIDW interpolation algorithm, especially for large values of $R$. 

\begin{figure*}[h]
\centering
\includegraphics[width=0.35\textwidth]{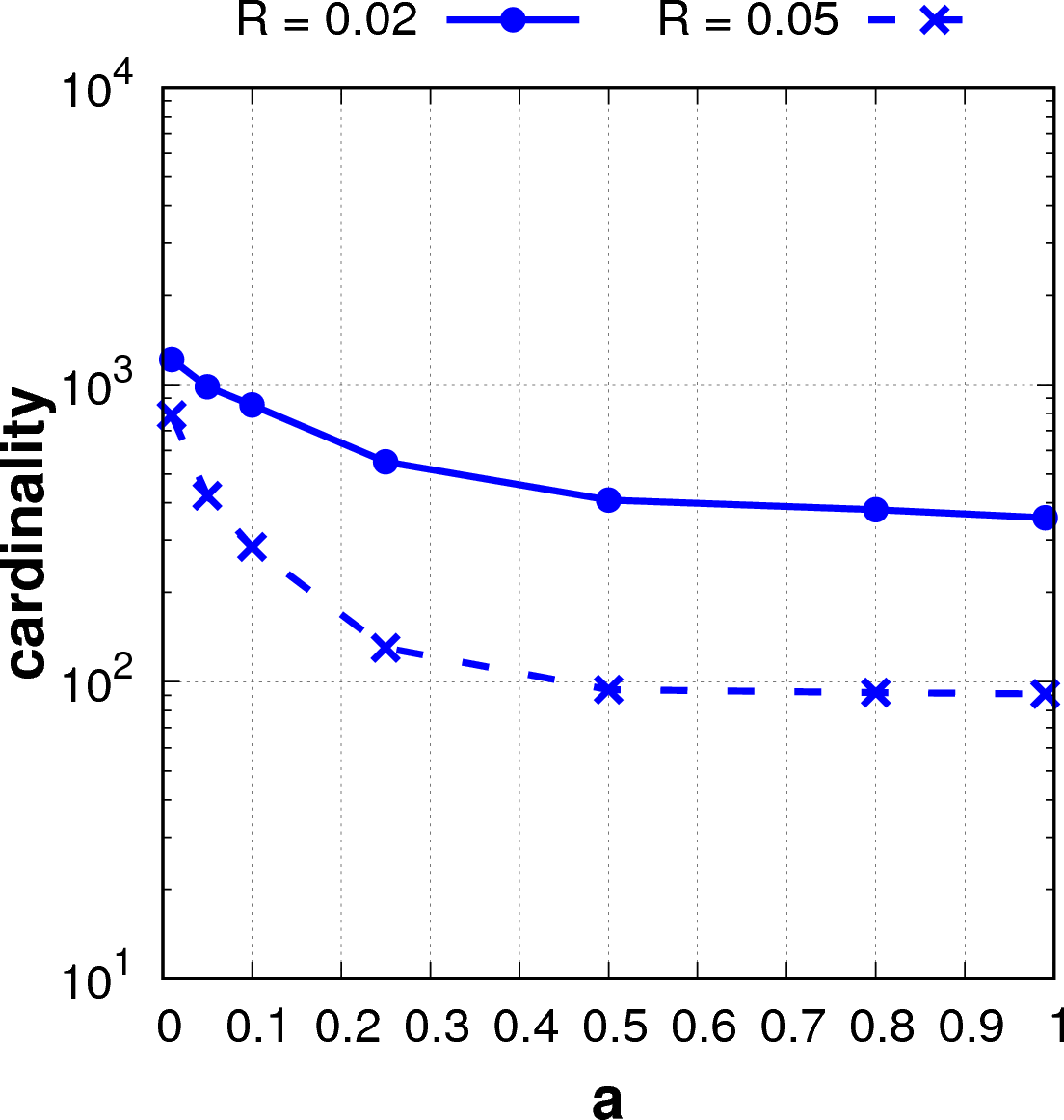}\quad
\includegraphics[width=0.35\textwidth]{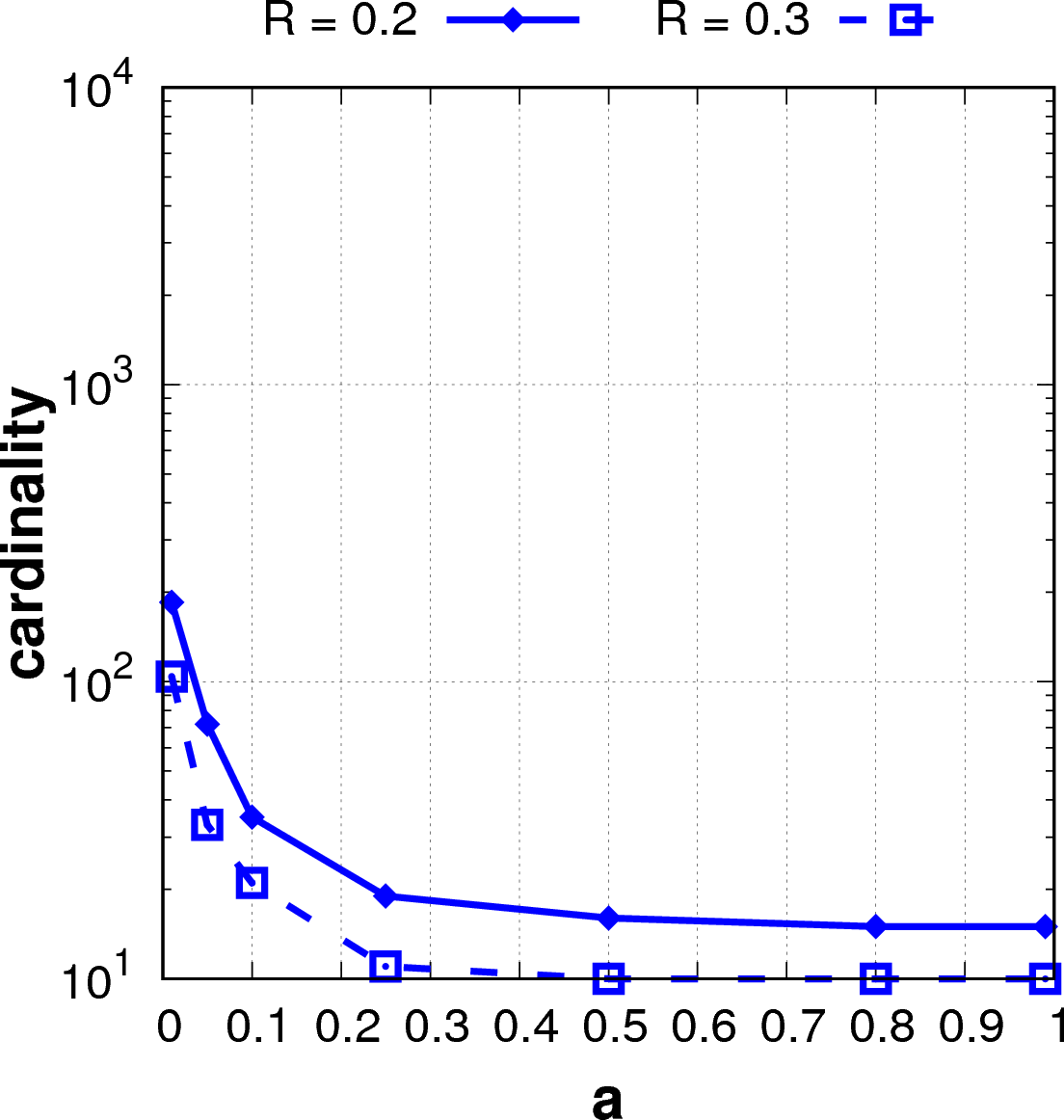}\\
\includegraphics[width=0.35\textwidth]{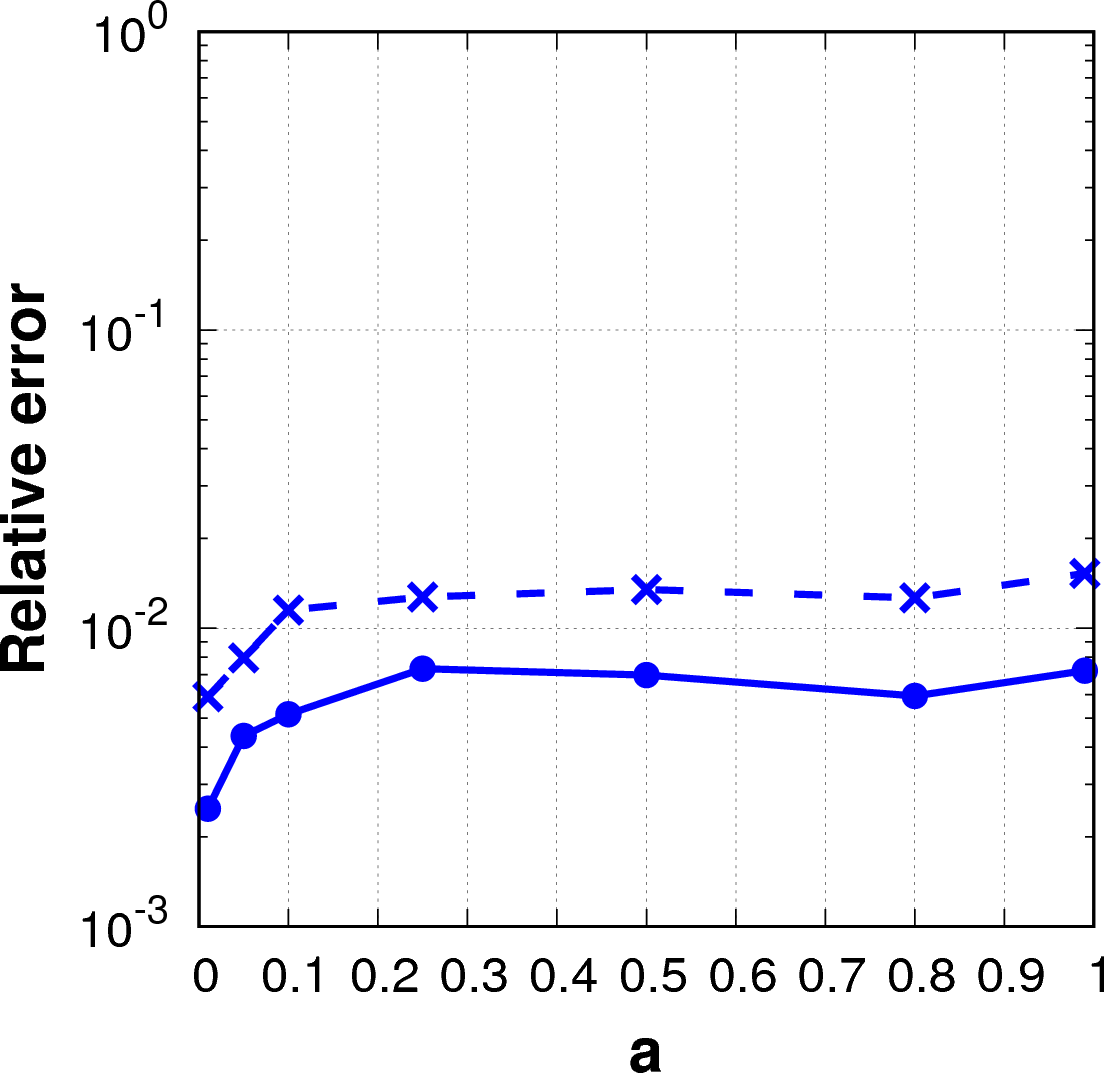}\quad
\includegraphics[width=0.35\textwidth]{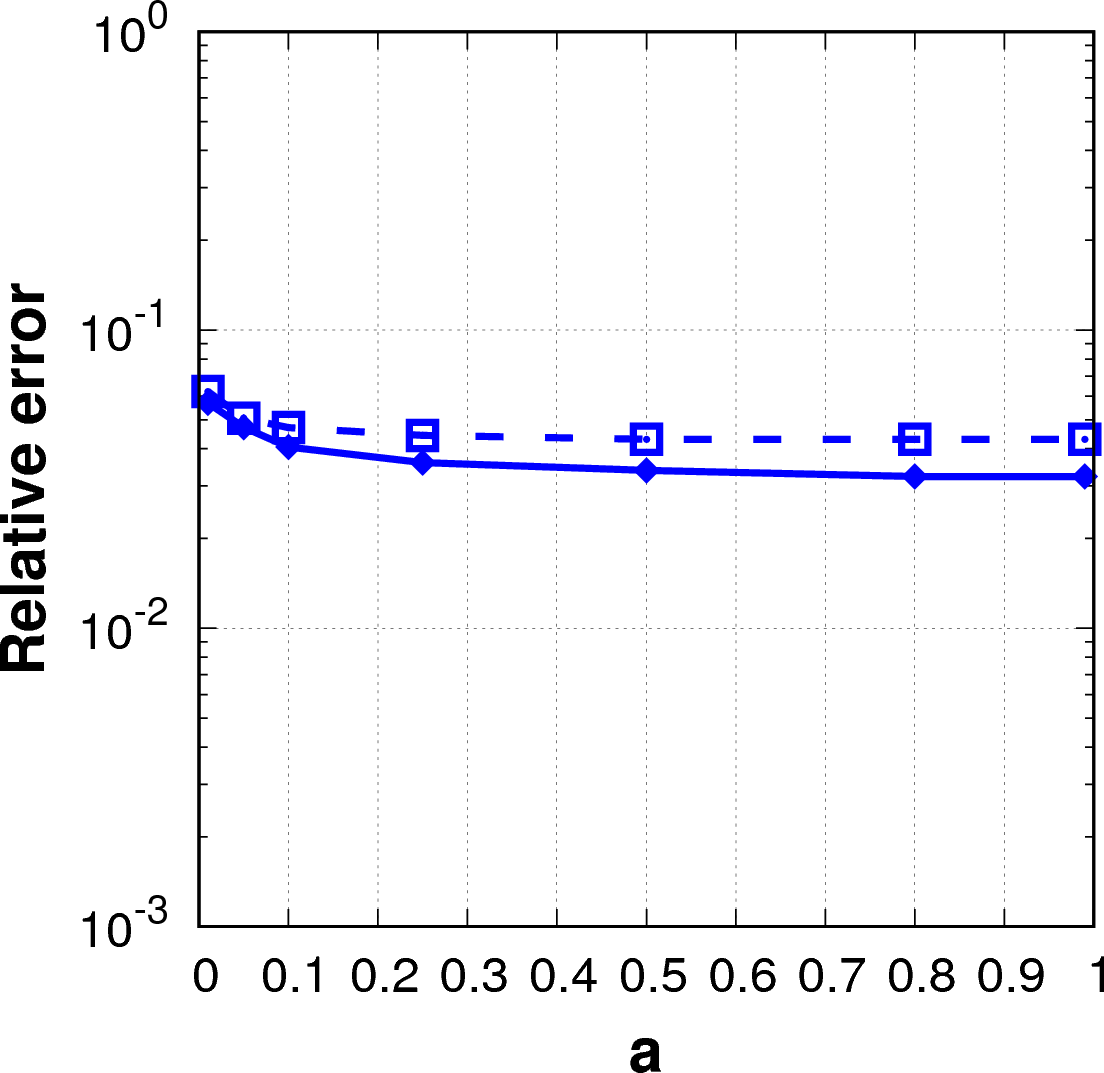}\\
\caption{Structural deformation of a wing (SIDW interpolation): cardinality of the control point set (top) and $L^2(\Omega_h)$-norm of the relative error (bottom) as a function of $a$, for $R = 0.02, 0.05$ (left) and for $R = 0.2, 0.3$ (right).}\label{sensitivity_to_ab}
\end{figure*}

\paragraph*{Sensitivity to $a$ and $b.\ $}
We also study the sensitivity of the control point cardinality and of the accuracy of SIDW interpolation
on the parameters $a$ and $b$. 
Similar conclusions can be drawn for the ESIDW approach.
To simplify this analysis, we relate $a$ and $b$ so that $b = 1/a$, we pick $a < 1$ (i.e., $b > 1$), and we
consider four different values of the selection radius, namely, $R=0.02$, $0.05$, $0.2$ and $0.3$. 
Figure~\ref{sensitivity_to_ab} collects the results of this check. 

For $a \ll 1$ (i.e., $b \gg 1$), the number of selected control points becomes larger and larger, as expected, with a consequent increment of the 
CPU time. Nevertheless, this does not necessarily entail an improvement in terms of accuracy, especially for large values of $R$.
Indeed, for $R = 0.2, 0.3$, choosing $a \approx 1$ (i.e., $b \approx 1$) 
reduces the number of control points of about one order of magnitude with respect to the choice $a \ll 1$ (i.e., $b \gg 1$), while only 
a slightly lower relative error is guaranteed (see Figure~\ref{sensitivity_to_ab}, right). On the contrary, when $R$ is small (see Figure~\ref{sensitivity_to_ab}, left), decreasing $a$ (i.e., increasing $b$) improves the accuracy. This is due the fact that, for these choices of $a$ ($b$) and $R$, the SIDW procedure selects almost all the available control points. 

Thus, since the actual goal of SIDW procedures, is to reduce the computational burden, we are essentially interested in sufficiently large values of $R$. This 
suggests us that the choice $a \approx 1$ (i.e., $b \approx 1$) ensures a reasonable trade-off between efficiency and accuracy to the morphing procedure. 

\paragraph*{Comparison with a random selection.}
We conclude this first test case by comparing the proposed selection procedure to a random control point selection, in terms of relative error. The random selection results from replacing step (i) in \textbf{Algorithm 2} with the following step:
\begin{enumerate}
\item[(i')] extract a random subset $\widehat {\mathcal C}_{\text{rand}} \subset \mathcal{C}$.
\end{enumerate}

\begin{figure*}[h]
\centering
\includegraphics[width=0.35\textwidth]{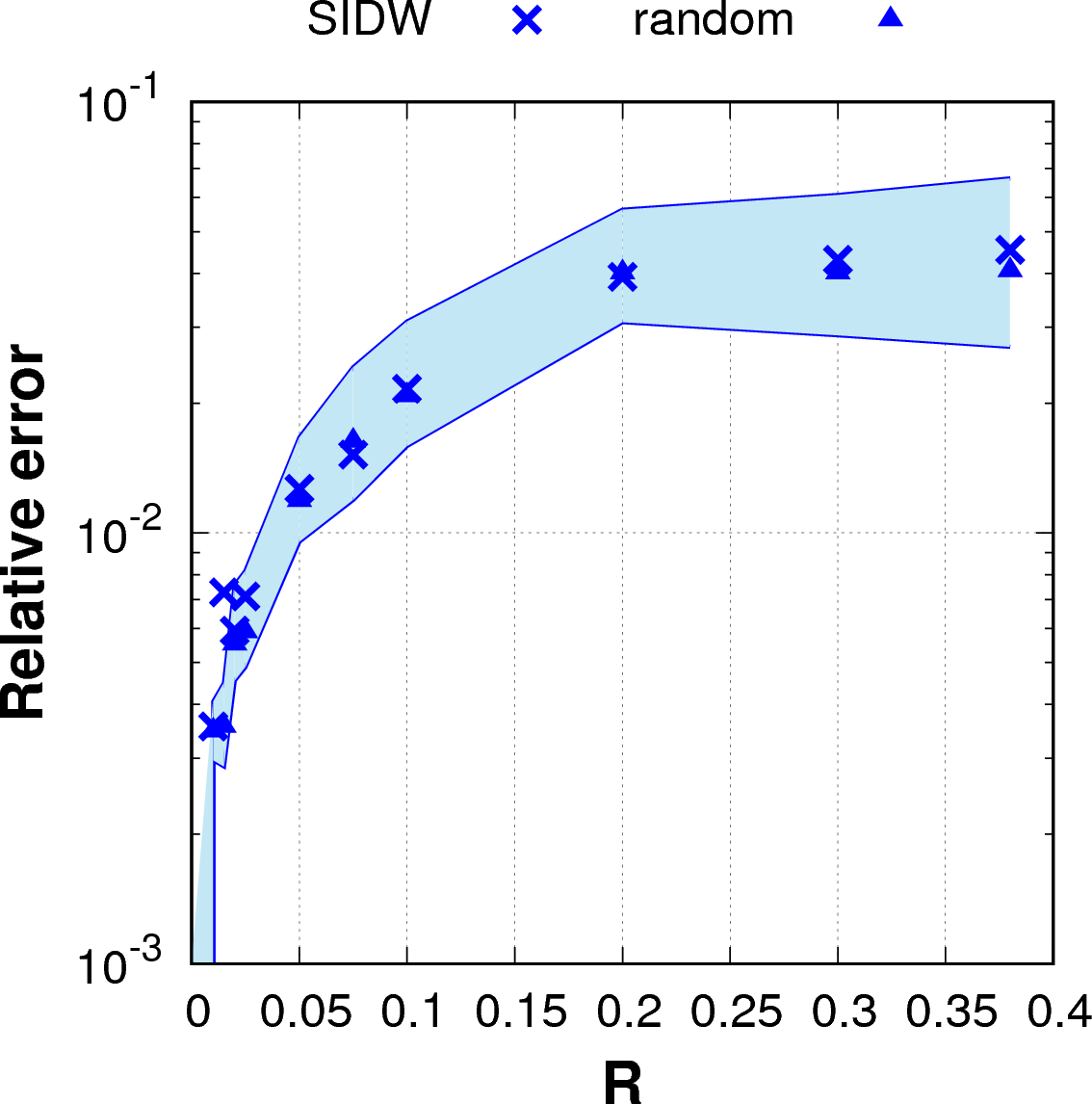}\quad
\includegraphics[width=0.35\textwidth]{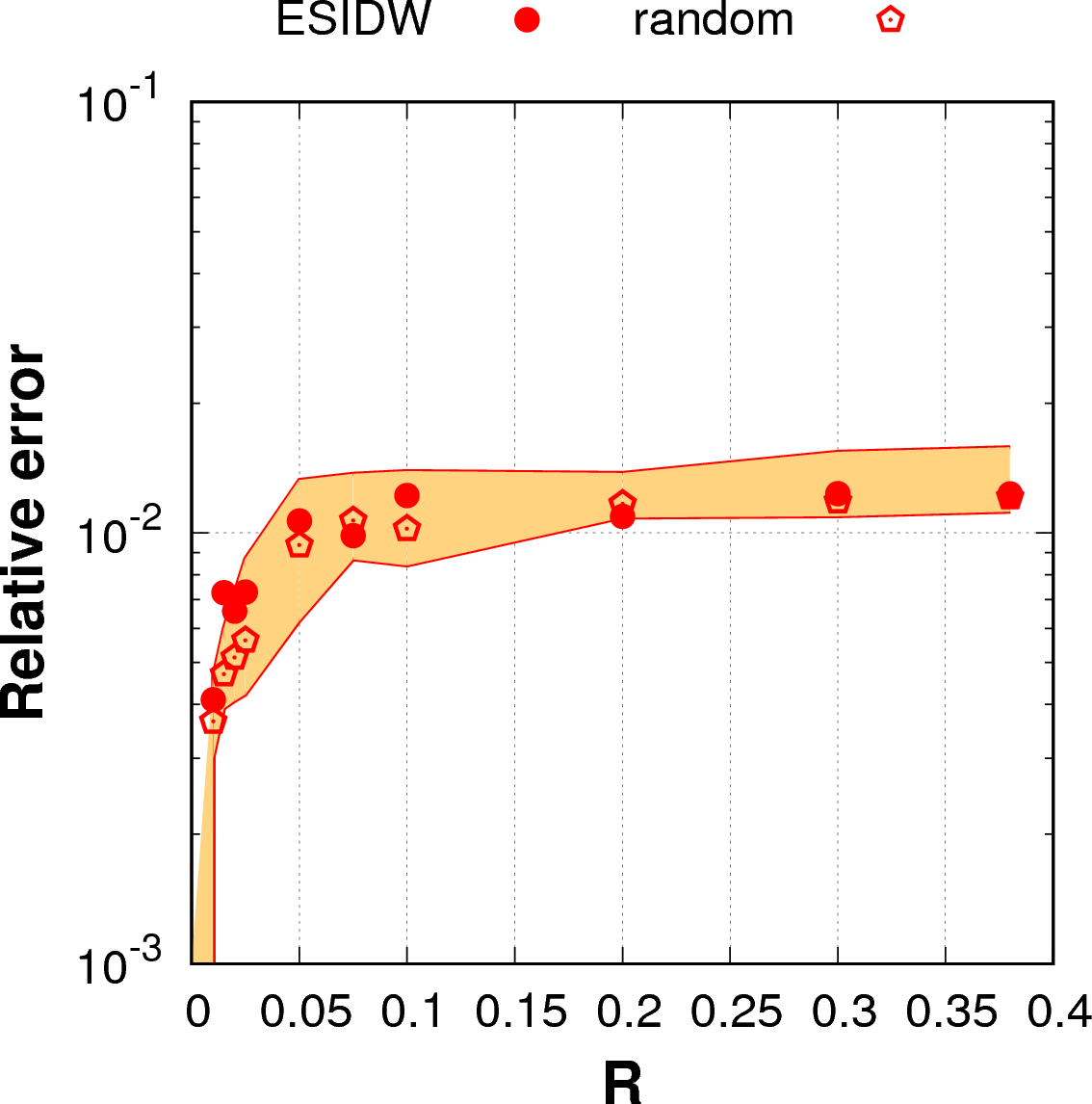}\\
\caption{Structural deformation of a wing: comparison of relative errors for SIDW (left) and ESIDW (right) selection procedures to a random control point selection.}\label{random}
\end{figure*}

The results of the comparison are shown in Figure~\ref{random}, for SIDW (left) and ESIDW (right) procedures. We compare (as a function of $R$) the relative error attained by the SIDW (ESIDW) selection with the average relative error over 100 random selections, each one such that $\text{card}(\widehat {\mathcal C}_{\text{rand}}) = \text{card}(\widehat {\mathcal C})$. In particular, the area between the minimum and maximum error over the random selection is shaded. The error due to SIDW (ESIDW) is very close to the average error, thus resulting in a numerically equivalent way to compute the average approximation that can be achieved with $\text{card}(\widehat {\mathcal C})$ control points. We also remark that, for large values of $R$ (i.e., for a very small number of control points), and especially for SIDW procedure, a random selection suffers from a large standard deviation $\sigma$ (see Table~\ref{tab_random}). A large number of random selections (and correspondingly of deformations) may thus be required: SIDW and ESIDW turn out to be consequently more efficient procedure, as they require only one selection (i.e., a unique deformation).
Finally, for what concerns the random selection, we denote by $\Delta_{\min}$ ($\Delta_{\max}$) the difference between the minimum (maximum) realized error and the average one, normalized to the standard deviation. As shown in Table~\ref{tab_random}, the resulting distribution is not symmetric, so that some random realizations return better configurations than the average one (to which the SIDW/ESIDW configuration is close to), since $|\Delta_{\min}| > 0$. However, they might also provide far worse configurations, being $\Delta_{\max} > |\Delta_{\min}|$. In particular, improvements are limited to the order of 1 standard deviation, while deteriorations can reach up to 3 standard deviations.

\begin{table*}
\centering
\begin{tabular}{|c|c|c|c|}
\hline
\multicolumn{4}{|c|}{\emph{Without enrichment}}\\\hline
$R$&$\sigma$&$\Delta_{\min}$&$\Delta_{\max}$\\\hline
0.01 & 0.00031 &-1.57 & 1.86\\\hline
0.02 & 0.00073 &-1.32 & 2.76\\\hline
0.05 & 0.0018 &-1.26 & 2.57\\\hline
0.1 & 0.0043 &-1.13 & 2.34\\\hline
0.2 & 0.0072 &-1.28 & 2.23\\\hline
0.3 & 0.0091 &-1.24 & 2.27\\\hline
\end{tabular}
\quad
\begin{tabular}{|c|c|c|c|}
\hline
\multicolumn{4}{|c|}{\emph{With enrichment}}\\\hline
$R$&$\sigma$&$\Delta_{\min}$&$\Delta_{\max}$\\\hline
0.01 & 0.00047 &-1.36 & 2.09\\\hline
0.02 & 0.00084 &-1.26 & 2.53\\\hline
0.05 & 0.0013 &-1.76 & 2.21\\\hline
0.1 & 0.0017 &-1.06 & 2.15\\\hline
0.2 & 0.0014 &-0.94 & 2.38\\\hline
0.3 & 0.0012 &-0.78 & 3.04\\\hline
\end{tabular}
\caption{Structural deformation of a wing: statistics for random selection procedure, without (left) and with (right) enrichment.}
\label{tab_random}
\end{table*}

\subsubsection{Fluid mesh motion around a wing}\label{sec:ass2}
\begin{figure*}
\centering
\includegraphics[width=0.48\textwidth]{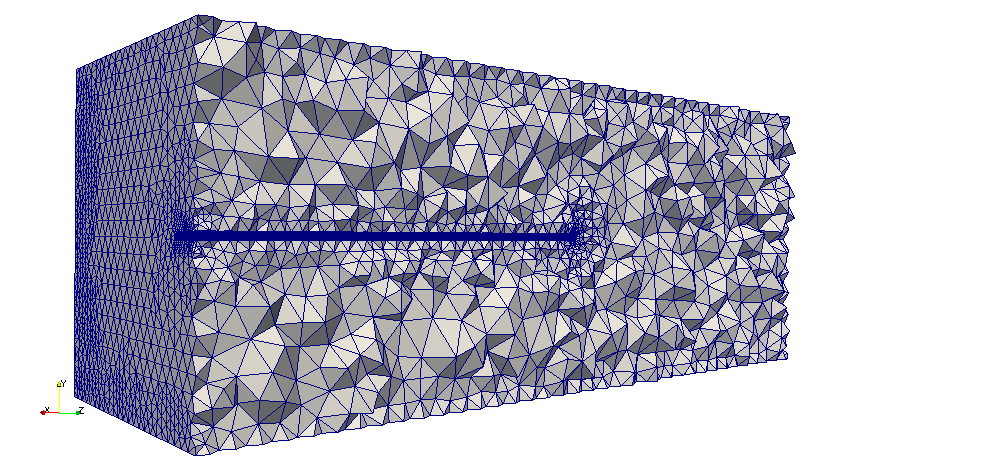}\;%
\includegraphics[width=0.48\textwidth]{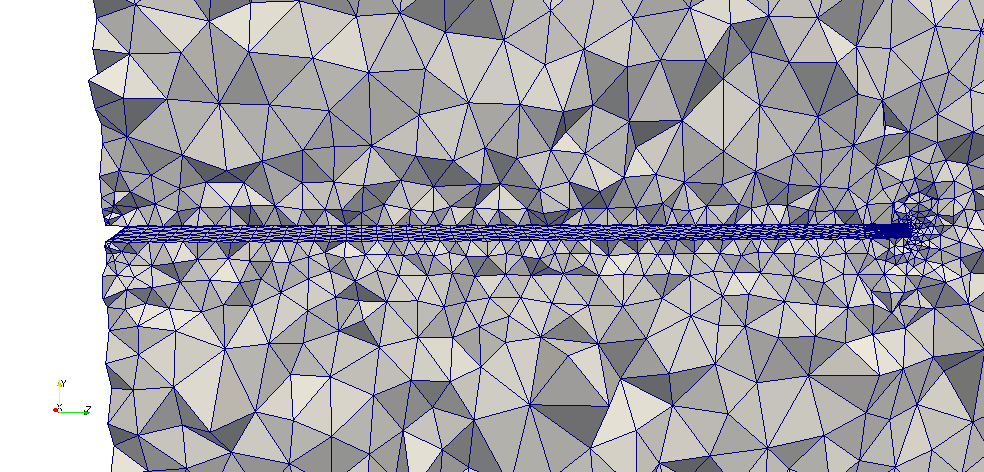}
\includegraphics[width=0.48\textwidth]{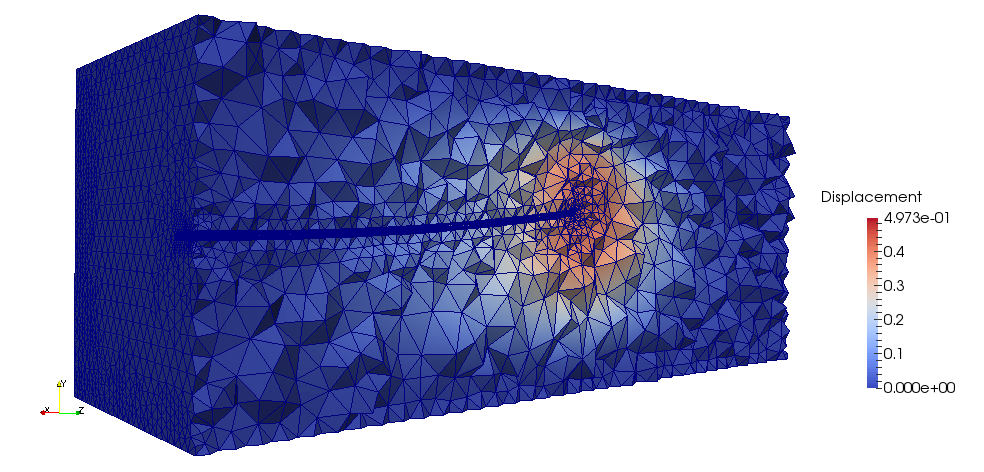}\;
\includegraphics[width=0.48\textwidth]{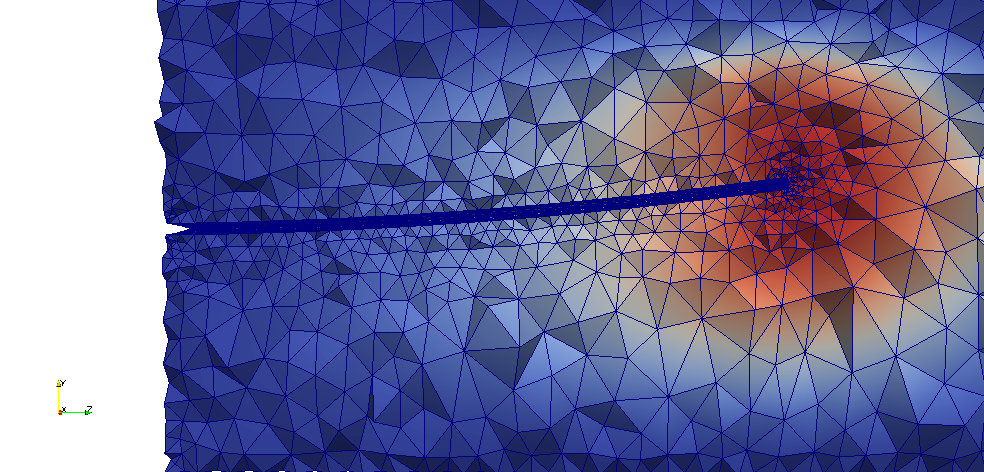}
\caption{Fluid mesh motion around a wing: rest (top) and deformed (bottom) configuration; 3D (left) and section (right) view.}\label{fluid_wing}
\end{figure*}

In this section we move to a FSI setting, by considering the fluid mesh motion around the NACA0012 profile in
Figure~\ref{fluid_wing}, top. 
This test case mimics a typical study performed in a wind tunnel, where the wing is clamped on one side while, on the other side, it is deformed by a vertical displacement, such as the one in \eqref{eq:wing_displacement}. Table~\ref{fluid_wing_table} gathers the main properties of the physical domain and of the corresponding discretization.

\begin{figure*}
\centering
    \includegraphics[width=0.6\textwidth]{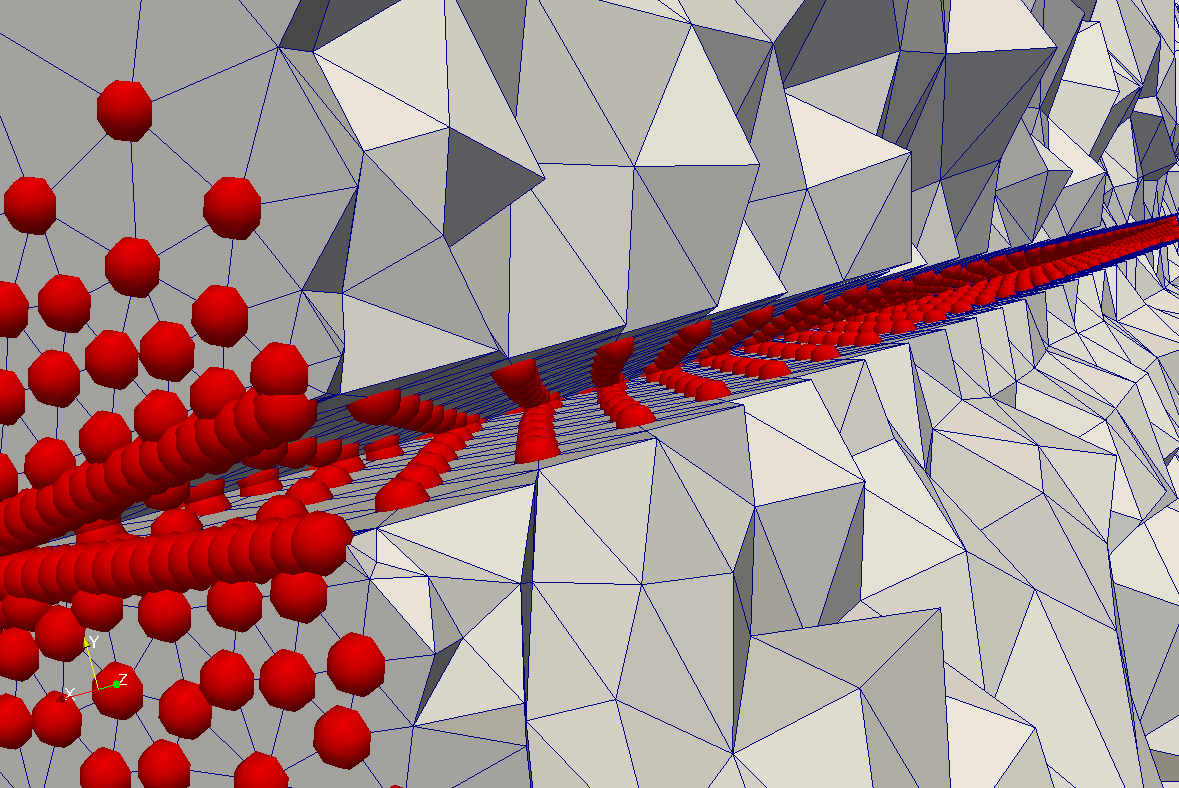}\\[2mm]
    \includegraphics[width=0.6\textwidth]{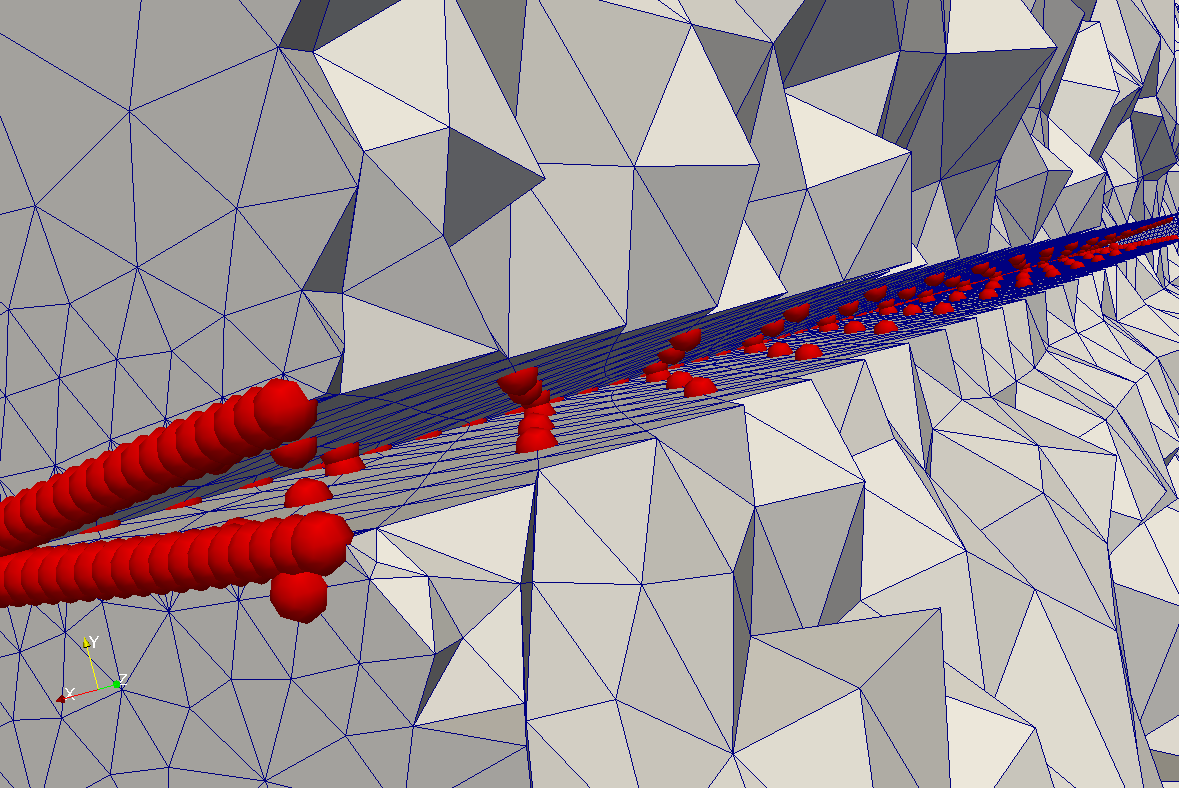}
\caption{Fluid mesh motion around a wing: control point distribution
associated with the IDW (top) and ESIDW (bottom) interpolation in correspondence with the clamped side.}\label{fluid_wing_idw}
\end{figure*}

To investigate the deformation in the fluid mesh, we start by applying the standard IDW approach, thus 
identifying the set $\mathcal{C}$ with all boundary nodes (see Figure~\ref{fluid_wing_idw}, top for a detail in correspondence with the clamped side). 
The corresponding interpolation matrix in \eqref{IDW_matrix_form} belongs to $\mathbb{R}^{21910 \times 14126}$. The resulting deformed configuration, shown in Figure~\ref{fluid_wing}, bottom, is obtained after $83.09$ [s].

In order to reduce the computational costs, we resort to the ESIDW interpolation, by enriching the selection of control points in 
$\widehat{\mathcal{C}}$ with all the nodes on the left and right edge profiles as well as along the horizontal edges of the NACA profile (we refer to the previous test case for a motivation to this choice).
Different values of radius $R$ are selected 
for the faces of the (outer) box and for the sides of the (inner) wing. 
In particular, with similar considerations as in the previous section, we choose:
$R_{top, b} = R_{bottom, b} = R_{front, b} = R_{rear, b} = 0.25$ [m] and $R_{right, b} = R_{left, b} = 0.125$ [m] for the box; 
$R_{top, w} = R_{bottom, w} = R_{right, w} = 0.025$ [m] for the wing.
Parameters $a$ and $b$ are set to 0.8 and 1.3, respectively. 
\textbf{Algorithm 2} provides a subset $\widehat{\mathcal{C}}$ consisting of $9339$ control points (see the enlarged view in Figure~\ref{fluid_wing_idw}, bottom), 
so that matrix $\widehat W$ in \eqref{SIDW_matrix_form} is now in $\mathbb{R}^{21910 \times 9339}$.
As shown in the figure, control points are essentially located on the structure profile only.\\
The selection process reduces the computational time to $57.07$ [s], and further reduction is possible taking a larger value of $R$ as discussed below.

Finally, the $L^2(\Omega_h)$-norm of the relative error between the ESIDW and the IDW deformation is approximately $5.86\%$.
\begin{table*}
\centering
\begin{tabular}{|cc|}
\hline
wind tunnel dimensions & $10 \times 5 \times 4\pi$ [m$^3$] \\
\hline
wing longitudinal dimension & $2\pi$ [m]	\\
\hline
wing chord length & $1.01$ [m]	\\
\hline\hline
\# elements & $169598$	\\
\hline
\# nodes & $36036$ \\
\hline  
\# internal nodes & $21910$ \\
\hline 
\# boundary nodes & $14126$ \\
\hline
\end{tabular}
\caption{Fluid mesh motion around a wing: main properties of $\Omega_r$ and of $\Omega_h$.} \label{fluid_wing_table}
\end{table*}

\paragraph*{Sensitivity to $R.\ $}
We numerically check the sensitivity of IDW and ESIDW techniques to the selection radius, by mimicking the investigation in the previous section. 
For this purpose, we choose $R_{top, b} = R_{bottom, b} = R_{front, b} = R_{rear, b} = R$, $R_{right, b} = R_{left, b} = 0.5 R$ for the box, and
$R_{top, w} = R_{bottom, w} = R_{right, w} = 0.1 R$ for the wing, while preserving the values previously adopted for $a$ and $b$. 

\begin{figure*}[h]
\centering
\includegraphics[width=0.31\textwidth]{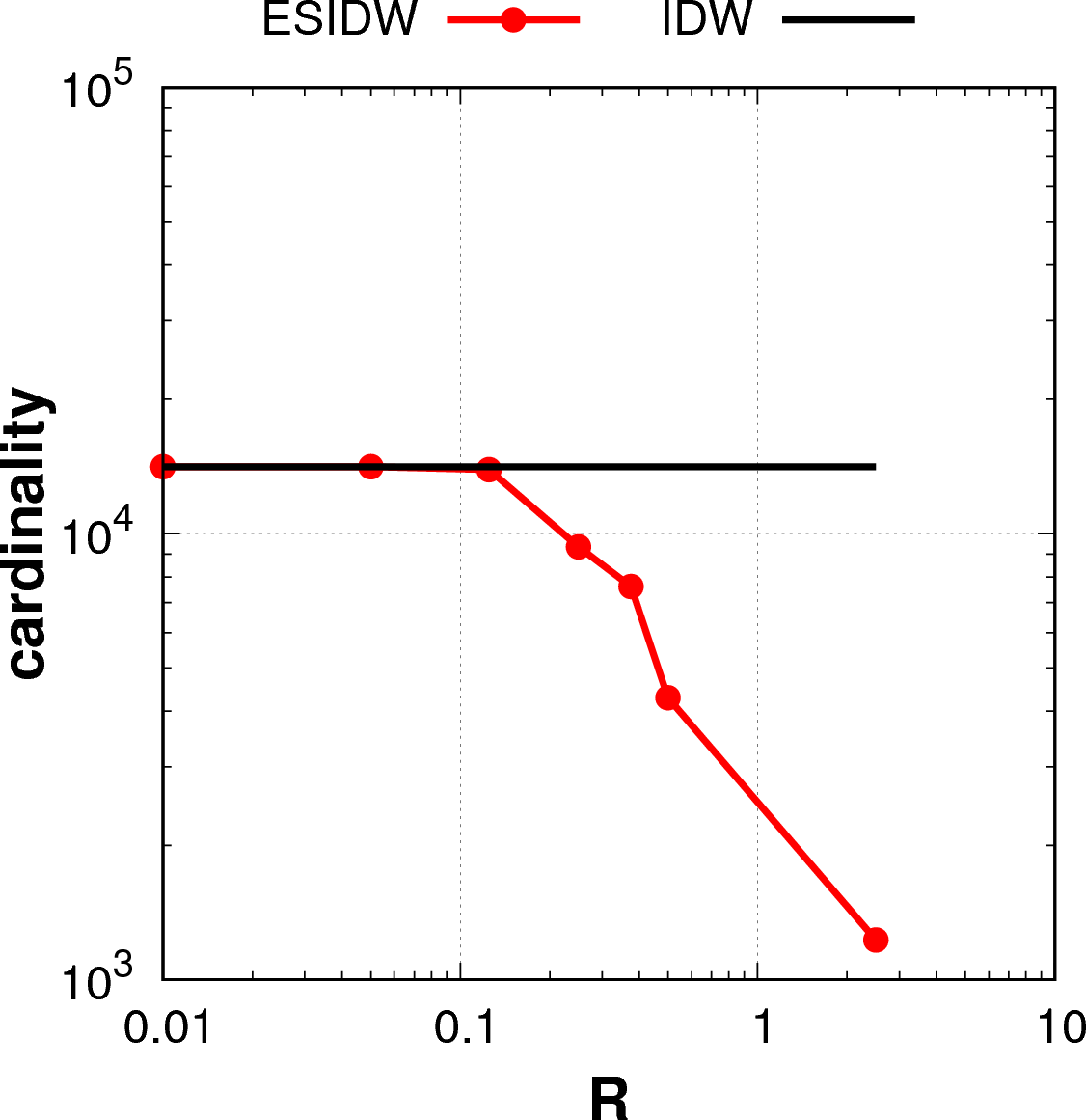}
\includegraphics[width=0.31\textwidth]{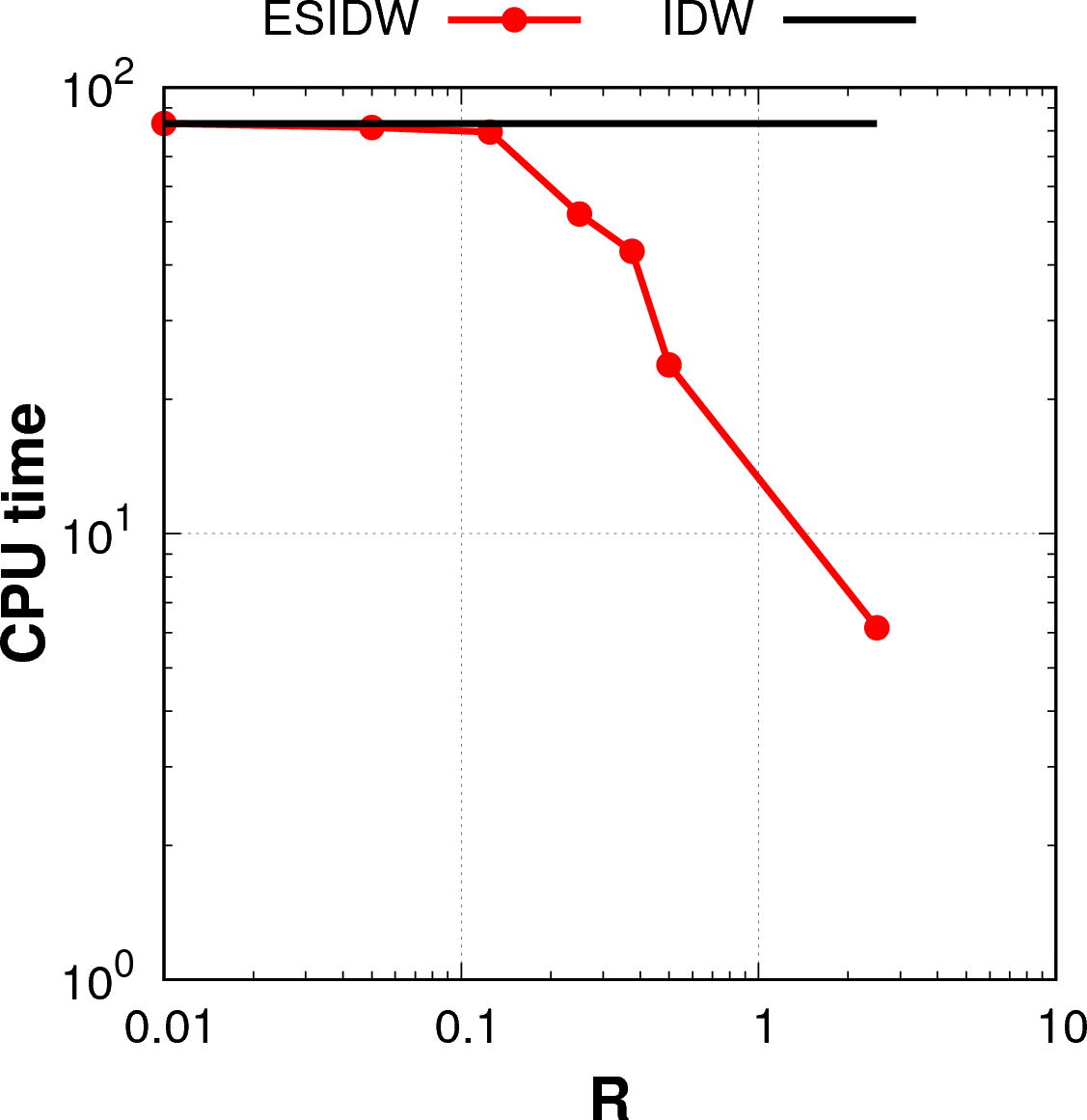}
\includegraphics[width=0.31\textwidth]{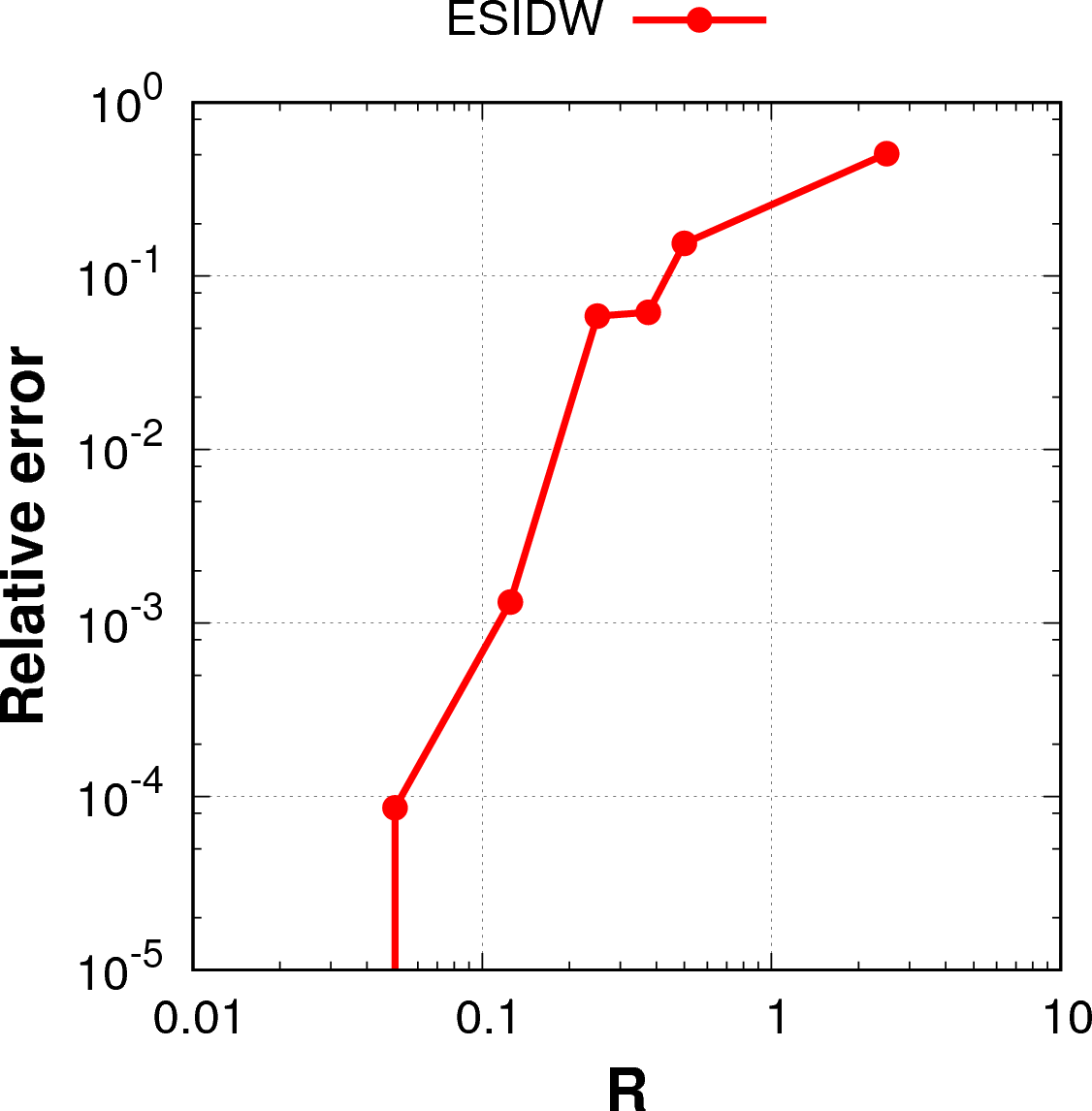}
\caption{Fluid mesh motion around a wing: cardinality of the control point set (left), total CPU time (center), and 
the relative error (right) as a function of R.}\label{sensi2}
\end{figure*}

The trend exhibited by the number of control points in Figure~\ref{sensi2}, left is completely different with respect to the one in Figure~\ref{wing_cpu_analysis}, left. The cardinality of $\widehat{\mathcal{C}}$ is essentially the same as for the IDW approach until 
$R$ sufficiently increases. Then, for $R>0.1$, the cardinality reduces more and more ensuring, for instance, a gain of one order of magnitude 
for $R>2$. \\
The central panel in Figure~\ref{sensi2} computes the associated CPU time (in seconds), where now we have quantified the assembly and the deformation time, altogether. As expected, the trend of the computational time follows the one of the control point cardinality, as well as 
the $L^2(\Omega_h)$-norm of the relative error between ESIDW and IDW deformation,  as shown in Figure~\ref{sensi2}, right.
In contrast to Figure~\ref{wing_error_quality_analysis}, left no stagnation of the error is detected when varying $R$.

Finally, concerning the quality characterizing the meshes yielded by IDW and ESIDW procedures, we have that 
the maximum and the mean value of $\mathcal Q$, independently of the adopted interpolation, is equal to $4.31$ and $1.48$, respectively 
with negligible variations (on the second decimal digit only) for ESIDW,  as long as $R \leq 1$. 

\subsubsection{Fluid mesh motion around a rotating hull}\label{sec:ass3}
\begin{figure*}
\centering
\includegraphics[width=0.9\textwidth]{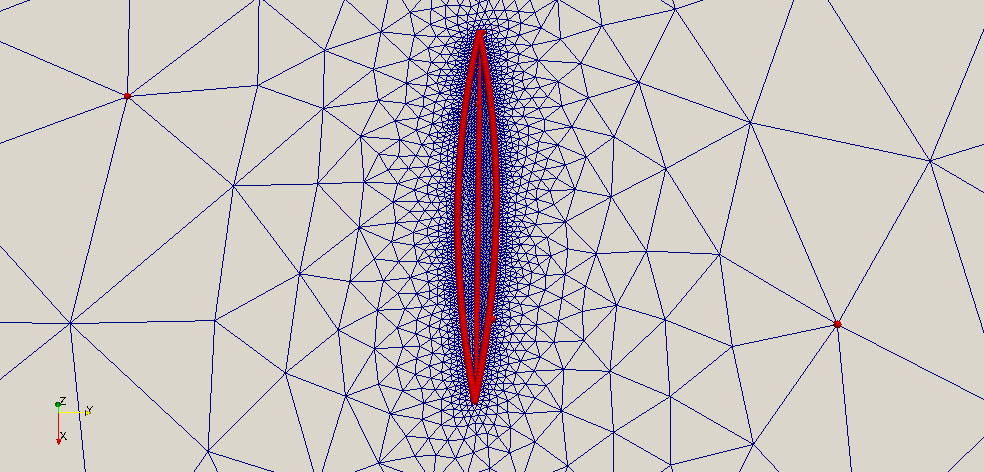}
\includegraphics[width=0.9\textwidth]{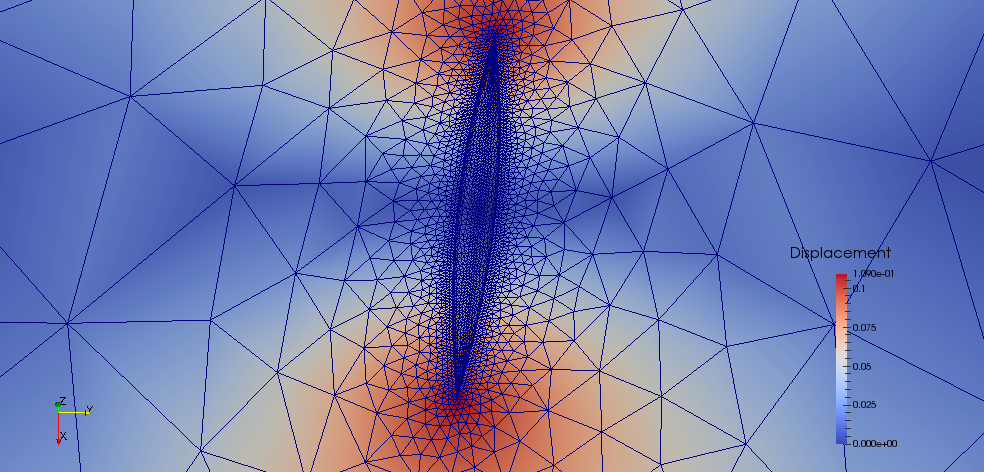}
\caption{Fluid mesh motion around a rotating hull: control points (top) and deformed configuration (bottom) provided by the ESIDW interpolation.}\label{fluid_hull}
\end{figure*}

We conclude this section by testing the performances of the ESIDW interpolation when dealing with a FSI setting
where a large deformation occurs. In more detail, we model a fluid mesh motion which results from the rotation of a structural domain in a naval engineering context. The initial configuration coincides with an outer box containing an inner 
Wigley hull \cite{Wigley1934} (see Figure~\ref{fluid_hull}, top). A rotation of $-5^\circ$ with respect to the $z$-axis is successively applied. We refer to 
Table~\ref{fluid_hull_table} for a summary of the main properties of the physical domain and of the corresponding discretization.
\begin{table*}
\centering
\begin{tabular}{|cc|}
\hline
mesh dimension & $50 \times 25 \times 10$ [m$^3$] \\
\hline
hull length & $2.5$ [m]	\\
\hline\hline
\# elements & $30265$	\\
\hline
\# nodes & $7186$ \\
\hline  
\# internal nodes & $3322$ \\
\hline 
\# boundary nodes & $3864$ \\
\hline
\end{tabular}
\caption{Fluid mesh motion around a rotating hull: main properties of $\Omega_r$ and of $\Omega_h$.} \label{fluid_hull_table}
\end{table*}
\begin{figure*}
\includegraphics[width=0.9\textwidth]{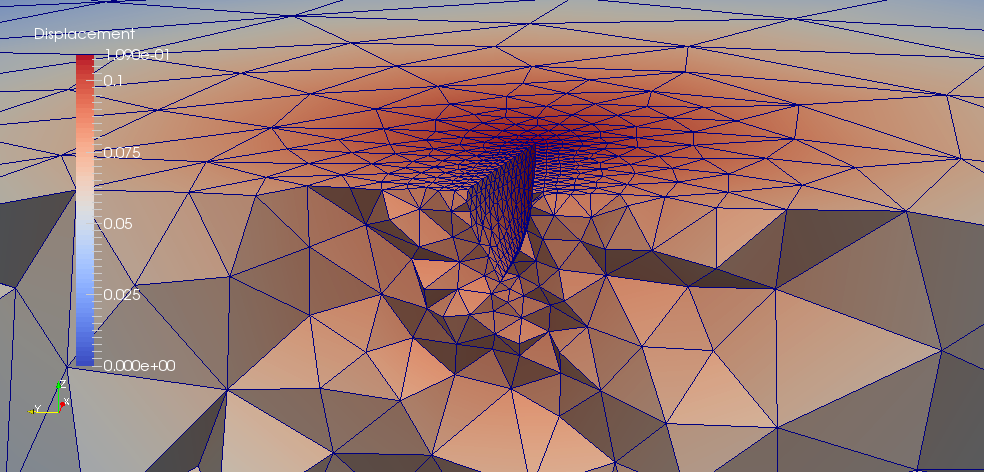}
\caption{Fluid mesh motion around a rotating hull: detail of the deformed configuration provided by ESIDW interpolation.}\label{fluid_hull_2}
\end{figure*}

The standard IDW approach results in solving the linear system \eqref{IDW_matrix_form}, matrix $W$ being in 
$\mathbb{R}^{3322 \times 3864}$. This leads to a CPU time of $4.51$ [s].

To reduce the computational burden, we apply \textbf{Algorithm 2}. 
The selection is performed by considering a radius $R = 1.5$ [m], $a=0.8$ and
$b=1.3$. Notice that the value of $R$ is comparable with the hull length. 
This allows a considerable reduction in the number of control points along the boundaries of the fluid domain. 
Moreover, to ensure that interior mesh nodes accurately capture the rotation of the hull, we switch on the enriching step {\tt (ii)} of the selection procedure,
by adding to $\widehat{\mathcal{C}}$ all the nodes belonging to the edge of the hull. 
The resulting distribution of control points is shown in Figure~\ref{fluid_hull}, top. 
Most of the points are identified by the constraint step, while very few nodes (only two in the specific case) are retained by the selection at step {\tt (i)}. As a result, the hull turns out to be sharply described by the control points.
ESIDW procedure essentially halves the computational time, the CPU time of the deformation phase being now equal to $2.56$ [s].
Figure~\ref{fluid_hull}, bottom displays the deformed fluid mesh, obtained by solving system \eqref{SIDW_matrix_form}, with 
$\widehat W \in \mathbb{R}^{3322 \times 1057}$. We also refer to Figure~\ref{fluid_hull_2} for a detail of the deformed grid.

The computational cost reduction provided by ESIDW does not compromise the accuracy of the approximation. Indeed, 
the $L^2(\Omega_h)$-norm of the relative error between the ESIDW and the IDW deformation is $2.42\%$.

\begin{figure*}[h]
\centering
\includegraphics[width=0.31\textwidth]{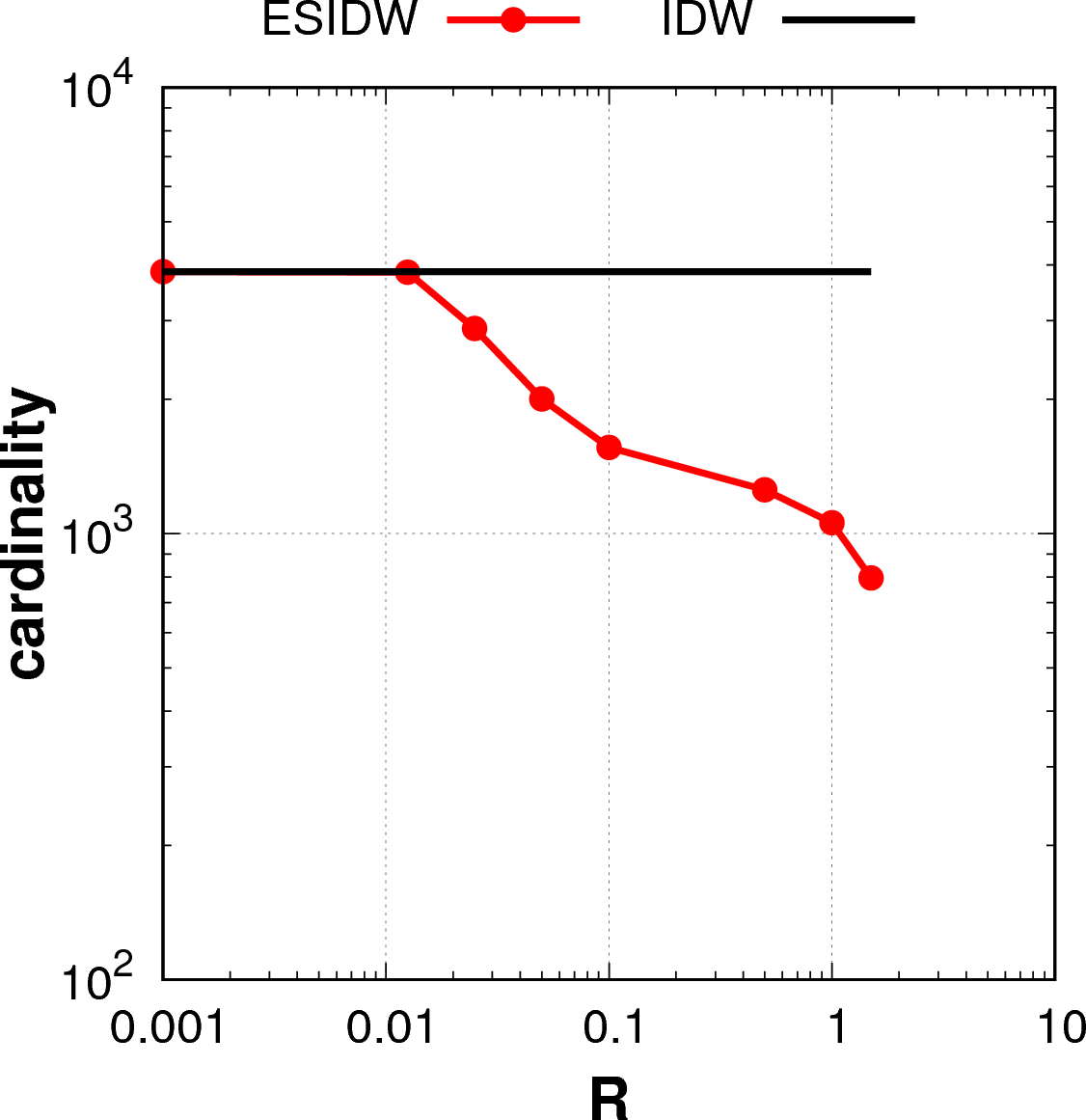}
\includegraphics[width=0.31\textwidth]{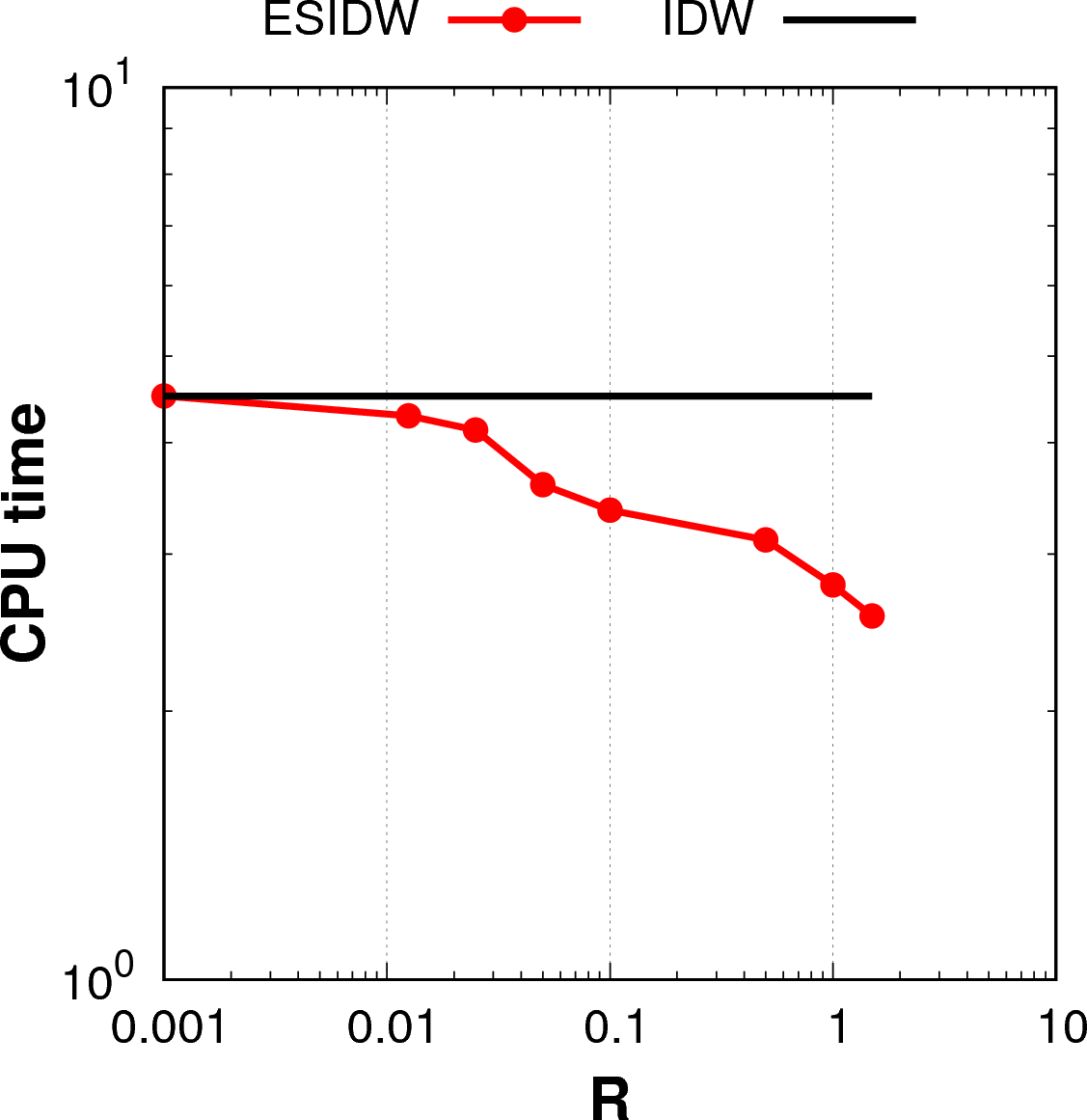}
\includegraphics[width=0.31\textwidth]{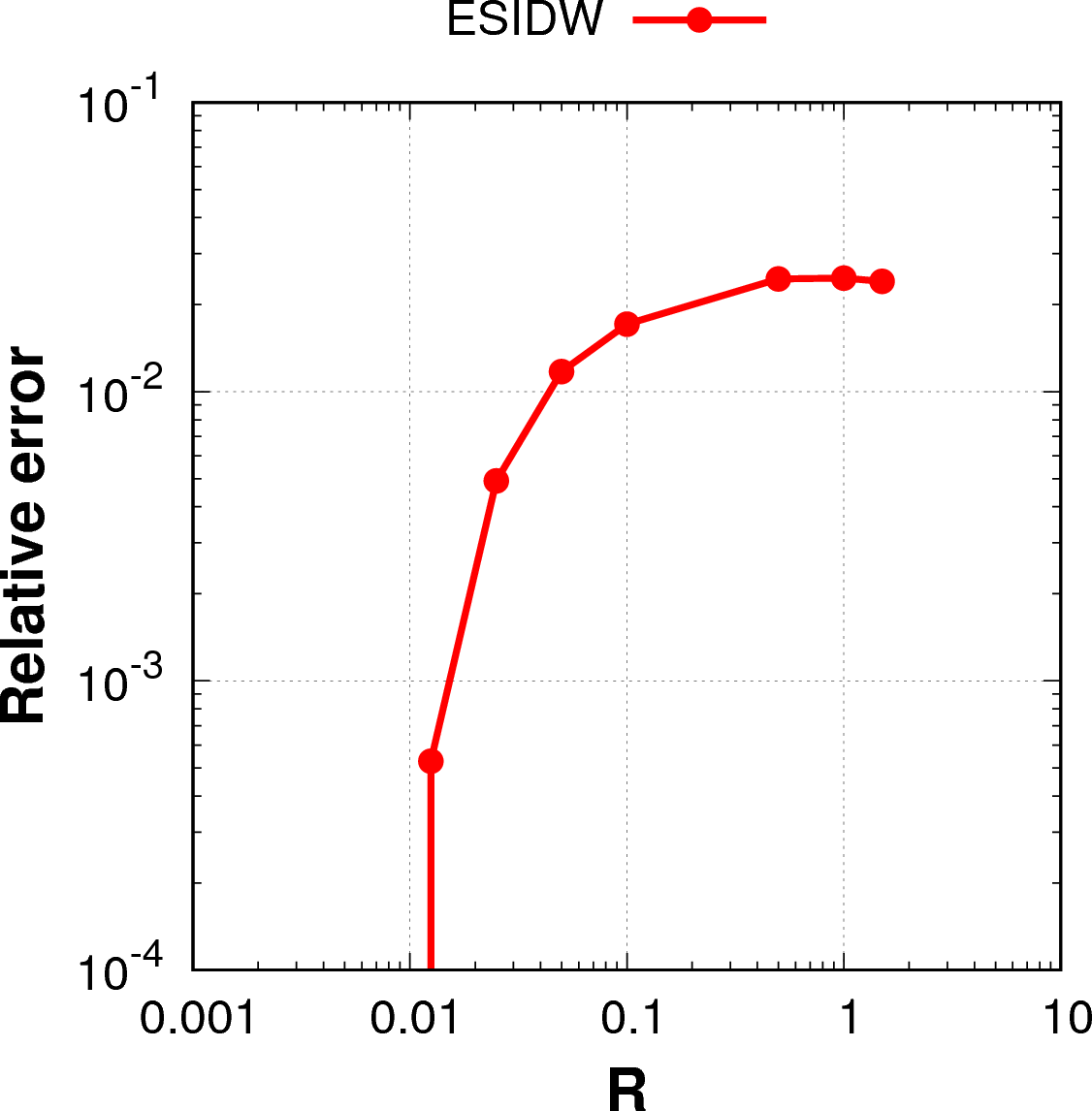}
\caption{Fluid mesh motion around a rotating hull: cardinality of the control point set (left), total CPU time (center), and 
the relative error (right) as a function of R.}\label{sensi3}
\end{figure*}

We have investigated the sensitivity to $R$ also for this configuration, by mimicking the analysis in the previous section. 
Figure~\ref{sensi3} provides the same results as in Figure~\ref{sensi2}. In contrast to this last case, we remark a less significant
reduction in the number of the control points and of the total CPU time. Moreover, while for the wing configuration
the ESIDW interpolation allows to gain six orders of magnitude in terms of accuracy for a sufficiently small $R$, only 2 orders less
are obtained for the rotating hull configuration. \\
Finally, the rotation causes a slight increment in the maximum mesh quality (which changes from $2.47$ in the reference configuration
to $2.56$ for the deformed one), while essentially preserving the mean value of $\mathcal Q$ (varying from $1.62$ to $1.63$ before
and after the deformation). These values are essentially independent of the adopted interpolation and of $R$.

\section{POD for SIDW shape morphing}\label{sec:model_reduction}
We merge now the SIDW techniques presented in the previous section with a dimensionality reduction technique, following an offline-online paradigm
\cite{hesthavenrozzastamm2016}. 
The reduced basis will be determined via 
a \emph{Proper Orthogonal Decomposition} (POD) technique (see, e.g., \cite{hesthavenrozzastamm2016,burkardtgunzburgerlee2006,volkwein2012proper,benner2013survey}). 
POD reduces the dimensionality of a system by transforming the original variables into a new set of uncorrelated variables (called POD modes, or principal components), so that, ideally, the first few modes retain most of the `energy' of the original system.

Generalizing notation, let $\boldsymbol{\mu} \in {\mathbb D}$ 
denote the generic parameter identifying the displacement $\mathbf{d}_{\widehat c}(\boldsymbol{\mu})$
of the control points $\widehat{\mathbf{c}}_k$ 
in $\widehat{\mathcal C}$.
As in the test cases of the previous section, we assume $\mathbf{d}_{\widehat c}(\boldsymbol{\mu})$ to be an input for the shape morphing process. More in general, the map $\mathbf{d}_{\widehat c}(\boldsymbol{\mu})$ can be analytically provided (as in our synthetic test cases) or can be the output of a solver. For instance, when dealing with FSI problems, $\mathbf{d}_{\widehat c}(\boldsymbol{\mu})$ coincides with the evaluation of the structural displacements (computed by solving the elastodynamic equations) on the selected control points $\widehat{\mathcal C}$ along the fluid-structure interface.

According to \eqref{SIDW_matrix_form}, the deformation 
of the internal nodes $\mathbf{x}_j$, for $j=1, \ldots, \mathcal{N}_h^i$, 
of the discretized reference domain $\Omega_h$ can be computed as 
\begin{equation}\label{POD_displace}
\mathbf{d}^i(\boldsymbol{\mu})= {\widehat W} \mathbf{d}_{\widehat c}(\boldsymbol{\mu}),
\end{equation}
the dependence on the parameter being highlighted. 
Notice that, while both the input displacement $\mathbf{d}_{\widehat c}(\boldsymbol{\mu})$ and the output 
deformation $\mathbf{d}^i(\boldsymbol{\mu})$ depend on $\boldsymbol{\mu}$, the SIDW matrix $\widehat W$ is parameter independent. 

During the POD offline stage, we consider a training set, $\Xi_{\text{train}} = \{\boldsymbol{\mu}_n\}_{n=1}^{n_{\text{train}}}$, of $n_{\text{train}}$ values for the parameter $\boldsymbol{\mu}$.
With each parameter $\boldsymbol{\mu}_n$, we associate a certain displacement
$\mathbf{d}_{\widehat c}(\boldsymbol{\mu}_n)$ of $\widehat{\mathcal C}$, and we compute the corresponding deformation $\mathbf{d}^i(\boldsymbol{\mu}_n)$ in \eqref{POD_displace}. This allows us to assemble the
so-called \emph{snapshot} matrix
\begin{equation}\label{matrixU}
\mathcal{U} = [\mathbf{d}^i(\boldsymbol{\mu}_1), \cdots , \mathbf{d}^i(\boldsymbol{\mu}_{n_{\text{train}}})]
\in \mathbb{R}^{ \mathcal{N}_h^i \times n_{\text{train}}},
\end{equation}
whose columns store the different scenarios induced by the training set of deformations.
To extract the desired reduced basis, we perform now the Singular Value Decomposition (SVD) of 
matrix $\mathcal{U}$, so that
\begin{equation}\label{eq:snapshots_matrix}
 \mathcal{U} = \mathcal{W} \Sigma \mathcal{V}^T,
\end{equation}
where 
$$
\mathcal{W} = [\boldsymbol{\zeta}_1, \boldsymbol{\zeta}_2, \cdots, \boldsymbol{\zeta}_{\mathcal{N}_h^i}] \in \mathbb{R}^{\mathcal{N}_h^i \times \mathcal{N}_h^i}
$$
and
$$
\mathcal{V} = [\boldsymbol{\Psi}_1, \boldsymbol{\Psi}_2, \cdots, \boldsymbol{\Psi}_{n_{\text{train}}}] \in \mathbb{R}^{n_{\text{train}} \times n_{\text{train}}}
$$ 
are the orthogonal matrices of the \emph{left} and of the \emph{right singular vectors} of $\mathcal{U}$, while 
$
\Sigma \in \mathbb{R}^{\mathcal{N}_h^i \times n_{\text{train}}}
$ is a rectangular matrix of the same size as $\mathcal{U}$, whose only non-zero entries are the positive real diagonal entries,
$\sigma _1 \geq \sigma _2 \geq \cdots \geq \sigma_r > 0$, which represent the \emph{singular} (or \emph{principal}) 
values of $\mathcal{U}$, $r \leq n_{\text{train}}$ being the rank of $\mathcal{U}$~\cite{Golub}.\\
The reduced \emph{POD basis} is thus identified by a subset of the left singular vectors of $\mathcal{U}$ collected in the 
\emph{POD basis matrix} 
\begin{equation}\label{eq:basis_matrix}
\mathcal{Z} = [\boldsymbol{\zeta}_1 , \cdots, \boldsymbol{\zeta}_N] \in \mathbb{R}^{\mathcal{N}_h^i \times N}.
\end{equation}
with $N \leq r$. The integer $N$ can be selected via different criteria. 
In particular, since the energy contained in the discarded (i.e., the last $r - N$) POD modes is provided by
\begin{equation}\label{nrg}
E_\mathcal{U}(N) = \frac{\sum_{n=N+1}^{r} \sigma_n^2}{\sum_{i=1}^{r} \sigma_n^2},
\end{equation}
we set a tolerance $\varepsilon$ and choose $N$ as the first integer such that  $E_\mathcal{U}(N) \leq \varepsilon$.
The identification of the POD basis concludes the offline phase. 

The online phase starts from the choice of a new value $\boldsymbol{\mu}^* \in \mathbb{D}$ for the parameter, i.e., with the assignment of a new displacement $\mathbf{d}_{\widehat c}(\boldsymbol{\mu}^*)$
to the control points in $\widehat{\mathcal C}$.
The POD basis is then exploited to approximate the corresponding deformation $\mathbf{d}^i(\boldsymbol{\mu}^*)$ of the internal nodes in $\Omega_h$, defined as
\begin{equation}\label{eq:rev2}
\mathbf{d}^i(\boldsymbol{\mu}^*) = {\widehat W} \mathbf{d}_{\widehat c}(\boldsymbol{\mu}^*).
\end{equation}

In particular, we look for 
a suitable linear combination of the POD basis functions $\{\boldsymbol{\zeta}_l\}_{l=1}^{N}$, such that
\begin{equation}\label{eq:ls2}
\mathcal{Z}\; \boldsymbol{\beta}(\boldsymbol{\mu}^*) = {\widehat W} \mathbf{d}_{\widehat c}(\boldsymbol{\mu}^*)
\end{equation}
and where $\boldsymbol{\beta}(\boldsymbol{\mu}^*) \in \mathbb{R}^N$, with $[\boldsymbol{\beta}(\boldsymbol{\mu}^*)]_l=\beta_l$ and $l=1, \ldots, N$, is a vector of unknown coefficients depending on parameter
$\boldsymbol{\mu}^*$. 
We manipulate \eqref{eq:ls2} as
\begin{equation}\label{per_le_normali}
{\widehat W}^+ \mathcal{Z}\; \boldsymbol{\beta}(\boldsymbol{\mu}^*) = \mathbf{d}_{\widehat c}(\boldsymbol{\mu}^*)
\end{equation}
where ${\widehat W}^+\in \mathbb{R}^{\mathcal{N}_{\widehat c} \times \mathcal{N}_h^i}$ denotes the pseudo-inverse of $\widehat W$~\cite{Golub}.
The vector $\boldsymbol{\beta}(\boldsymbol{\mu}^*)$ is then computed in the least squares sense, i.e., by solving the normal equation system
\begin{equation}\label{eq:ls3}
\mathcal{Z}^T ({\widehat W}^+)^T \, {\widehat W}^+ \mathcal{Z}\; \boldsymbol{\beta}(\boldsymbol{\mu}^*) = 
\mathcal{Z}^T ({\widehat W}^+)^T  \mathbf{d}_{\widehat c}(\boldsymbol{\mu}^*),
\end{equation}
\reviewerB{assuming that ${\widehat W}^+ \mathcal{Z}$ has full column rank}.

The system we are lead to solve has size $N$, being $\mathcal{Z}^T ({\widehat W}^+)^T \, {\widehat W}^+ \mathcal{Z}\in  \mathbb{R}^{N \times N}$ and $\mathcal{Z}^T ({\widehat W}^+)^T\in \mathbb{R}^{N \times \mathcal{N}_{\widehat c}}$.
The deformation of the domain is finally computed via the product 
in \eqref{eq:rev2}.

\begin{remark}
\reviewerB{As alternative to system \eqref{eq:ls3}, one can derive from \eqref{eq:ls2}
\begin{equation}
\boldsymbol{\beta}(\boldsymbol{\mu}^*) = \mathcal{Z}^+ {\widehat W} \mathbf{d}_{\widehat c}(\boldsymbol{\mu}^*),
\label{eq:ls4}
\end{equation}
i.e., solve the normal equation system
\begin{equation*}
\mathcal{Z}^T \mathcal{Z}\; \boldsymbol{\beta}(\boldsymbol{\mu}^*) = 
\mathcal{Z}^T {\widehat W}  \mathbf{d}_{\widehat c}(\boldsymbol{\mu}^*),
\end{equation*}
assuming $\mathcal{Z}$ to have full column rank.

A comparison between this approach and the one proposed above is beyond the scope of this work. Nevertheless, we remark that the term $({\widehat W}^+)^T \, {\widehat W}^+$ in \eqref{eq:ls3} can be conceived as a weighting factor which takes into account the mesh connectivity.

Indeed, \eqref{eq:ls4} results from the solution of the minimization problem
\begin{equation*}
\mathcal{Z}\boldsymbol{\beta}(\boldsymbol{\mu}^*) = \argmin_{\mathbf{x} \in \text{range}(\mathcal{Z})} \|\mathbf{x} - {\widehat W} \mathbf{d}_{\widehat c}(\boldsymbol{\mu}^*)\|_2,
\end{equation*}
while \eqref{eq:ls3} stems from the minimization problem
\begin{equation*}
\mathcal{Z}\boldsymbol{\beta}(\boldsymbol{\mu}^*) = \argmin_{\mathbf{x} \in \text{range}(\mathcal{Z})} \|\mathbf{x} - {\widehat W} \mathbf{d}_{\widehat c}(\boldsymbol{\mu}^*)\|_w,
\end{equation*}
where the weighted norm $\|\cdot\|_w$ is defined as
\begin{equation*}
\|\mathbf{x}\|_w = \sqrt{\mathbf{x}^T ({\widehat W}^+)^T \, {\widehat W}^+ \mathbf{x}}.
\end{equation*}

We note that the withed approach in \eqref{eq:ls3} entails stronger assumptions compared with \eqref{eq:ls4}. In particular, problem \eqref{eq:ls4} only requires $\mathcal{Z}$ to be full rank, while \eqref{eq:ls3} necessitates ${\widehat W}^+ \mathcal{Z}$ to be full rank. As a result, condition $\mathcal{N}_{\widehat c} \geq N$ has to be added for the weighted approach, while such a condition is not required for the other approach. In practice, we stress that condition $\mathcal{N}_{\widehat c} \geq N$ is trivially verified since we are interested in performing a dimensionality reduction.
}
\end{remark}

We observe that both the matrices in \eqref{eq:ls3} can be computed once and for all at the end of the offline stage.
Thus, the online phase reduces to evaluate the matrix-vector product on the right-hand side of \eqref{eq:ls3},
to solve a linear system of size $N$ and then to compute the linear combination $\mathcal{Z}\; \boldsymbol{\beta}(\boldsymbol{\mu}^*)$. In terms of computational burden, the most effort involves
the internal nodes. In more detail, we are comparing the full problem \eqref{eq:rev2}
demanding $O(\mathcal{N}_h^i \cdot \mathcal{N}_{\widehat c})$ floating point operations with the POD approach
characterized by $O(N^3 + N \cdot \mathcal{N}_{\widehat c}+\mathcal{N}_h^i \cdot N)$ operations. This suggests that, 
if $N\ll \mathcal{N}_{\widehat c}$, we expect a computational saving via the reduced approach. This will be numerically verified in the next section.

The complete POD-SIDW interpolation procedure is itemized in \textbf{Algorithm~\ref{alg:pod_sidw_shape_morphing}}
when applied to a generic shape morphing context. 
The POD-ESIDW variant can be set in a straightforward way
simply by switching on the optional step {\tt (ii)} of \textbf{Algorithm~2} at the first item
of the offline phase.
With a view to a FSI problem, the same setting as for 
\textbf{Algorithm~2} is adopted.
\begin{algorithm*}
\normalsize
\caption{POD-SIDW interpolation for shape morphing}\label{alg:pod_sidw_shape_morphing}
\begin{algorithmic}[999]
\State \textsc{Offline phase}: 
    \State \qquad (a) apply the \textsc{Selection phase} of {\bf Algorithm~\ref{alg:sidw_shape_morphing}} to extract 
    $\widehat{\mathcal C}$ and to assemble the SIDW matrix $\widehat W$;
    \State \qquad (b) \textbf{for each ${\boldsymbol \mu}_n \in \Xi_{\text{train}}$}
    \State \qquad\qquad apply the \textsc{Deformation phase} of {\bf Algorithm~\ref{alg:sidw_shape_morphing}} to compute    
    $\mathbf{d}^i(\boldsymbol{\mu}_n)$ via \eqref{POD_displace};
    \State \qquad \qquad \textbf{end for}
    \State \qquad (c) assemble the snapshot matrix $\mathcal{U}$ in \eqref{matrixU};
    \State \qquad (d) extract the POD basis matrix $\mathcal{Z}$ in \eqref{eq:basis_matrix};
    \State \qquad (e) compute matrices $\mathcal{Z}^T ({\widehat W}^+)^T \, {\widehat W}^+ \mathcal{Z}$ and 
    $\mathcal{Z}^T ({\widehat W}^+)^T$ in \eqref{eq:ls3}.\\

\State \textsc{Online phase}: chosen $\boldsymbol{\mu}^* \in \mathbb{D}$:
    \State \qquad (f) solve system \eqref{eq:ls3} to derive the weights in $\boldsymbol{\beta}(\boldsymbol{\mu}^*)$;
    \State \qquad (g) compute the linear combination $\mathcal{Z}\; \boldsymbol{\beta}(\boldsymbol{\mu}^*)$.
\end{algorithmic}
\end{algorithm*}

\subsection{POD-SIDW algorithms in action}\label{combined_Sec_num}
We come back to the test cases in Section~\ref{sec:app} to investigate possible benefits on SIDW interpolation techniques due to POD. In particular, we quantify the computational improvements in terms of CPU time and accuracy of the approximation.

\subsubsection{Structural deformation of a wing}\label{sec:ass1_POD}
We apply the POD reduction procedure to a parameter dependent variant of the configuration in Section~\ref{sec:ass1}. To this aim, we impose the parametrized vertical displacement  
\begin{equation}\label{eq:par_ver_disp}
\delta y=\delta y(z; \mu) = \mu z^2
\end{equation}
to the wing, where the parameter $\mu$ is a scalar varying in the interval $\mathbb{D} = [0, 1.3]$. 

We focus on the standard IDW and on the ESIDW interpolation techniques. With reference to the ESIDW approach,
we preserve the two different choices done for the selection radius in 
Section~\ref{sec:ass1}, by picking 
$R=R_{lr}=0.05$ [m] and $R=R_{tb}=0.5$ [m] for the left and right and for the top and bottom surfaces of the wing, respectively.
Then, accordingly to Section~\ref{sec:ass1}, the enrichment is performed by adding in 
$\widehat{\mathcal C}$ all the nodes along the left, the right and the 
horizontal edges of the NACA profile.

To build the snapshot matrix $\mathcal{U}$, we randomly select $n_{\text{train}} = 100$ values in ${\mathbb D}$
which identify the set $\Xi_{\text{train}}$. Then, to extract the POD basis
we fix $\varepsilon = 10^{-5}$ as tolerance on the energy $E_{\mathcal{U}}(N)$ in \eqref{nrg}. Finally, we choose $\mu^* = 0.65$ as parameter value during the online phase. 
\begin{table*}
\centering
\begin{tabular}{|c|c|c|c|c|c|}
\hline
\textbf{method}&{\rm card}($\widehat{\mathcal C}$)&\textbf{CPU time}(\textsc{Offline})&\textbf{{\rm card}(${\mathcal Z}$)}&\textbf{CPU time}(\textsc{Online})&\textbf{Error}\\ 
\hline
IDW & 1666 & - & - & 0.25 [s] & - \\
\hline
ESIDW & 390 & - & - & 0.034 [s] & 1.06 \% \\
\hline
POD-IDW & 1666 & 90.93 [s] & 1 & 0.013 [s] & negligible \\
\hline
POD-ESIDW & 390 & 76.84 [s] & 1 & 0.013 [s] & 1.10 \% \\
\hline
\end{tabular}
\caption{Structural deformation of a wing: comparison between the basic IDW and ESIDW techniques and the corresponding POD variants.} \label{summarizing_costs_table_ass1}
\end{table*}

Table \ref{summarizing_costs_table_ass1} compares the performances of the plain IDW and ESIDW interpolations
with the corresponding POD variants. 
The second column provides the number of control points employed for morphing the original structure, being understood that 
$\widehat{\mathcal C}\equiv {\mathcal C}$ when dealing with the IDW approach.
The third column gathers the CPU time required by the offline phase of \textbf{Algorithm~\ref{alg:pod_sidw_shape_morphing}} to construct 
the response matrix, to extract the POD basis, whose cardinality is furnished in the fourth column, and to assemble matrices
$\mathcal{Z}^T ({\widehat W}^+)^T \, {\widehat W}^+ \mathcal{Z}$ and $\mathcal{Z}^T ({\widehat W}^+)^T$ in \eqref{eq:ls3}.
 The fifth column summarizes the CPU time required to perform the shape morphing via \eqref{IDW_matrix_form} and \eqref{SIDW_matrix_form}
 in the case of the basic IDW and SIDW interpolations, respectively; for the POD variants, this coincides with the CPU 
 time demanded by the online phase,
 i.e., by the resolution of system \eqref{eq:ls3} and by the computation of the linear combination $\mathcal{Z}\boldsymbol{\beta}(\boldsymbol{\mu}^*)$.
Finally, the last column investigates the accuracy of the provided deformation by collecting the value of the $L^2(\Omega_h)$-norm of the relative error 
between the computationally cheaper deformations  and the reference IDW shape.

\begin{figure*}[h]
\centering
\includegraphics[height=0.28\textwidth]{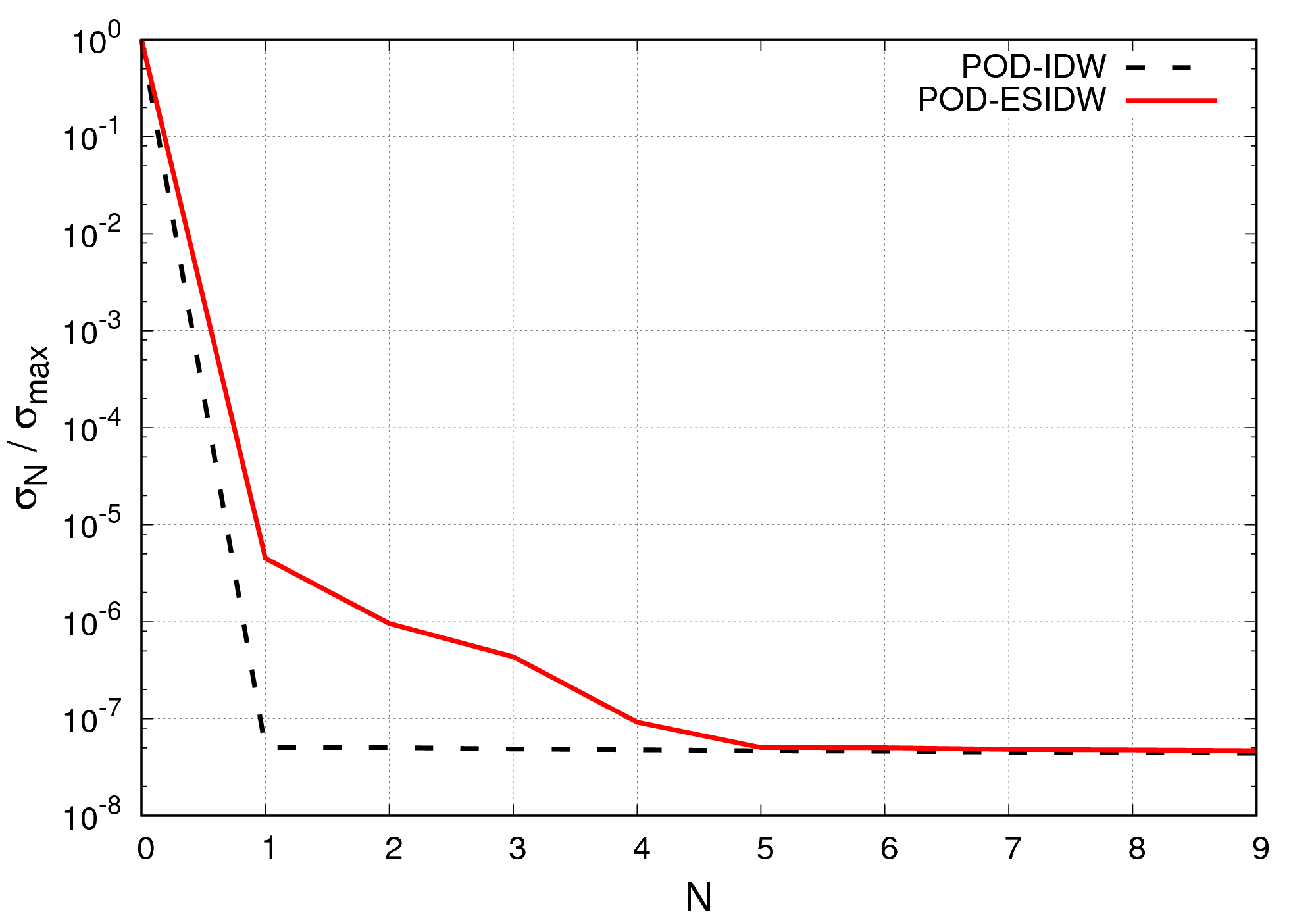}\quad
\includegraphics[height=0.29\textwidth]{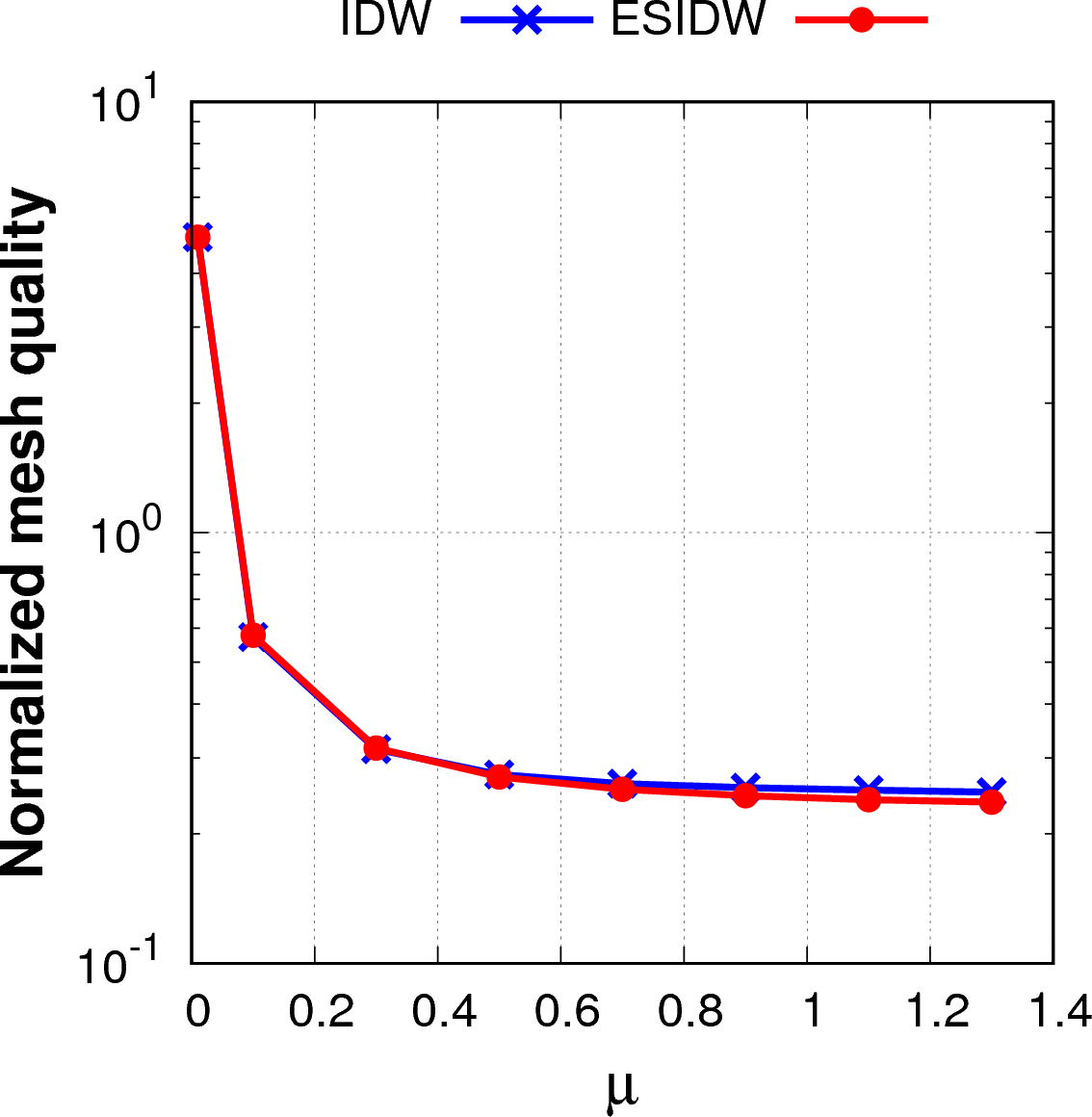}
\includegraphics[height=0.29\textwidth]{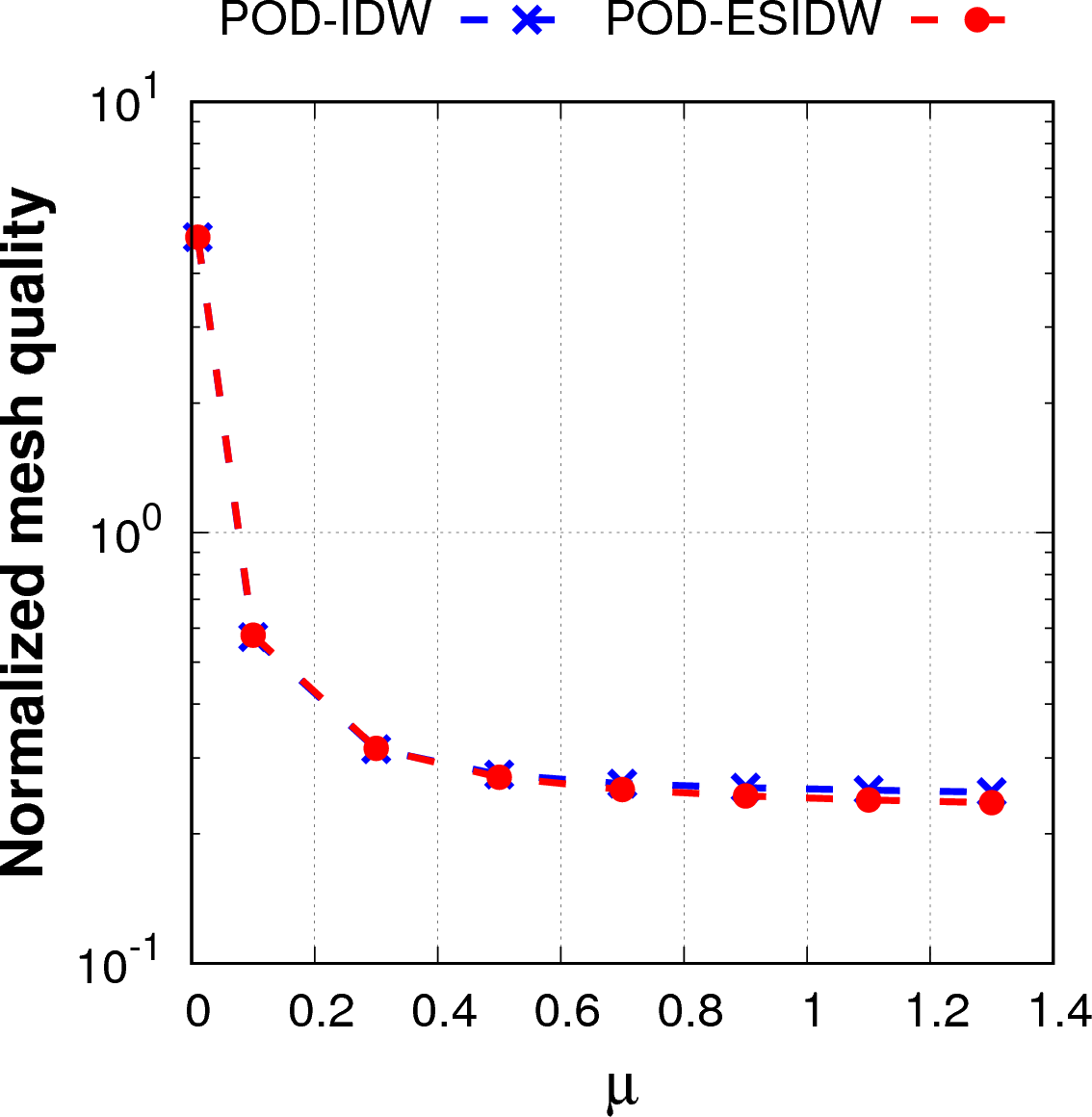}
\caption{Structural deformation of a wing: spectrum decay (of the first ten singular values) for POD-IDW and POD-ESIDW procedures (left);
comparison in terms of mesh quality between IDW and ESIDW interpolation (center) and POD-IDW and POD-ESIDW reduction (right) in the presence of large deformations.}\label{pod1}
\end{figure*}

Concerning the specific values in Table~\ref{summarizing_costs_table_ass1}, we observe that 
both the POD variants identify a minimal reduced basis, a single POD mode being sufficient to ensure the desired tolerance. 
Figure~\ref{pod1}, left shows the decay of the spectrum, normalized to the maximum singular value, for both the POD-IDW and the POD-ESIDW reduction. The first approach exhibits a more evident drop so that a single mode ensures actually an accuracy of $10^{-7}$  (considerably higher than the one demanded). The POD-ESIDW procedure in this case requires a larger number of modes to guarantee the same accuracy, for instance, two modes are demanded for $\varepsilon=10^{-6}$ whereas four modes are required for $\varepsilon=10^{-7}$.
The computational time demanded by the offline phase of both the POD-IDW and POD-ESIDW procedures is not negligible, being equal to $90.93$ [s] and to $76.84$ [s], respectively. Nevertheless, this
phase takes place just once and, as expected, allows a considerable saving in the actual computation of the shape morphing
($0.013$ [s] to be compared with $0.25$ [s] and with $0.034$ [s], respectively), 
thus becoming the ideal tool, for instance, for a multi-query context. 
Moreover, the POD procedure does not compromise the accuracy of the approximation, 
the error (with respect to standard IDW interpolation) being negligible (less than $10^{-8}$) 
in the POD-IDW case, while increasing from 1.06\% to 1.10\% when dealing with the enriched procedure.

Finally, we check how the POD procedure does influence the quality of the deformed mesh, with emphasis on large deformations. In particular, since the mesh quality is expected to deteriorate for larger and larger displacements (of magnitude up to $1.3 \cdot (2 \pi)^2$, i.e. more than $800\%$ of the longitudinal dimension $2\pi$), we consider a normalized mesh quality index, defined as 
the ratio between the mean mesh quality and the maximum displacement measured at the side of the wing which is not clamped. A cross comparison of the plots in Figure~\ref{pod1}, center and right, shows that the POD approach essentially has
no influence on the mesh quality as well as the ESIDW interpolation, the selected index preserving a constant value of about $0.28$ for $\mu>0.3$.

\subsubsection{Fluid mesh motion around a wing}
We move to the parametric version of the FSI test case in Section~\ref{sec:ass2}, dealing with the fluid mesh motion around a NACA0012 profile. The wing is now deformed by the parametrized vertical displacement \eqref{eq:par_ver_disp}, with $\mu \in \mathbb{R}$ varying in $\mathbb{D} = [0, 0.05]$. 

We compare POD-IDW with POD-ESIDW and with the corresponding procedures without any dimensionality reduction.
The parameters $R$, $a$ and $b$ characterizing the selection procedure are the same as in Section~\ref{sec:ass2}, as well as the constraints 
driving the enrichment step of \textbf{Algorithm~\ref{alg:sidw_shape_morphing}}.
The training set $\Xi_{\text{train}}$ is now constituted of $20$ samples randomly distributed in $\mathbb{D}$, 
and the tolerance $\varepsilon$ employed in the POD offline phase is set to $10^{-5}$
as in the previous test case.
The actual deformation is then identified by $\mu^*=0.01$.
\begin{table*}
\centering
\begin{tabular}{|c|c|c|c|c|c|}
\hline
\textbf{method}&{\rm card}($\widehat{\mathcal C}$)&\textbf{CPU time}(\textsc{Offline})&\textbf{{\rm card}(${\mathcal Z}$)}&\textbf{CPU time}(\textsc{Online})&\textbf{Error}\\ 
\hline
IDW & 14126 & - & - & 83.09 [s] & - \\
\hline
ESIDW & 9339 & - & - & $57.07$ [s] & 5.86\% \\
\hline
POD-IDW & 14126 &  2479.51 [s] & 1 & 0.55 [s] & negligible \\
\hline
POD-ESIDW & 9339 & 1613.65 [s] & 1 & 0.55 [s] & 5.94 \% \\
\hline
\end{tabular}
\caption{Fluid mesh motion around a wing: comparison between the basic IDW and ESIDW techniques and the corresponding POD variants.} \label{summarizing_costs_table_ass2}
\end{table*}

\begin{figure*}[h]
\centering
\includegraphics[height=0.3\textwidth]{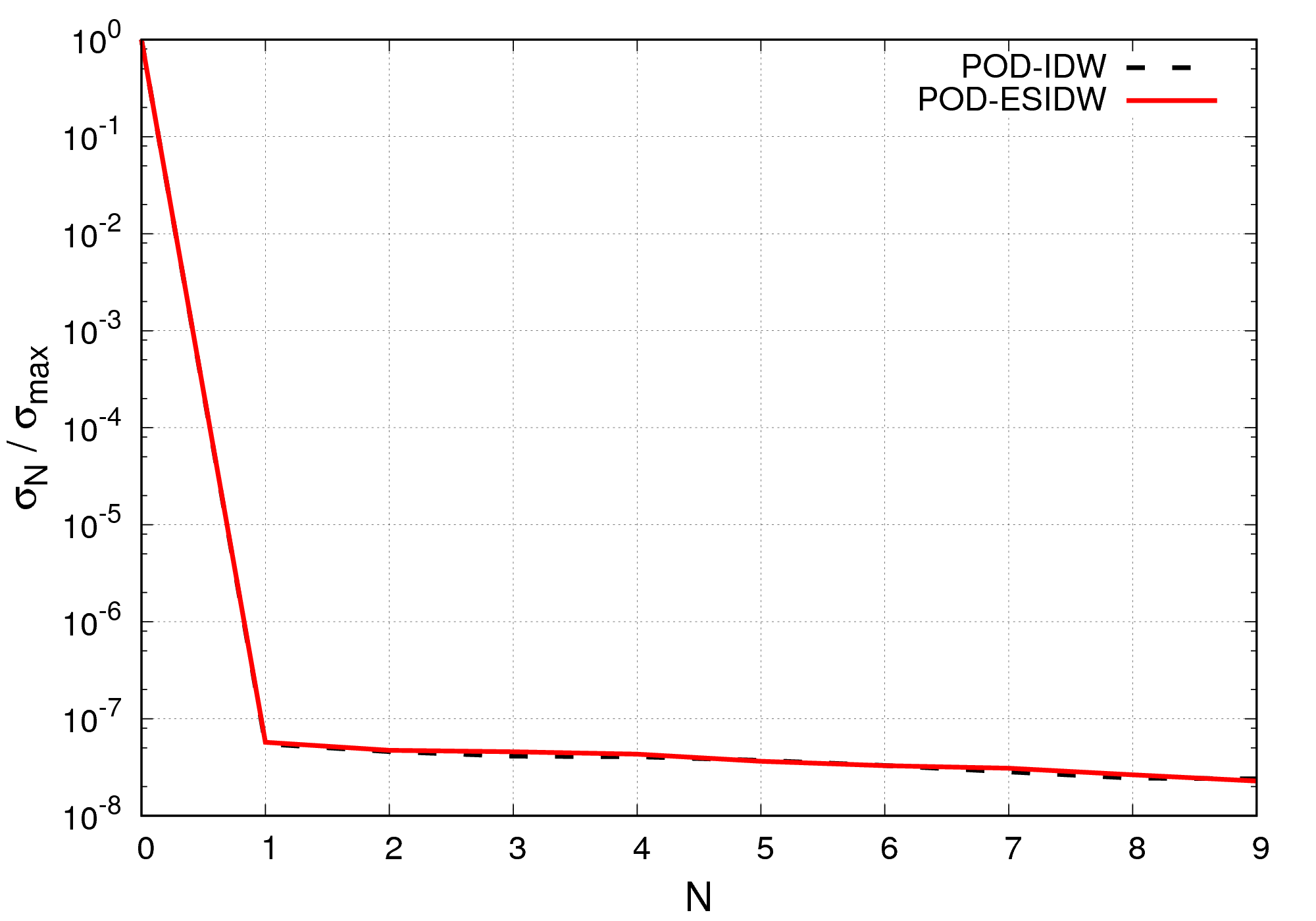}\caption{Fluid mesh motion around a wing: spectrum decay (for the first ten singular values) for POD-IDW and POD-ESIDW procedures.}\label{pod2}
\end{figure*}

Table~\ref{summarizing_costs_table_ass2} offers a summary of the performances of the proposed methods. The successive columns collect the
same quantities as in Table~\ref{summarizing_costs_table_ass1}. The advantages due to the enriched selective procedure have been already highlighted in Section~\ref{sec:ass2}, both in terms of reduction of the control points and of the CPU time.\\
A further saving in the computational time demanded for the actual structure deformation is obtained via POD, 
the online CPU time reducing to $0.55$ [s] for both the POD-IDW and the POD-ESIDW procedures.
The global CPU time reduction due to a selection of the control points combined with a POD procedure
is of two orders of magnitude with respect to the standard IDW approach, for each new morphing. 
This considerable saving is justified by the small dimension of the POD basis, consisting of a unique mode (notice that, for this test case, the trend of the first ten singular values associated with the two POD procedures is exactly the same, as Figure~\ref{pod2}, left shows). The POD offline phase 
remains the most time consuming part of the whole procedure, requiring $2479.51$ [s] and $1613.65$ [s] in the POD-IDW and POD-ESIDW
case, respectively.

Finally, we remark that the accuracy of the POD-ESIDW interpolation is essentially limited by the selection applied to control points. Indeed, the POD-ESIDW approach entails 
a relative error of $5.94\%$ to be compared with an error of $5.86\%$ for the basic ESIDW method.
 
\subsubsection{Fluid mesh motion around a rotating hull}
We conclude the numerical assessment by studying the benefits provided by the POD reduction onto the large deformation setting
in Section~\ref{sec:ass3}. The rotation of the hull with respect to the $z$ axis is now parametrized, $\mu$ coinciding with 
the rotation angle assuming values in the range $\mathbb{D} = \left[-36^\circ,\ 0^\circ\right]$.

\begin{figure*}[h]
\centering
\includegraphics[height=0.28\textwidth]{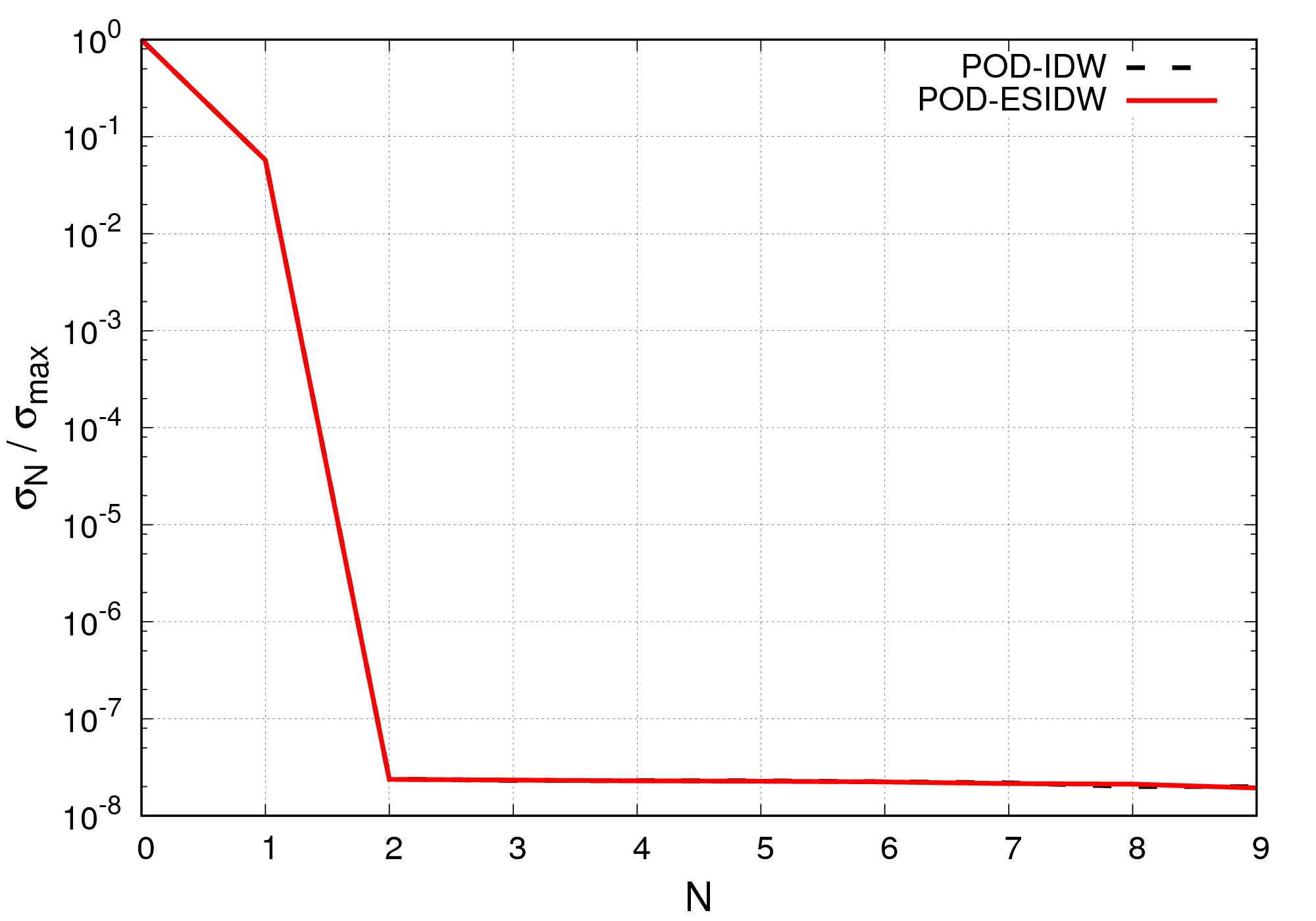}
\includegraphics[height=0.28\textwidth]{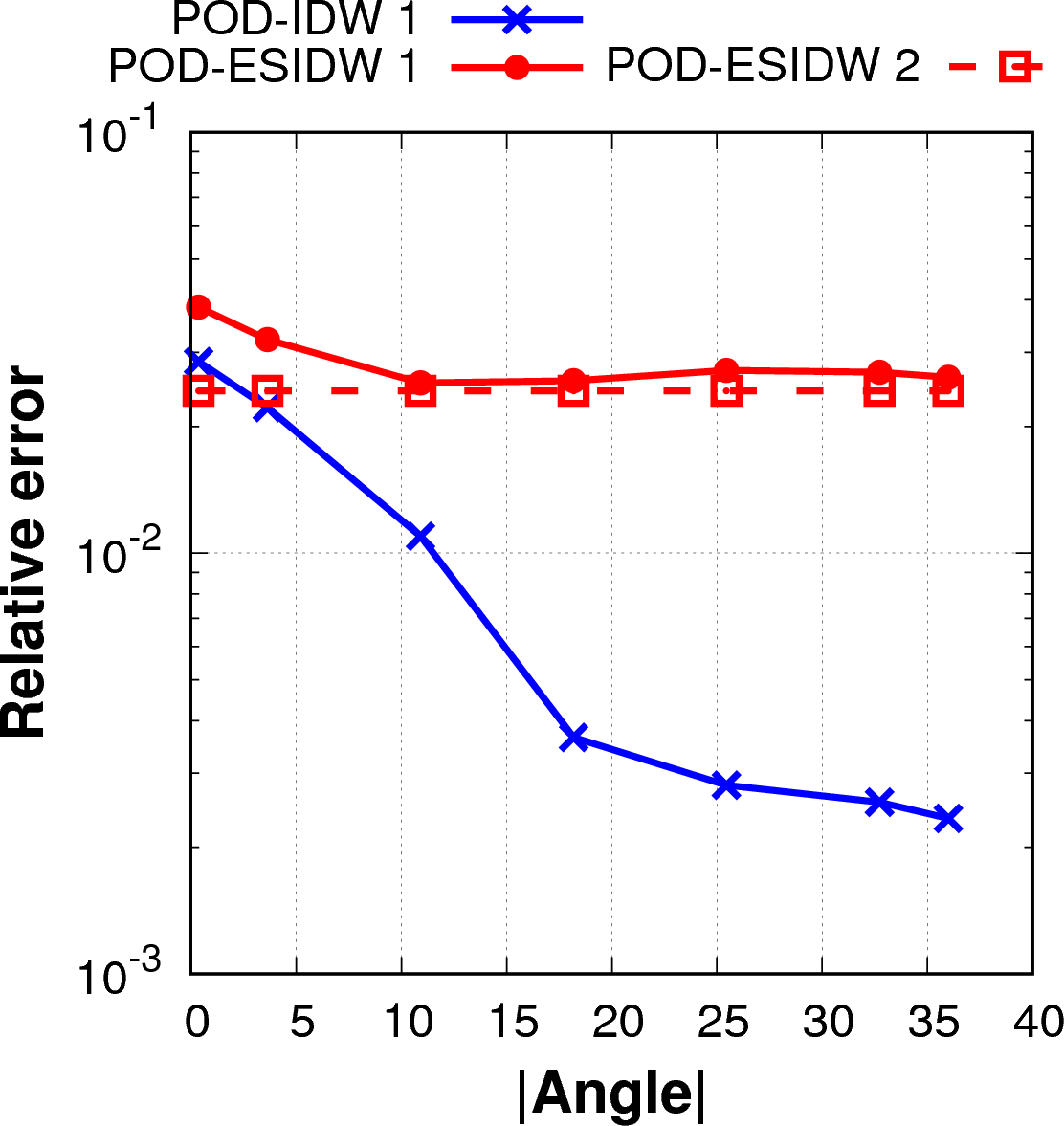}
\includegraphics[height=0.28\textwidth]{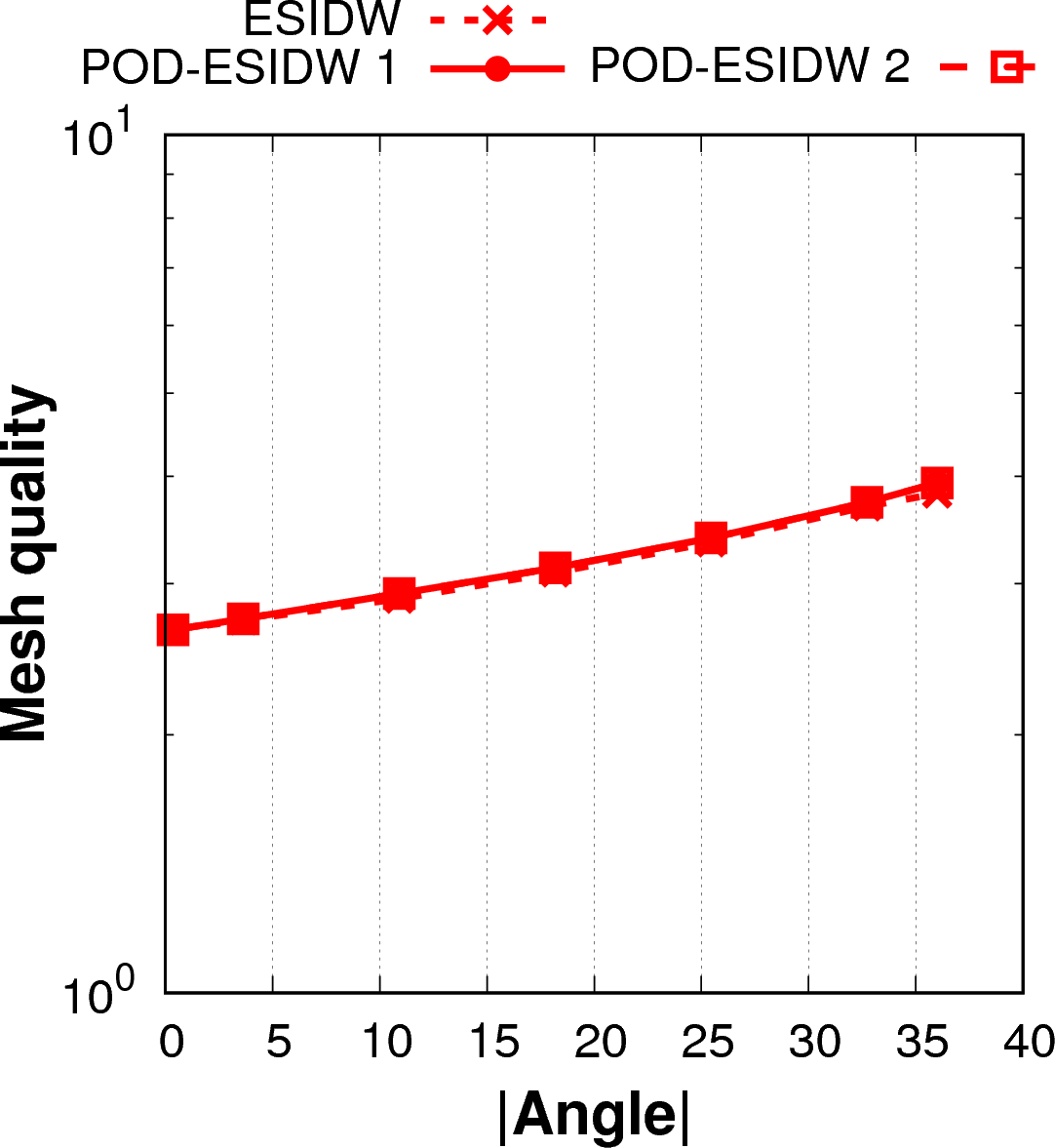}
\caption{Fluid mesh motion around a rotating hull: spectrum decay (for the first ten singular values) for POD-IDW and POD-ESIDW procedures; relative error as a function of the rotation angle for POD-IDW and POD-ESIDW reductions 
(center); comparison in terms of mesh quality between ESIDW and POD-ESIDW procedures (right).}\label{poderror3}
\end{figure*}

As in the two previous sections, POD is combined with IDW and ESIDW interpolations. For the control point selection procedure, 
we adopt the same values for $R$, $a$ and $b$ as in Section~\ref{sec:ass3}. To build the snapshot matrix and to extract the POD basis, 
we exploit $n_{\text{train}} = 50$ values for parameter $\mu$, randomly distributed in $\mathbb{D}$.
The tolerance $\varepsilon$ for the energy check is set to $10^{-5}$. Finally, the target configuration is identified by the parameter $\mu^* = -5^\circ$.
From Figure~\ref{poderror3}, left, we realize that few POD modes suffice to describe the new rotation with a good accuracy.
In particular, tolerance $\varepsilon$ is reached by resorting only to two modes. 
\begin{table*}
\centering
\begin{tabular}{|c|c|c|c|c|c|}
\hline
\textbf{method}&{\rm card}($\widehat{\mathcal C}$)&\textbf{CPU time}(\textsc{Offline})&\textbf{{\rm card}(${\mathcal Z}$)}&\textbf{CPU time}(\textsc{Online})&\textbf{Error}\\ 
\hline
IDW & 3864 & - & - & 4.51 [s] & - \\
\hline
ESIDW & 1057 & - & - & 2.56 [s] & 2.42\% \\
\hline
\multirow{2}{*}{POD-IDW} & \multirow{2}{*}{3864} & \multirow{2}{*}{203.07 [s]} & 1 & 0.38 [s] & 0.23\% $\sim$ 2.85\% \\
\cline{4-6}  &  &  & 2 & 0.38 [s] & negligible \\
\hline
\multirow{2}{*}{POD-ESIDW} & \multirow{2}{*}{1057} & \multirow{2}{*}{150.26 [s]} & 1 & 0.38 [s] & 2.62\% $\sim$ 3.84\% \\
\cline{4-6}
&  & & 2 & 0.38 [s] & 2.43\% \\
\hline
\end{tabular}
\caption{Fluid mesh motion around a rotating hull: comparison between the basic IDW and ESIDW techniques and the corresponding POD variants.} \label{summarizing_costs_table_ass3}
\end{table*}

The performances of POD-IDW and POD-ESIDW methods are summarized in Table~\ref{summarizing_costs_table_ass3}, whose columns are organized as in Tables~\ref{summarizing_costs_table_ass1} and~\ref{summarizing_costs_table_ass2}. 
The computational gain yielded by ESIDW with respect to the standard IDW interpolation is evident, 
following the analysis in Section~\ref{sec:ass3}.\\
Concerning the combined effect of selecting the control points and resorting to a POD reduction, Table~\ref{summarizing_costs_table_ass3}
quantitatively confirms what already remarked for the two previous test cases.
In particular, while ESIDW manages almost to halve the CPU time required for a shape morphing
compared to the standard IDW approach, a successive reduction of about one seventh is reached via 
POD, two basis functions being enough to reach the prescribed tolerance $\varepsilon$.
Similar conclusions hold when comparing IDW with POD-IDW, the CPU time being reduced of about eleven times when resorting to a POD space 
constituted of two basis functions only.
The price to pay for this computationally cheap shape morphing is represented, according to an offline/online paradigm, by the offline phase in
\textbf{Algorithm~\ref{alg:pod_sidw_shape_morphing}}, which demands $203.07$ [s] and $150.26$ [s] for the POD-IDW and the POD-ESIDW approach, respectively.

Figure \ref{poderror3}, center compares the trend of the $L^2(\Omega_h)$-norm of the relative error between the POD-IDW (POD-ESIDW) and the standard IDW deformation
as a function of the rotation angle, when either one or two POD modes are adopted to predict the new deformation. The plot of the POD-IDW procedure associated with two POD modes is omitted, the corresponding relative error being essentially negligible (less than $10^{-8}$).
The accuracy guaranteed by the POD-IDW approach strongly depends on the selected angle when a single POD mode is employed, with a significant improvement in the presence of large rotations. 
On the contrary, a low sensitivity to the rotation angle is exhibited by the POD-ESIDW reduction, the error always being of the order of $10^{-2}$. 
Moreover, from the last column in Table \ref{summarizing_costs_table_ass3}, we deduce that 
the control point selection still represents the principal responsible for an accuracy deterioration
regardless of the selected angle, the relative error remaining essentially the same when combining the ESIDW approach with a POD reduction.

Finally, as shown in Figure \ref{poderror3}, right, the POD reduction does not perturb essentially the quality of the deformed mesh with respect to the standard ESIDW approach, and exhibits a contained dependence on the applied rotation.

\section{Conclusions and perspectives}\label{sec:ROMIDW conclusions}
We have proposed two strategies (and a possible combination of these) to reduce the dimensionality of shape morphing techniques based on IDW.
The first approach (the SIDW/ESIDW interpolation), based on a geometric selection of the control points,
is very general and can be applied to any dimension and to arbitrary meshes. SIDW/ESIDW variants, tested on 
configurations of interest in engineering applications, such as airfoils and hulls, proved to be very effective since they considerably 
reduce the set of control points (i.e., the CPU time) without excessively compromising the accuracy of the approximation.\\
A further reduction of the computational burden is then carried out by means
of a POD (offline phase) and a least squares regression (online phase) techniques.
This allows us to convert any shape morphing into the resolution 
of a linear system of size $N$, with $N$ the number of selected POD basis functions.
In the considered test cases, the number $N$ turns out to be very small. Indeed, one or, at most, two POD basis functions
were enough to guarantee a tolerance of $10^{-5}$ to the energy contained in the discarded POD modes.
We have combined POD dimensionality reduction with IDW and ESIDW interpolations, event though any 
shape morphing algorithm can be alternatively considered. \\
The numerical verification in Section~\ref{combined_Sec_num} shows that
the combined POD-IDW and POD-ESIDW techniques lead to a large computational saving, 
up to two orders of magnitude on the CPU time. The offline phase remains the most time-consuming part of the procedure, 
according to an offline/online paradigm.
Finally, an accuracy analysis highlights that SIDW/ESIDW entail very small relative errors (few percentage points) 
with respect to the standard IDW procedure, while POD
reduction does not significantly contribute further. 

A possible future development of this work might concern the integration of the proposed methods into
a FSI solver, or the application to several optimization contexts. An adaptive selection of the control points driven by some quantity of interest, combination with reduction procedures for parametrized problems 
(e.g., \cite{HiPOD,PGD,hesthavenrozzastamm2016}), as well as the use of active subspaces method \cite{constantine2015active,Marco} as pre-processing,
represent further topics of interest for the following of the current work.

\section*{Acknowledgments}
We gratefully acknowledge Dr. Andrea Mola (SISSA mathLab) for helpful discussions and his support with mesh generation. We acknowledge the support by European Union Funding for Research and Innovation -- Horizon 2020 Program -- in the framework of European Research Council Executive Agency (H2020 ERC CoG 2015 AROMA-CFD project 681447 ``Advanced Reduced Order Methods with Applications in Computational Fluid Dynamics''). We also acknowledge the INDAM-GNCS project ``Metodi Numerici Avanzati Combinati con Tecniche di Riduzione Computazionale per PDEs Parametrizzate e Applicazioni''.
Alessandro D'Amario has been supported by a SISSA scholarship (Mathematics Area).

\bibliography{bib/bibliography}
\bibliographystyle{abbrv}

\end{document}